\newif\ifSMF
\newif\ifBOOK
\newif\ifABSTRACT
\makeatletter
\@ifundefined{SMFstyle}{\SMFfalse}{\SMFtrue}
\@ifundefined{BOOKlayout}{\BOOKfalse}{\BOOKtrue}
\@ifundefined{ABSTRACTonly}{\ABSTRACTfalse}{\ABSTRACTtrue}
\makeatother

\ifSMF
  \documentclass[a4paper,reqno,twoside]{smfart}
\else
  \ifABSTRACT
    \documentclass[a4paper,reqno]{article}
  \else
    \documentclass[a4paper,reqno,twoside]{article}
  \fi
\fi

\ifSMF
  \usepackage[francais,english]{babel}
\else
  \usepackage[english,francais]{babel}
\fi

\usepackage[latin1]{inputenc}
\usepackage[T1]{fontenc}
\usepackage{textcomp}

\ifSMF
  \usepackage{amsmath,amssymb,amsfonts,smfthm,bull,smfenum}
\else
  \usepackage{amsmath,amssymb,amsfonts,amsthm}
\fi

\usepackage{graphicx}

\usepackage{multicol}

\ifSMF
  \author[V. Feuvrier]{Vincent Feuvrier}
  \address{Département de Mathématiques\\Bâtiment 425\\Faculté des Sciences d'Orsay\\Université Paris-Sud 11\\F-91405 Orsay Cedex\\France}
  \email{vincent.feuvrier@normalesup.org}
  \urladdr{http://www.math.u-psud.fr/~feuvrier}
  \title[Remplissage polyédrique de rotondité uniforme]{Remplissage de l'espace Euclidien par des complexes polyédriques d'orientation imposée et de rotondité uniforme}
  \alttitle{Filling Euclidean space using polyhedric complexes of imposed orientation with uniform rotondity}
\else
  \author{Vincent Feuvrier}
  \title{Remplissage de l'espace Euclidien par des complexes polyédriques d'orientation imposée et de rotondité uniforme}
  \def\frontmatter{}
  \def\mainmatter{}
  \def\backmatter{}
  \newenvironment{altabstract}{%
    \selectlanguage{english}%
    \begin{abstract}%
  }{%
    \end{abstract}%
    \selectlanguage{francais}%
  }
  \def\subjclass#1{}%
  \def\keywords#1{}%
  \def\altkeywords#1{}%
\fi

\ifBOOK
\fi

\newtheorem{definition}{Définition}
\newtheorem{theorem}{Théorème}
\newtheorem{lemma}{Lemme}
\newtheorem{property}{Propriété}

\def\open#1{\overset{\circ}{#1}}
\def\dist{\mathbf{d}}
\DeclareMathOperator\vect{Vect}
\DeclareMathOperator\affine{Affine}
\DeclareMathOperator\extrem{Extrem}
\def\fracd#1#2{\frac{\displaystyle #1}{\displaystyle #2}}
\def\A{\mathcal{A}}
\def\C{\mathcal{C}}
\def\F{\mathcal{F}}
\def\G{\mathcal{G}}
\def\H{\mathcal{H}}
\def\K{\mathcal{K}}
\def\M{\mathcal{M}}
\def\O{\mathcal{O}}
\def\P{\mathcal{P}}
\def\R{\mathcal{R}}
\def\S{\mathcal{S}}
\def\T{\mathcal{T}}
\def\U{\mathcal{U}}
\def\V{\mathcal{V}}

\def\itemx#1#2{%
  \item[$\mathbf{\open{F_#1}\cap\open{G_#2}}$]%
}

\usepackage{verbatim,listings}
\usepackage[noend]{algorithm,algorithmic}

\makeatletter
\@ifundefined{algorithmicprint}{
  \def\TRUE{\textbf{vrai}}
  \def\FALSE{\textbf{faux}}
  \def\RETURN{\item\textbf{renvoyer}}
}{}
\makeatother

\floatname{algorithm}{Algorithme}

\makeatletter
\@ifundefined{algsetup}{}{
  \algsetup{
    linenosize=\footnotesize,
    linenodelimiter=.,
    indent=20pt
  }
}
\makeatother

\begin{document}

\frontmatter

\ifSMF
\else
  \maketitle
\fi

\begin{abstract}
Nous donnons une méthode de construction de complexes polyédriques dans $\mathbb{R}^n$ permettant de relier entre elles des grilles dyadiques d'orientations différentes tout en s'assurant que les polyèdres utilisés ne soient pas trop plats, y compris leurs sous-faces de toutes dimensions. Pour cela, après avoir rappelé quelques définitions et propriétés simples des polyèdres euclidiens compacts et des complexes, on se dote d'un outil qui permet de remplir de polyèdres $n$-dimensionnels un ouvert en forme de tube dont la frontière est portée par un complexe $n-1$-dimensionnel. Le théorème principal est démontré par induction sur $n$ en reliant les complexes dyadiques couche par couche, en remplissant des tubes disposés autour des différentes couches et en utilisant le théorème en dimension inférieure pour construire les morceaux manquants de la frontière des tubes. Une application possible de ce résultat est la recherche de solutions à des problèmes de minimisation de la mesure en dimension et codimension quelconques dans certaines classes topologiques.
\end{abstract}

\begin{altabstract}
We build polyhedral complexes in $\mathbb{R}^n$ that coincide with dyadic grids with different orientations, while keeping uniform lower bounds (depending only on $n$) on the flatness of the added polyhedrons including their subfaces in all dimensions. After the definitions and first properties of compact Euclidean polyhedrons and complexes, we introduce a tool allowing us to fill with $n$-dimensionnal polyhedrons a tubular-shaped open set, the boundary of which is a given $n-1$-dimensionnal complex. The main result is proven inductively over $n$ by completing our dyadic grids layer after layer, filling the tube surrounding each layer and using the result in the previous dimension to build the missing parts of the tube boundary. A possible application of this result is a way to find solutions to problems of measure minimization over certain topological classes of sets, in arbitrary dimension and codimension.
\end{altabstract}

\subjclass{}
\keywords{Polyèdres euclidiens, Complexes polyédriques, Pavages polyédriques}
\altkeywords{Euclidean polyhedrons, Polyhedric complexes, Polyhedric tesselations}

\ifSMF
  \maketitle
  \ifABSTRACT
    \end{document}
  \fi
  \tableofcontents
\else
  \ifABSTRACT
    \end{document}
  \fi
  \def\contentsname{Sommaire}
  \pagebreak
  \null
  \vfill
  \tableofcontents
  \vfill
  \vfill
  \pagebreak
\fi

\mainmatter

\section{Introduction}\label{sectionA}

Le résultat principal de ce papier (le théorème~\ref{theoremfusion}, dit <<~de fusion~>>) peut s'énoncer simplement:
\begin{quotation}
<<~\'Etant donnés deux complexes dyadiques $S_1$ et $S_2$ tels qu'un morceau de la frontière de $S_1$ forme la frontière d'un ouvert borné $O$ disjoint de $S_1$ qui contient $S_2$, si la distance séparant $S_1$ et $S_2$ est suffisamment grande devant la taille des cubes dyadiques considérés alors on peut construire un complexe $S_3$ tel que $S_2\cup S_3$ remplisse $\overline{O}$, avec une borne inférieure uniforme sur la rotondité des polyèdres construits et leurs sous-faces.~>>
\end{quotation}

Le lecteur qui se risquerait à effectuer un dessin en dimension $2$ n'aurait probablement aucun mal à compléter deux complexes dyadiques d'orientations différentes en s'imposant une borne inférieure raisonnable sur les angles des segments qui s'intersectent. Il obtiendrait d'ailleurs vraisemblablement des bornes inférieures sur la rotondité et la distance entre les deux complexes meilleures que celles du théorème~\ref{theoremfusion}. Bien évidemment les choses se compliquent en dimension plus grande, en particulier faute d'outils descriptifs <<~morphologiques~>> efficaces. En outre des problèmes supplémentaires surviennent lorsque $n\geq4$, rendant les représentations beaucoup plus difficiles.

\begin{figure}
\begin{center}\includegraphics[width=\textwidth]{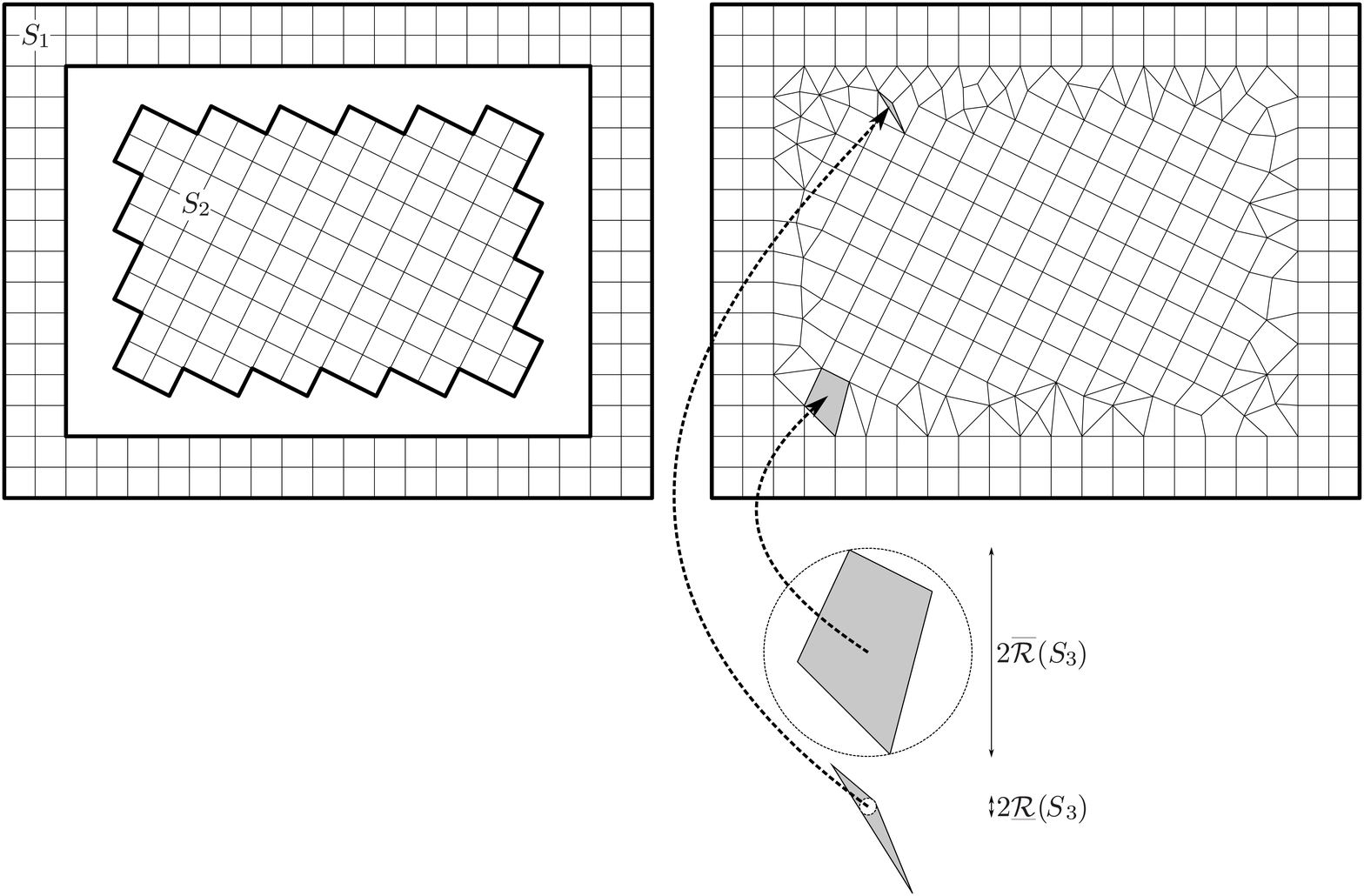}\end{center}
\caption{Fusion à la main de deux complexes bidimensionnels et les constantes de forme obtenues}
\label{figureintroductionA}
\end{figure}

Nous donnons d'abord une définition simple et intuitive (définition~\ref{definitionpolyedre}) des polyèdres euclidiens convexes de $\mathbb{R}^n$ en tant qu'intersection compacte de demi-espaces affines, équivalente à celle des polytopes (propriété~\ref{propertypolytope}) et amenant naturellement la définition des faces et des sous-faces (définition~\ref{definitionsousfaces}). On pourra lire à ce sujet l'article d'Andrée Bastiani \cite{bastiani} pour des définitions plus générales dans des espaces topologiques. Nous introduisons une quantité (la rotondité, définition~\ref{definitionrotondite}) permettant de contrôler la forme d'un polyèdre donné en considérant le rapport compris entre $0$ et $1$ des rayons d'une boule inscrite et d'une boule circonscrite; plus ce rapport est proche de zéro, plus le polyèdre est aplati. Enfin pour formaliser l'idée intuitive de familles de polyèdres de même dimension qui se raccordent bien entre eux nous définissons la notion de complexe (définition~\ref{definitioncomplexe}), en imposant que les polyèdres et leurs sous-faces de dimensions inférieures soient d'intérieurs disjoints deux à deux (hormis ceux qui sont confondus); c'est le cas par exemple des complexes dyadiques, abordés dans la section~\ref{sectionC}. La notion de rotondité est généralisée aux complexes en considérant la rotondité minimale parmi les polyèdres et les sous-faces, de façon à pouvoir minorer la rotondité de la sous-face la plus aplatie.

Avant d'énoncer le plan de l'article donnons rapidement et sans démonstration une application possible de ce résultat inspirée de Reifenberg \cite{reifenberg} pour trouver des ensembles de mesure minimale parmi certaines classes topologiques. Par exemple, trouver parmi une classe $\mathfrak{F}$ stable par des déformations lipschitziennes un ensemble $E\in\mathfrak{F}$ tel que
\begin{equation}
\H^d(E)=\inf_{F\in\mathfrak{F}}\H^d(F),
\end{equation}
la mesure utilisée ici étant la mesure de Hausdorff $d$-dimensionnelle $\H^d$ avec $0\leq d<n$ (on pourra trouver plus de détails dans le livre de Mattila \cite{mattila}). Notons que la technique utilisée reste valable pour la minimisation de fonctionnelles ensemblistes plus générales.

Considérons un polyèdre $n$-dimensionnel $\delta$ (convexe par définition, de rotondité $R(\delta)$) et $c\in\open{\delta}$. Notons $\Pi_{\delta,c}$ la projection radiale sur $\partial\delta$ qui à $x\in\delta\setminus\{c\}$ associe l'unique intersection de la demi-droite $[c,x)$ avec $\partial\delta$. Posons $d=n-1$ et soit $E\subset\delta$ une sous-partie fermée $\H^d$-mesurable telle que $\H^d(E)<\infty$. En calculant la valeur moyenne de $\H^d(\Pi_{\delta,c}(E))$ lorsque $c$ parcourt $\open{\delta}\setminus E$ on peut montrer en utilisant Fubini qu'il existe une constante $K>0$ ne dépendant que de $d$ et $n$ telle que
\begin{equation}
\exists c\in\open{\delta}\colon\H^d(\Pi_{\delta,c}(E))\leq KR(\delta)^{-2d}\H^d(E).
\end{equation}
En utilisant par exemple le théorème d'extension lipschitzienne de Kirzbraun \cite{kirszbraun} on peut montrer que pour tout complexe $S$ de rotondité $\R(S)$ et toute sous-partie fermée $E$ telle que $E\subset\U(S)$ on peut trouver une application lipschitzienne $\phi$ telle que $\phi(E)$ est inclus dans les faces de $S$ et
\begin{equation}\label{equationintroductionA}
\H^d(\phi(E))\leq K\R(S)^{-2d}\H^d(E).
\end{equation}
En continuant les projections radiales dans les sous-faces de dimension inférieure qui ne sont pas entièrement recouvertes on peut même effectuer cette construction en codimension $n-d\geq1$ quelconque, et imposer que $\phi(E)$ soit une union finie de sous-faces de dimension au plus $d$ de $S$.

Supposons que $E\in\mathfrak{F}$. Lorsque $E$ est rectifiable, par un lemme de type Vitali on peut recouvrir $E$ à une partie de mesure arbitrairement petite près par une union finie de complexes dyadiques disjoints dont les orientations suivent la direction des plans tangents approximatifs de $E$. D'après le théorème~\ref{theoremfusion} il est alors possible de relier tous ces complexes dyadiques en un complexe plus grand $S$ de façon à avoir à la fois $E\subset\U(S)$ et $\R(S)>C$ où $C$ ne dépend que de $n$. En projetant préalablement $E$ sur ses plans tangents approximatifs et en composant avec les projections radiales mentionnées plus haut on peut construire une application lipschitzienne $\psi$ telle que cette fois $E'=\psi(E)$ soit une union finie de sous-faces de dimension au plus $d$ de $S$ et
\begin{equation}
\H^d(E')\leq (1+\epsilon)\H^d(E).
\end{equation}
En minimisant parmi les éléments de $\mathfrak{F}$ qui sont des unions de sous-faces de dimension au plus $d$ de $S$ (il y en a un nombre fini) on peut trouver un ensemble polyédrique optimal $E''$ qui vérifie en particulier
\begin{equation}
\H^d(E'')\leq\H^d(E')\leq (1+\epsilon)\H^d(E).
\end{equation}
Par ailleurs pour toute déformation lipschitzienne $F$ de $E''$ à l'intérieur de $\U(S')$, d'après~\eqref{equationintroductionA} et puisque $E''$ est optimal on a
\begin{equation}\label{equationintroductionB}
\H^d(E'')\leq KC^{-2d}\H^d(F)
\end{equation}
c'est à dire que $E''$ est $M$-quasiminimal avec $M=KC^{-2d}$ qui ne dépend que de $d$ et $n$.

Pour résumer en considérant une suite minimisante $E_k$ de $\mathfrak{F}$, c'est à dire telle que 
\begin{equation}
\lim_{k\rightarrow +\infty}\H^d(E_k)=\inf_{F\in\mathfrak{F}}\H^d(F)
\end{equation}
il est donc possible de construire automatiquement une suite minimisante d'ensembles quasiminimaux polyédriques de $\mathfrak{F}$. Dans ce cas un résultat de Guy David dans \cite{david} établit la semi-continuité inférieure de la mesure par passage à la limite, ce qui n'est généralement pas le cas et constitue une difficulté technique pour la recherche de minimiseurs.

Une version détaillée de ce processus d'optimisation polyédrale devrait faire l'objet d'un prochain article. Notre méthode pourrait par ailleurs permettre de généraliser en dimension et codimension quelconques un résultat de Thierry De Pauw dans \cite{depauw:acr} basé sur un théorème d'approximation polyédrale d'ensembles rectifiables de dimension $2$ dans $\mathbb{R}^3$, dans un cadre de recherche de minimiseurs de taille pour les courants entiers.

Le plan de l'article est le suivant.

La section~\ref{sectionB} est consacrée aux définitions et propriétés immédiates des polyèdres et complexes euclidiens. En particulier nous démontrons le lemme~\ref{lemmasuspensiontubulaire} de suspension tubulaire, un outil qui permet de remplir de polyèdres un ouvert en forme de tube dont la frontière est portée par un complexe $n-1$-dimensionnel en construisant un complexe $n$-dimensionnel qui remplit son adhérence en s'appuyant sur le complexe-frontière. Ce résultat suppose de prendre les précautions nécessaires pour s'assurer que les polyèdres construits sont d'intérieurs disjoints, comme indiqué dans le lemme~\ref{lemmaocclusionsimple}, et permet de contrôler la rotondité du complexe obtenu en fonction de la forme du tube.

La section~\ref{sectionC} est consacrée aux complexes dyadiques et à la démonstration du théorème~\ref{theoremfusion} par récurrence sur la dimension $n$. On commence par le démontrer en dimension $2$ (lemme~\ref{lemmafusion2}). Puis l'induction est prouvée par le lemme~\ref{lemmafusionn}. En supposant que les bases des deux complexes dyadiques à faire fusionner sont l'image l'une de l'autre par une rotation planaire parallèle aux cubes il est possible de travailler couche par couche. Il est en outre possible d'imposer que l'angle de rotation soit arbitrairement proche de zéro. Nous obtenons ces deux conditions en décomposant une isométrie affine de changement de base entre les deux complexes dyadiques en un produit de rotations planaires, puis en décomposant chaque rotation en un produit de rotations d'angle suffisamment proche de zéro. Il suffit alors de remplir des couches de transitions en <<~oignon~>> pour passer de $S_1$ à $S_2$ et supposant $\rho$ suffisamment grand, le nombre total de couches ne dépendant que de $n$ et du choix de l'angle de rotation maximal.

On effectue alors des suspensions tubulaires autour des différentes couches en utilisant le théorème en dimension inférieure pour compléter les parties manquantes de la frontière des couches. L'un des problèmes techniques est qu'on ne dispose pas de bornes uniformes sur le diamètre des tubes utilisés pour les suspensions. Il font donc creuser des canalisations (définition~\ref{definitioncanalisation}) à la surface des deux complexes à fusionner, qui peuvent s'imbriquer de manière complémentaire (lemme~\ref{lemmacanalisationsn}). Ces canalisations sont obtenues en étudiant des complexes dyadiques bidimensionnels. Nous prouvons le lemme~\ref{lemmalaboureur} dit <<~du laboureur~>> qui permet de creuser des sillons dans un complexe dyadique bidimensionnel de façon à ce que le talus de ces sillons (le complémentaire des cubes enlevés) forme lui aussi quasiment des sillons, avant de généraliser en dimension plus grande.

La dernière section est une étude de différents cas intervenant dans la démonstration du lemme du laboureur. \'Etant donné le nombre élevé de cas à considérer, nous utilisons un algorithme et son implémentation en langage C pour terminer la démonstration.

Je tiens à remercier Guy David pour son constant soutien, ses nombreux conseils et suggestions.

\section{Polyèdres et complexes euclidiens}\label{sectionB}

Un demi-espace affine $A$ est la somme directe d'un hyperplan affine $H$ avec une demi-droite $\mathbb{R}^+u$ où $u$ est une direction non parallèle à $\overrightarrow{H}$, c'est à dire
\begin{equation}
A=\{x+ru\colon x\in H\text{ et }r\geq0\}.\end{equation}
On dira qu'une intersection de demi-espaces affines est un polyèdre, au sens de la définition suivante.

\begin{definition}[Polyèdres]\label{definitionpolyedre}
Un polyèdre $\delta$ de dimension $n$ est une partie compacte de $\mathbb{R}^n$ d'intérieur non vide, obtenue par intersection finie de demi-espaces affines.
\end{definition}

En ne gardant que les demi-espaces affines dont la frontière intersecte $\delta$ sur une sous-partie de dimension de Hausdorff égale à $n-1$ on vérifie facilement que parmi toutes les familles de demi-espaces affines qui peuvent convenir il en existe une minimale pour l'inclusion; on la notera $\A(\delta)$.

En autorisant des parties compactes non vides mais d'intérieur vide on généralise aussi la définition à des polyèdres de dimension $k\leq n$ en considérant la dimension $k$ du plus petit sous-espace affine qui les contient, noté $\affine(\delta)$. Dans ce cas les différents opérateurs topologiques usuels (frontière, adhérence ou intérieur) seront pris relativement à ce sous-espace affine minimal, de même que les demi-sous-espaces affines dans $\A(\delta)$. Par convention, on considère que les singletons sont des polyèdres de dimension $0$, égaux à leur intérieur et de frontière vide.

Avec ces conventions, un argument simple de convexité permet d'établir la correspondance entre la dimension de $\affine(\delta)$ et la dimension de Hausdorff de $\delta$, qu'on note $\dim\delta$.

\subsection{Sous-faces et suspensions de polyèdres}

Les polyèdres tels qu'on les a définis sont convexes et possèdent des faces et des sous-faces.

\begin{definition}[Sous-faces]\label{definitionsousfaces}
Soit $\delta$ un polyèdre $n$-dimensionnel tel que $\A(\delta)=\{A_1,\ldots,A_p\}$ et $\{A'_1,\ldots,A'_p\}$ une famille de sous-parties de $\mathbb{R}^n$ telle que $A'_i=A_i$ ou $A'_i=\partial A_i$ pour $1\leq i\leq p$. En posant $\alpha=\bigcap_i A'_i$, si $\alpha\neq\emptyset$ on dira que $\alpha$ est une sous-face de $\delta$, et plus précisément:
\begin{itemize}
\item si $\dim\alpha<\dim\delta$ on dira que $\alpha$ est une sous-face stricte;
\item si $\dim\alpha=\dim\delta-1$ on dira que $\alpha$ est une face;
\item si $\dim\alpha=0$ (autrement dit si $\alpha$ est un singleton) on dira que $\alpha$ est un sommet, et on le confondra souvent avec le point qu'il contient.
\end{itemize}
On note $\F(\delta)$ l'ensemble des sous-faces de $\delta$ (dont $\delta$ lui-même) et pour $0\leq k\leq\dim\delta$ l'ensemble des sous-faces $k$-dimensionnelles
\begin{equation}
\F_k(\delta)=\{\alpha\in\F(\delta)\colon\dim\alpha=k\}.
\end{equation}
\end{definition}

Là encore on généralise naturellement cette définition à des polyèdres de dimension $k\leq n$. On peut vérifier aisément que les sous-faces sont elles-mêmes des polyèdres, et que les sous-faces des sous-faces de $\delta$ sont aussi des sous-faces de $\delta$. En outre les faces sont d'intérieur disjoint, et leur union forme la frontière du polyèdre. On peut même écrire que
\begin{equation}
\delta=\bigsqcup_{\alpha\in\F(\delta)}\open{\alpha}
\end{equation}
où $\sqcup$ désigne une union disjointe, l'intérieur des sous-faces étant pris à chaque fois relativement au sous-espace affine engendré correspondant.

\begin{figure}
\begin{center}\includegraphics[width=0.4\textwidth]{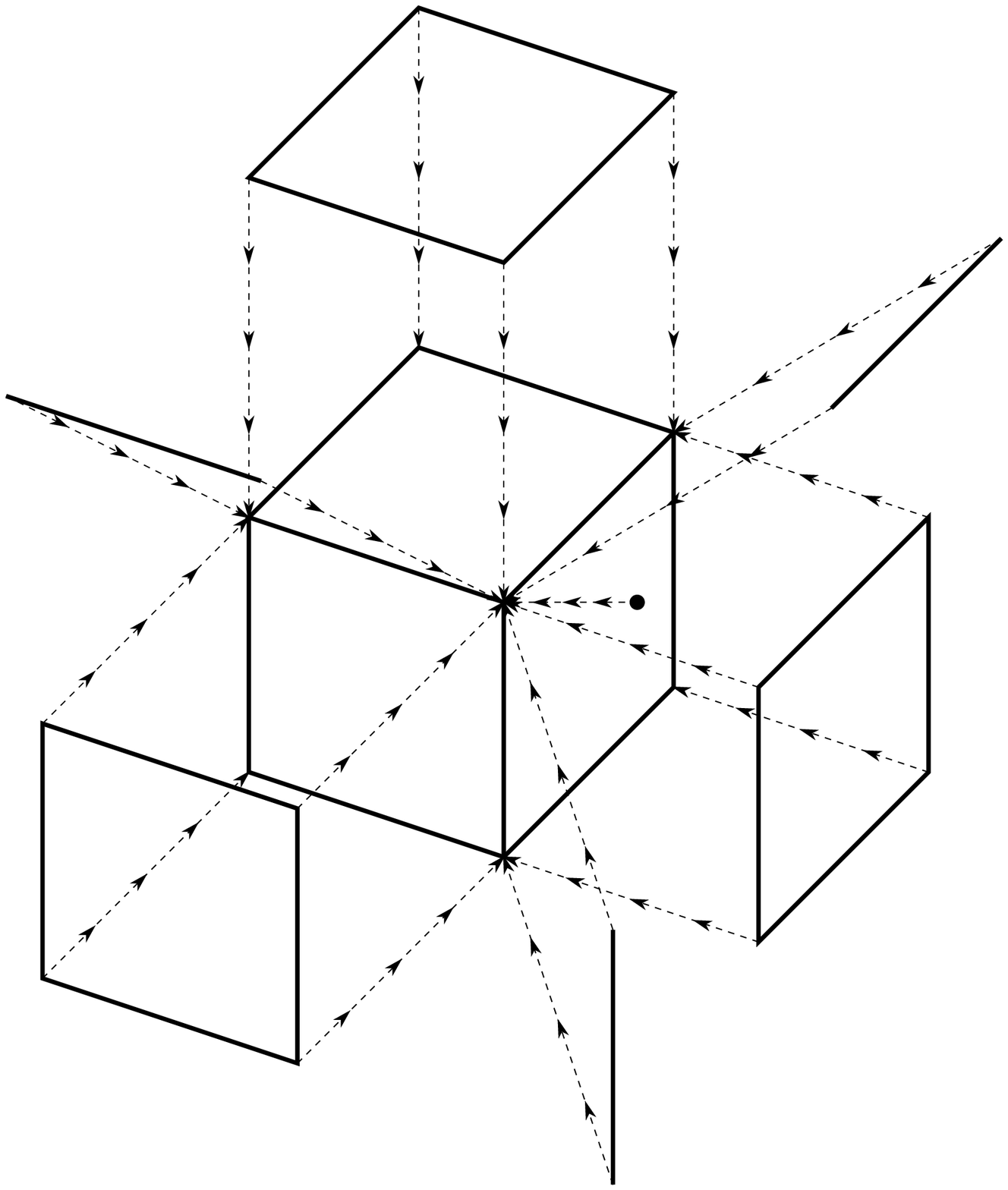}\end{center}
\caption{Un polyèdre (ici un cube) et ses sous-faces visibles en dimension $3$}
\label{figurepolyedreA}
\end{figure}

De façon à rendre plus lisibles certains des énoncés à venir on va encore se doter de la définition suivante.

\begin{definition}[Suspension]\label{definitionsuspension}
Pour une partie $A\subset\mathbb{R}^n$ et $x\in\mathbb{R}^n$ on définit la suspension de $A$ par rapport à $x$ par
\begin{equation}
\S(A,x)=\{ty+(1-t)x\colon y\in A\text{ et }t\in[0,1]\}.
\end{equation}
\end{definition}

Pour une partie $A$ on note $\left<A\right>$ son enveloppe convexe, c'est à dire l'intersection de tous les convexes qui la contiennent. Il est facile de vérifier que si $A$ est convexe, pour deux points quelconques $x$ et $y$ on a $\S(A,x)=\left<A\cup\{x\}\right>$ et $\S(\S(A,x),y)=\S(\S(A,y),x)$

\begin{figure}
\begin{center}\includegraphics[width=0.7\textwidth]{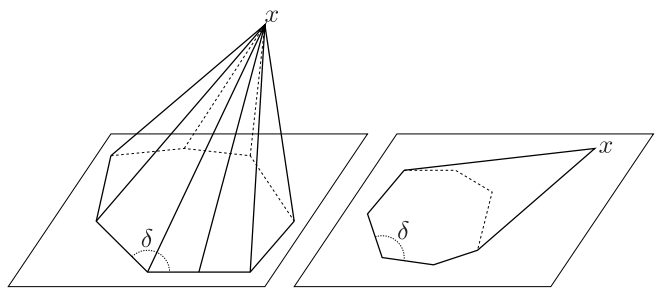}\end{center}
\caption{Suspension d'un polyèdre $\delta$ par rapport à un point $x$ coplanaire ou non}
\label{figuresuspensionA}
\end{figure}

Pour une partie $A\subset\mathbb{R}$ on note l'ensemble de ses points extrémaux
\begin{equation}
\extrem(A)=\left\{x\in A\colon x\notin\left<A\setminus\{x\}\right>\right\}.
\end{equation}
Avec ces notations on peut donner la caractérisation suivante, qui établit l'équivalence entre les polyèdres de la définition~\ref{definitionpolyedre} et les polytopes. Il est possible de la démontrer (ce qu'on ne fera pas ici) par récurrence sur le nombre de points de $\extrem(\delta)$ pour l'implication $1\Rightarrow2$, par récurrence sur la dimension de $\delta$ pour l'implication réciproque, et en utilisant le théorème de Krein-Milman dans les espaces euclidiens \cite{kreinmilman} pour le dernier point.

\begin{property}[Polytopes]\label{propertypolytope}
Pour toute partie convexe compacte non vide $\delta\subset\mathbb{R}^n$, les deux énoncés suivants sont équivalents:
\begin{enumerate}
\item $\extrem(\delta)$ est finie;
\item $\delta$ est un polyèdre.
\end{enumerate}
En particulier, si $\delta$ est un polyèdre alors $\extrem(\delta)=\F_0(\delta)$.
\end{property}

\'Etant donné un compact $A\subset\mathbb{R}^n$, donnons-nous trois quantités permettant de contrôler sa forme, et en particulier dans le cas d'un polyèdre de donner une borne inférieure implicite sur les angles que font ses faces entre elles.

\begin{definition}[Rotondité]\label{definitionrotondite}
Les constantes de forme d'un compact $A\subset\mathbb{R}^n$ sont:
\begin{itemize}
\item le supremum des rayons des boules incluses relativement au sous-espace affine engendré (avec la convention $\sup\emptyset=0$) appelé \emph{rayon intérieur}
\begin{equation}
\underline{R}(A)=\sup\{r>0\colon\exists x\in\mathbb{R}^n,A\supset B(x,r)\cap\affine(A)\};
\end{equation}
\item l'infimum des rayons des boules qui le contiennent (avec la convention $\inf\emptyset=0$) appelé \emph{rayon extérieur}
\begin{equation}
\overline{R}(A)=\inf\{r>0\colon\exists x\in\mathbb{R}^n,A\subset B(x,r)\};
\end{equation}
\item le rapport des deux (avec la convention $R(A)=1$ lorsque $\overline{R}(A)=0$) appelé \emph{rotondité}
\begin{equation}
R(A)=\frac{\underline{R}(A)}{\overline{R}(A)}\in[0,1].
\end{equation}
On dira que plus $R(A)$ est proche de $1$, plus $A$ est arrondi.
\end{itemize}
\end{definition}

Dans le cas d'une partie convexe et lorsqu'il est non nul, le supremum dans le calcul de $\underline{R}(A)$ est atteint par compacité, on parle alors d'une boule inscrite dans $A$, centrée sur un orthocentre. Lors de la suspension d'un polyèdre $\delta$ par rapport à un point $x$ non coplanaire, il est possible de donner des bornes sur les constantes de forme du polyèdre obtenu en fonction de la distance de $x$ à un orthocentre de $\delta$ et au sous-espace $\affine(\delta)$. La propriété suivante pourrait être donnée sous une forme plus précise mais elle suffira amplement pour la suite.

\begin{property}\label{propertysuspension}
Lorsque $\delta$ est un polyèdre, les sous-faces de $\S(\delta,x)$ sont de trois sortes:
\begin{itemize}
\item $\{x\}$ lui-même (lorsque $x\notin\delta$ ou $x\in\F_0(\delta)$);
\item des sous-faces de $\delta$;
\item des suspensions de sous-faces de $\delta$ par rapport à $x$.
\end{itemize}

Pour tout compact $K\subset]0,+\infty[^2$ il existe des constantes $c_1$ et $c_1$ strictement positives telles que pour tous polyèdre $\delta$ avec un orthocentre $o$ et $x\in\mathbb{R}^n\setminus\affine(\delta)$, si
\begin{equation}
\left(\frac{\dist(x,\affine(\delta))}{\dist(x,o)},\frac{\dist(x,o)}{\overline{R}(\delta)}\right)\in K
\end{equation}
alors
\begin{equation}
\overline{R}(\S(\delta,x))\leq c_1\overline{R}(\delta)\text{ et }\underline{R}(\S(\delta,x))\geq c_2\underline{R}(\delta).
\end{equation}
\end{property}

\begin{proof}
Le premier point est évident si l'on se réfère à la propriété~\ref{propertypolytope} et en particulier à la minimalité pour la convexité des sommets de $\S(\delta,x)$.

Pour vérifier le second point, notons:
\begin{itemize}
\item $\delta'=\S(\delta,x)$;
\item $H=\affine(\delta)$;
\item $o$ un orthocentre de $\delta$;
\item $\underline{B}$ une boule inscrite dans $\delta$ de centre $o$;
\item $\overline{B}$ la boule de centre $o$ et de rayon $2\overline{R}(\delta)$ (dès lors $\delta\subset\overline{B}$);
\item $o'$ un point du segment $[x,o]$.
\end{itemize}
Le problème peut se ramener à trouver une boule $B(o',\underline{r})$ incluse dans $\S\left(\underline{B},x\right)$, et une boule $B(o,\overline{r})$ contenant $\S\left(\overline{B},x\right)$. Appelons $\C$ le cône de sommet $x$ engendré par $\underline{B}$ (par hypothèse $\dist(x,H)>0$ donc $x\notin H$), et $\underline{B}'$ la plus grande boule de centre $o$ contenue dans $\C$. Son rayon $\underline{R}'$ vaut $\underline{R}(\delta)\cos\alpha$ où $\alpha$ est l'angle (non orienté) entre la normale à $H$ et la droite $(x,o)$. Or
\begin{equation}
\cos\alpha=\frac{\dist(x,H)}{\dist(x,o)}
\end{equation}
donc
\begin{equation}
\underline{R}'=\underline{R}(\delta)\frac{\dist(x,H)}{\dist(x,o)}.
\end{equation}
Les boules images de $\underline{B}'$ par l'homothétie de centre $x$ et de rapport $\gamma>0$ sont toutes contenues dans le cône $\C$, en particulier la boule $\underline{B}''$ de centre $c'$ obtenue avec $\gamma=\frac{1}{2}$ par exemple. Choisissons $o'$ comme le milieu du segment $[x,o]$ et posons
\begin{align}
\underline{r}_1&=\frac{1}{2}\cdot\underline{R'}=\frac{1}{2}\cdot\frac{\dist(x,H)}{\dist(x,o)}\cdot\underline{R}(\delta)&\underline{r}_2&=\dist(o',H)=\frac{\dist(x,H)}{2}.
\end{align}
Puisque par définition $\underline{R}(\delta)\leq\overline{R}(\delta)$ il vient encore
\begin{equation}
\underline{r}_2\geq\frac{1}{2}\cdot\frac{\dist(x,H)}{\dist(x,o)}\cdot\frac{\dist(x,o)}{\overline{R}(\delta)}\cdot\underline{R}(\delta).
\end{equation}

\begin{figure}
\begin{center}\includegraphics[width=0.5\textwidth]{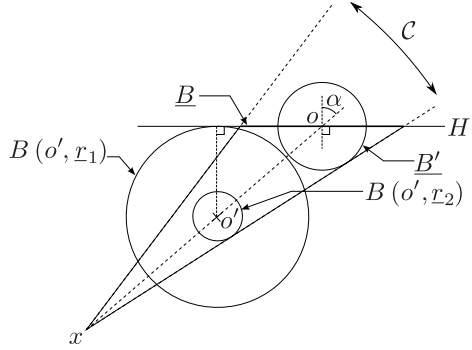}\end{center}
\caption{Rayon intérieur d'une suspension de polyèdre}
\label{figuresuspensionB}
\end{figure}

Par construction $B(o',\underline{r}_1)\subset\C$ et $B(o',\underline{r}_2)\subset\S(H,x)$, donc $B(o',\underline{r}_1)\cap B(o',\underline{r}_2)\subset\C\cap\S(H,x)\subset\delta'$, et finalement $\underline{r}=\min(\underline{r}_1,\underline{r}_2)$ convient. On obtient donc les bornes suivantes en fonction de $K$:
\begin{equation}
\begin{aligned}
\underline{R}(\delta')&\geq\frac{1}{2}\cdot\frac{\dist(x,H)}{\dist(x,o)}\min\left(1,\frac{\dist(x,o)}{\overline{R}(\delta)}\right)\underline{R}(\delta)\\
\overline{R}(\delta')&\leq\max\left(2,\frac{\dist(x,o)}{\overline{R}(\delta)}\right)\overline{R}(\delta).
\end{aligned}
\end{equation}
\end{proof}

\subsection{Complexes polyédriques et graphes}

Dans ce qui va suivre, pour en ensemble $S$ fini de polyèdres $k$-dimensionnels de $\mathbb{R}^n$ on notera:
\begin{itemize}
\item l'union des polyèdres
\begin{equation}
\U(S)=\bigcup_{\delta\in S}\delta;
\end{equation}
\item l'ensemble des sous-faces
\begin{equation}
\F(S)=\bigcup_{\delta\in S}\F(\delta);
\end{equation}
\item l'ensemble des sous-faces $k'$-dimensionnelles (pour $0\leq k'\leq k$)
\begin{equation}
\F_{k'}(S)=\bigcup_{\delta\in S}\F_{k'}(\delta);
\end{equation}
\item l'ensemble des faces de la frontière
\begin{equation}
\F_\partial(S)=\{\alpha\in\F_{k-1}(S)\colon\forall(\beta,\gamma)\in S^2,\alpha\neq\beta\cap\gamma\}.
\end{equation}
\end{itemize}
La définition suivante permet de formaliser l'idée intuitive de grilles de polyèdres qui se raccordent bien entre eux.

\begin{definition}[Complexes]\label{definitioncomplexe}
Lorsque toutes les sous-faces de $S$ sont d'intérieurs (relativement au sous-espace affine engendré correspondant) disjoints deux à deux, autrement dit si
\begin{equation}
\forall(\alpha,\beta)\in(\F(S))^2\colon\alpha\neq\beta\Rightarrow\open{\alpha}\cap\open{\beta}=\emptyset
\end{equation}
on dira que $S$ est un complexe $k$-dimensionnel.

Les constantes de forme de $S$ sont les extrema de celles de ses sous-faces de toute dimension et seront notées avec des lettres rondes:
\begin{align}
\overline{\R}(S)&=\max_{\delta\in\F(S)}\overline{R}(\delta)&\underline{\R}(S)&=\min_{\delta\in\F(S)\setminus\F_0(S)}\underline{R}(\delta)&\R(S)&=\min_{\delta\in\F(S)}R(\delta).
\end{align}
\end{definition}

On peut par exemple vérifier que pour tout polyèdre $\delta$ et $0\leq k\leq\dim\delta$ l'ensemble $\F_k(\delta)$ est un complexe, de même que $\F_{k'}(S)$ lorsque $S$ est un complexe $k$-dimen\-sion\-nel et $0\leq k'\leq k$. Dans le cas d'un complexe $n$-dimensionnel $S$, on a aussi par définition $\partial\U(S)=\U(\F_\partial(S))$.

Pour les besoins des constructions à venir on devra par ailleurs utiliser des ensembles finis de points munis d'une structure connective non orientée. Selon la terminologie usuelle on appellera de tels ensembles graphes, au sens de la définition suivante. On notera qu'il ne s'agit pas seulement ici de graphes abstraits, mais bien du plongement des objets correspondants dans $\mathbb{R}^n$.

\begin{definition}[Graphes]\label{definitiongraphe}
Un graphe $G=(T,A)$ est la donnée d'un couple formé de deux ensembles:
\begin{itemize}
\item l'ensemble des sommets, une partie finie non vide de $\mathbb{R}^n$ appelée support de $G$;
\item un ensemble d'arêtes $A\subset\P_2(T)$ qui contient des doublons de sommets distincts formant des segments ouverts disjoints deux à deux, c'est à dire que
\begin{equation}
\forall (\{a,b\},\{c,d\})\in A^2\colon\{a,b\}\neq\{c,d\}\Rightarrow ]a,b[\cap]c,d[=\emptyset.
\end{equation}
\end{itemize}
\end{definition}

Puisque les arêtes d'un graphe ainsi défini forment un complexe de dimension $1$ (car supposées d'intérieurs disjoints), on se réserve le droit d'user parfois d'une terminologie identique pour les deux types d'objets. Lorsqu'une arête $\{x,y\}\in A$ on dira que les sommets $x$ et $y$ sont voisins, et on confondra souvent cette arête avec le segment $[x,y]$. Pour tout sommet $x\in T$ on appelle ordre de $x$ le nombre de ses voisins:
\begin{equation}
\O(x)=\#\{y\in T\colon\{x,y\}\in A\}.
\end{equation}
Cette définition peut encore s'étendre au graphe tout entier:
\begin{equation}
\O(G)=\max_{x\in T}\O(x).
\end{equation}

Dotons-nous encore de la terminologie qui suit afin de décrire la structure des graphes:
\begin{itemize}
\item pour $k>1$ on dira qu'un chemin de longueur $k$ est un $k+1$-uplet de sommets $(x_1,\ldots,x_{k+1})\in T^{k+1}$ tel que $\forall i\in\{1,\ldots,k\}\colon\{x_i,x_{i+1}\}\in A$;
\item un $k$-cycle est un chemin de longueur $k$ dont les deux extrémités sont égales, et les $k-1$ autres sommets le composant distincts deux à deux. Par exemple, si $\{x,y\}\in A$ le triplet $(x,y,x)$ est un $2$-cycle;
\item on dira que $G$ est connexe s'il existe un chemin qui le parcourt en entier;
\item on dira que $G$ est linéaire si $\O(G)\leq2$, s'il ne possède aucun $3$-cycle et s'il est connexe;
\item on dira que $G$ est cyclique s'il existe un cycle qui passe par tous ses sommets.
\end{itemize}

\subsection{Suspension de complexes par rapport à un graphe linéaire}

On suppose qu'on dispose d'un complexe $k$-dimensionnel $S$, d'un graphe linéaire $G=(T,A)$ et d'une application $p$ de $S$ dans $T$, appelée choix de suspension. Posons
\begin{equation}
\begin{split}
S'&=\left\{\alpha\cap\beta\colon\{p(\alpha),p(\beta)\}\in A\right\}\\
S^*&=S\cup S'
\end{split}
\end{equation}
et pour $\delta\in S^*$ on notera
\begin{equation}
p^*(\delta)=
\begin{cases}
\{p(\delta)\}\text{ si $\delta\in S$}\\
\{p(\alpha),p(\beta)\}\text{ si $\delta\in S'$}.
\end{cases}
\end{equation}
Avec ces notations on va donner une définition pour la suspension de $S$ par rapport à $G$.

\begin{definition}[Suspension de complexe]\label{definitionsuspensioncomplexe}
La suspension $\S(S,G,p)$ du complexe $k$-dimensionnel $S$ par rapport au graphe linéaire $G$ selon le choix $p$ est l'ensemble des polyèdres $\left<\delta\cup p^*(\delta)\right>$ de dimension $k+1$ obtenus lorsque $\delta$ parcourt $S^*$:
\begin{equation}
\S(S,G,p)=\left\{\delta'=\left<\delta\cup p^*(\delta)\right>\colon\delta\in S^*\text{ et }\dim\delta'=k+1\right\}.
\end{equation}
Lorsque $\S(S,G,p)$ est un complexe on dira que $p$ est un choix adapté à la suspension.
\end{definition}

Il est clair qu'une suspension de complexe n'est en général pas un complexe, on peut donner l'exemple simple de deux polyèdres et d'un graphe à un seul sommet situé sur l'origine d'une demi-droite qui intersecte l'intérieur de chacun des polyèdres: dans ce cas les suspensions respectives des deux polyèdres ne sont pas d'intérieurs disjoints.

En fait, il est nécessaire que toute demi-droite dont l'origine est l'un des sommets du graphe ne rencontre au maximum qu'une seule fois l'union des polyèdres mis en correspondance avec ce sommet par le choix $p$. Pour formaliser cette idée, pour deux parties $A$ et $B$ de $\mathbb{R}^n$ on dira que $A$ est en occlusion simple par rapport à $B$ si
\begin{equation}
\forall (x,y)\in A\times B\colon[x,y]\cap A=\{x\}.
\end{equation}
En particulier lorsque $A$ est la frontière d'un ouvert borné $U$ cette propriété est équivalente au fait que $U$ soit étoilé par rapport à tout point de $A$. De manière plus générale on se donne aussi la définition suivante dans le cas de la suspension d'un complexe $k$-dimensionnel $S$ par rapport à un graphe $G=(T,A)$.

\begin{figure}
\begin{center}\includegraphics[width=0.9\textwidth]{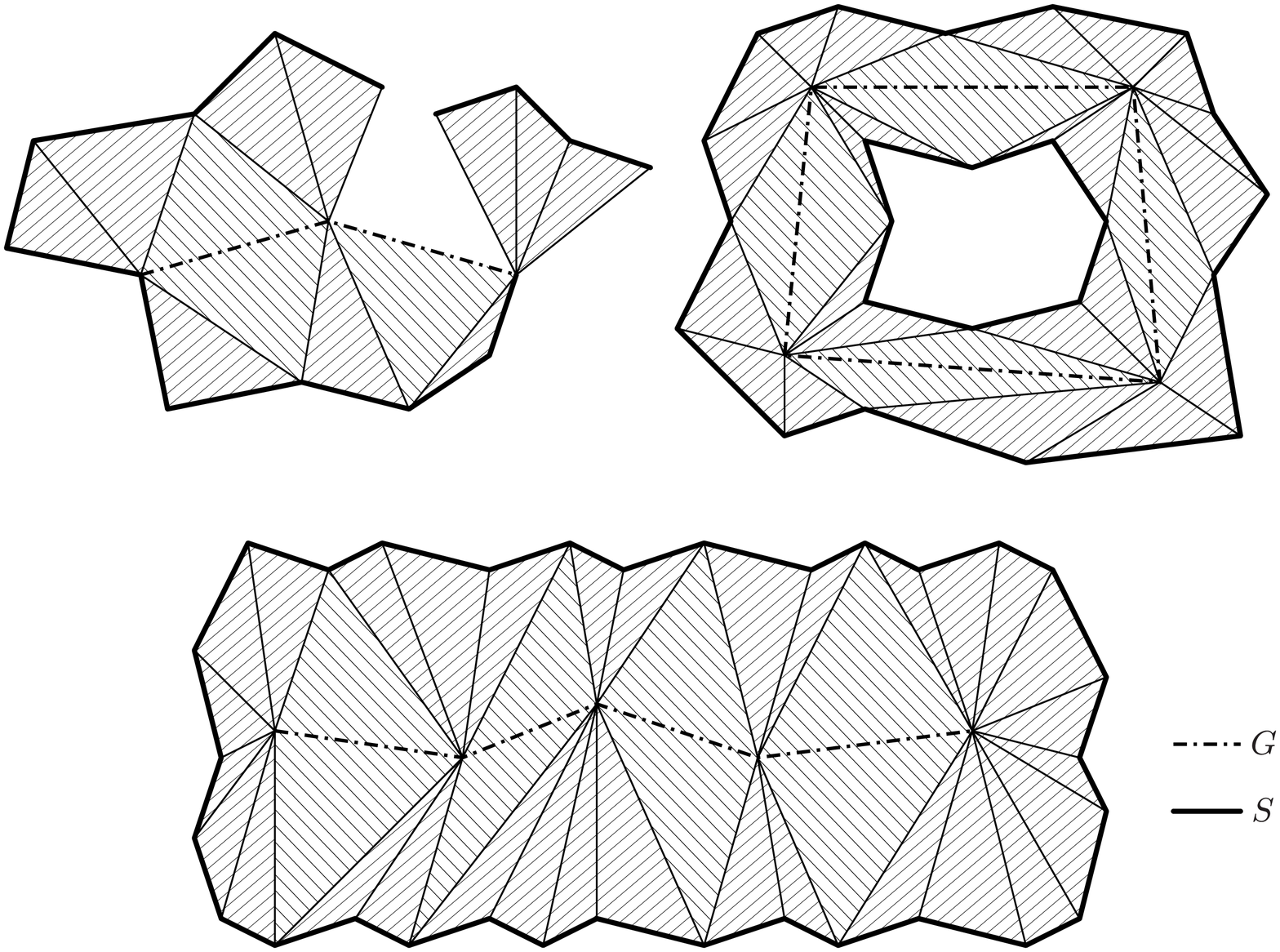}\end{center}
\caption{Exemples de suspensions de complexes avec choix adapté; les exemples du bas et de droite vérifient en plus les hypothèses du lemme~\ref{lemmaocclusiontotale} d'occlusion totale, le second avec un graphe cyclique}
\label{suspensioncomplexeA}
\end{figure}

\begin{definition}[Occlusion simple]\label{definitionocclusionsimple}
On dira que $S$ est en occlusion simple par rapport à $G$ avec le choix $p$ si les trois conditions suivantes sont réalisées:
\begin{enumerate}
\item les polyèdres de $S$ ne rencontrent pas les arêtes de $G$;
\item pour tout sommet $x\in T$, l'ensemble
\begin{equation}
\U(p^{-1}(x))=\bigcup_{\delta\in S\colon p(\delta)=x}\delta
\end{equation}
est en occlusion simple par rapport à $x$;
\item pour toute arête $\{x,y\}\in A$:
\begin{itemize}
\item $\U(p^{-1}(x))\cap\U(p^{-1}(y))$ est en occlusion simple par rapport au segment $[x,y]$;
\item $\U(p^{-1}(x))\cup\S(\U(p^{-1}(x))\cap\U(p^{-1}(y)),y)$ est en occlusion simple par rapport au point $x$.
\end{itemize} 
\end{enumerate}
\end{definition}

Avec cette définition on va être en mesure de donner des conditions suffisantes pour que la suspension d'un complexe en occlusion simple par rapport à un graphe linéaire soit encore un complexe.

\begin{lemma}[Condition suffisante de suspension adaptée]\label{lemmaocclusionsimple}
Soient $S$ un complexe $k$-dimensionnel, $G=(T,A)$ un graphe linéaire et $p$ un choix de suspension vérifiant les propriétés suivantes:
\begin{itemize}
\item $S$ est en occlusion simple par rapport à $G$ avec le choix $p$;
\item il existe une famille $(\kappa_x)_{x\in T}$ d'ouverts deux à deux disjoints de $\mathbb{R}^n$ tels que:
\begin{itemize}
\item $\forall x\in T$, $\kappa_x$ est étoilé par rapport à $x$ et $\U(p^{-1}(x))\subset\overline{\kappa_x}$;
\item si $x$ et $y$ sont deux sommets voisins de $G$ alors en notant
\begin{equation}
\kappa_{x,y}=\bigcup_{(u,v)\in\kappa_x\times\kappa_y}[u,v]
\end{equation}
pour tout $z\in T\setminus\{x,y\}$ on a $\kappa_z\cap\kappa_{x,y}=\emptyset$;
\item si $\{a,b\}$ et $\{c,d\}$ sont deux arêtes distinctes du graphe alors $\kappa_{a,b}\cap\kappa_{c,d}=\emptyset$ (cette condition impliquant la précédente lorsque $G$ est cyclique).
\end{itemize}
\end{itemize}

Alors $\S(S,G,p)$ est un complexe $k+1$-dimensionnel.
\end{lemma}

\begin{proof}
Par définition les polyèdres de $\S(S,G,p)$ sont de deux espèces:
\begin{enumerate}
\item ceux obtenus par suspension d'un polyèdre de $S$ par rapport à un sommet de $T$;
\item ceux obtenus par suspension successive d'une face de $S$ par rapport à deux sommets voisins de $T$.
\end{enumerate}
Vérifions que l'ensemble des polyèdres de ces deux espèces forme bien un complexe. Considérons $F$ une sous-face d'un polyèdre $\delta$ créé après suspension, elle peut être obtenue de quatre manières différentes:
\begin{enumerate}
\item $\delta$ est de première ou seconde espèce et $F$ est une sous-face de $S$;
\item $\delta$ est de première ou seconde espèce et $F=\left<F'\cup\{x\}\right>$ où $F'$ est une sous-face de $S$ et $x$ un sommet de $G$;
\item $\delta$ est de seconde espèce et $F=\left<F'\cup\{x,y\}\right>$ où $F'$ est une sous-face de $S$ et $x$ et $y$ deux sommets voisins de $G$;
\item $\delta$ est de seconde espèce et $F$ est une arête de $G$.
\end{enumerate}
Notons respectivement $S_1$, $S_2$, $S_3$ et $S_4$ les ensembles des sous-faces obtenues après suspension et correspondant respectivement à ces quatre sortes (les ensembles $S_i$ sont deux à deux disjoints, puisque par construction une sous-face ne peut être que d'une seule sorte à la fois), ainsi
\begin{equation}
\F(\S(S,G,p))=S_1\sqcup S_2\sqcup S_3\sqcup S_4.
\end{equation}
On va montrer que les polyèdres de $\F(\S(S,G,p))$ sont d'intérieurs disjoints deux à deux. Pour cela considérons (lorsque cela est possible\footnote{Il se peut que l'un des $S_i$ soit vide ou réduit à un seul polyèdre; dans ce cas-là encore on obtient ce qu'on cherchait, à savoir que les polyèdres de $\bigcup_i S_i$ sont d'intérieurs disjoints deux à deux.}) pour $1\leq i\leq4$ deux polyèdres distincts $F_i$ et $G_i$ de $S_i$ et démontrons que
\begin{equation}
\forall (i,j)\in\{1,2,3,4\}^2\colon\open{F_i}\cap\open{G_j}=\emptyset.
\end{equation}
Il y a un total de seize cas à considérer, qu'on peut ramener à dix à une permutation des notations près, et qu'on va traiter un à un en utilisant ces notations:
\begin{itemize}
\item $F_1$ est une sous-face d'un polyèdre $\alpha_1$ de $S$
\item $F_2=\left<F'_2\cup\{x\}\right>$ où $F'_2$ est une sous face de $\alpha_2\in S$;
\item $G_2=\left<G'_2\cup\{y\}\right>$ où $G'_2$ est une sous face de $\beta_2\in S$;
\item $F_3=\left<F'_3\cup\{a,b\}\right>$ où $F'_3$ est une sous face de $\alpha_3\in S$;
\item $G_3=\left<G'_3\cup\{c,d\}\right>$ où $G'_3$ est une sous face de $\beta_3\in S$;
\item $G_4=[u,v]$.
\end{itemize}

\begin{description}
\itemx11$F_1$ et $G_1$ sont d'intérieurs disjoints par définition puisque $S$ est un complexe.

\itemx12remarquons d'abord que si $F_1$ est une sous-face de $G'_2$ elle est aussi une sous-face de $G_2$ (car $y$ est en occlusion simple par rapport à $G'_2$) et donc $\open{F_1}\cap\open{G_2}=\emptyset$. Supposons alors que $F_1$ ne soit pas une sous-face de $G'_2$, posons $z=p(\alpha_1)$ et considérons les deux cas:
\begin{itemize}
\item si $z\neq y$ alors par hypothèse d'abord $F_1\subset\alpha_1\subset\overline{\kappa_z}$ et $G'_2\subset\beta_2\subset\overline{\kappa_y}$ donc $G_2\subset\overline{\kappa_y}$ car $\kappa_y$ est étoilé par rapport à $y$. De plus puisque $y\notin\overline{\kappa_z}$ alors $G_2\cap\partial\kappa_z\subset G'_2$. Dès lors puisque $\kappa_y$ et $\kappa_z$ sont disjoints il vient $F_1\cap G_2=F_1\cap\overline{\kappa_z}\cap G_2\cap\overline{\kappa_y}=F_1\cap G_2\cap\partial\kappa_y\cap\partial\kappa_z\subset F_1\cap G_2\cap\partial\kappa_y\subset F_1\cap G'_2$;
\item si $z=y$, puisque par hypothèse $y$ est en occlusion simple par rapport à $\beta_2\cup F_1$ donc aussi par rapport à $G'_2\cup F_1$ alors $F_1\cap G_2=F_1\cap G'_2$.
\end{itemize}
Dans les deux cas, $G'_2$ est une sous-face stricte de $G_2$ et puisque l'on a supposé que $F_1$ n'est pas une sous face de $G'_2$ alors $F_1\cap G'_2$ est soit vide, soit une sous-face stricte de $F_1$. Pour conclure, $F_1\cap G_2$ est donc soit vide, soit une sous-face stricte de $F_1$ et $G_2$ (donc incluse dans leur frontière) et par conséquent $\open{F_1}\cap\open{G_2}=\emptyset$.

\itemx13remarquons d'abord que si $F_1$ est une sous-face de $G'_3$, elle est aussi une sous-face de $G_3$ (car aucun point de $[c,d]$ n'est dans le sous-espace affine engendré par $G'_3$ par hypothèse d'occlusion simple) et donc $\open{F_1}\cap\open{G_3}=\emptyset$. Supposons donc que $F_1$ ne soit pas une sous-face de $G'_3$ et posons $z=p(\alpha_1)$. Quitte à permuter les notations, il y a deux cas à envisager:
\begin{itemize}
\item si $z\neq c$ et $z\neq d$ alors de la même façon que pour le cas précédent on a que $G_3\subset\overline{\kappa_{c,d}}$ et puisque $[c,d]\subset\kappa_{c,d}$ alors $G_3\cap\partial\kappa_{c,d}\subset G'_3$. Puisque par hypothèse $\kappa_{c,d}$ et $\overline{\kappa_z}$ sont disjoints alors $F_1\cap G_3=F_1\cap\overline{\kappa_z}\cap G_3\cap\overline{\kappa_{c,d}}=F_1\cap G_3\cap\partial\kappa_z\cap\partial\kappa_{c,d}\subset F_1\cap G_3\cap\partial\kappa_{c,d}\subset F_1\cap G'_3$;
\item si $z=c$, puisque $[c,d]$ est en occlusion simple par rapport à $\beta_3\cup\alpha_1$ donc aussi par rapport à $G'_3\cup F_1$ alors $F_1\cap G_3=F_1\cap G'_3$.
\end{itemize}
Dans les deux cas, $G'_3$ est une sous-face stricte de $G_3$ et puisque l'on a supposé que $F_1$ n'est pas une sous-face de $G'_3$ alors $F_1\cap G'_3$ est soit vide, soit une sous-face stricte de $F_1$. Pour conclure, $F_1\cap G_3$ est donc soit vide, soit une sous-face stricte de $F_1$ et $G_3$ (donc incluse dans leur frontière) et par conséquent $\open{F_1}\cap\open{G_3}=\emptyset$.

\itemx14par hypothèse d'occlusion simple les arêtes du graphe ne ren\-con\-trent pas les polyèdres du complexe donc a fortiori $\open{F_1}\cap\open{G_4}=\emptyset$.

\itemx22considérons les deux cas possibles:
\begin{itemize}
\item si $x\neq y$ alors comme vu dans le cas $(i,j)=(1,2)$, $F_2\cap\partial\kappa_x\subset F'_2$ et $G_2\cap\partial\kappa_y\subset G'_2$ avec par hypothèse $\kappa_x\cap\kappa_y=\emptyset$. Par conséquent, puisque $\kappa_x\cap\kappa_y=\emptyset$ alors $F_2\cap G_2=F_2\cap\overline{\kappa_x}\cap G_2\cap\overline{\kappa_y}=F_2\cap\partial\kappa_x\cap G_2\cap\partial\kappa_y\subset F'_2\cap G'_2$. Or, puisque $F'_2$ et $G'_2$ sont deux sous-faces strictes de $F_2$ et $G_2$, alors $F_2\cap G_2$ est soit vide, soit une sous-face stricte de $F_2$ et $G_2$ (donc incluse dans leurs frontières) et par conséquent $\open{F_2}\cap\open{G_2}=\emptyset$;
\item si $x=y$ alors $F_2\cap G_2=\left<(F'_2\cap G'_2)\cup\{x\}\right>$. Remarquons ensuite que par hypothèse $F'_2\neq G'_2$ (sinon $F_2=G_2$), et si $F'_2$ (respectivement $G'_2$) est une sous face stricte de $G'_2$ (respectivement $F'_2$) alors $F_2$ (respectivement $G_2$) est une sous-face stricte de $G_2$ (respectivement $F_2$) et alors $\open{F_2}\cap\open{G_2}=\emptyset$. Supposons alors que $F'_2$ et $G'_2$ ne soient pas sous-face de respectivement $G'_2$ et $F'_2$, alors $F'_2\cap G'_2$ est soit vide, soit une sous-face stricte de $F'_2$ et $G'_2$. Par conséquent, $F_2\cap G_2$ est une sous-face stricte de $F_2$ et $G_2$, donc incluse dans leurs frontières et par conséquent $\open{F_2}\cap\open{G_2}=\emptyset$.
\end{itemize}

\itemx23quitte à permuter les notations, il reste deux cas à étudier:
\begin{itemize}
\item si $x\neq c$ et $x\neq d$ alors comme vu dans les cas précédents, $F_2\cap\partial\kappa_x\subset F'_2$ et $G_3\cap\partial\kappa_{c,d}\subset G'_3$. De plus par hypothèse $\kappa_x\cap\kappa_{c,d}=\emptyset$ d'où $F_2\cap G_3=F_2\cap\overline{\kappa_x}\cap G_3\cap\overline{\kappa_{c,d}}=F_2\cap\partial\kappa_x\cap G_3\cap\partial\kappa_{c,d}\subset F'_2\cap G'_3$. Or $F'_2$ et $G'_3$ sont des sous-faces strictes de respectivement $F_2$ et $G_3$, par conséquent $F_2\cap G_3$ est soit vide, soit une sous-face stricte de $F_2$ et $G_3$ (donc incluse dans leurs frontières) et finalement $\open{F_2}\cap\open{G_3}=\emptyset$;
\item si $x=c$ alors d'après la dernière hypothèse de la définition de l'occlusion simple d'un complexe on sait que $F'_2\cup\left<G'_3\cup\{y\}\right>$ est en occlusion simple par rapport à $x$. On en tire que $F_2\cap G_3=\left<F'_2\cup\{x\}\right>\cap\left<G'_3\cup\{x,y\}\right>=\left<(F'_2\cap\left<G'_3\cup\{y\}\right>)\cup\{x\}\right>$. Par ailleurs, on sait aussi que $F'_2\subset\overline{\kappa_x}$ et $G'_3\subset\overline{\kappa_x}\cap\overline{\kappa_y}$, avec $\kappa_x\cap\kappa_y=\emptyset$ et $\kappa_{y}$ qui est étoilé par rapport à $y$. Il vient donc $F'_2\cap\left<G'_3\cup\{y\}\right>=F'_2\cap G'_3$, d'où on tire $F_2\cap G_3=\left<(F'_2\cap G'_3)\cup\{x\}\right>$. Remarquons encore que si $F'_2$ est une sous-face de $G'_3$ alors $F_2$ est une sous-face de $G_3$ et donc $\open{F_2}\cap\open{G_3}=\emptyset$. Supposons alors que $F'_2$ ne soit pas une sous-face de $G'_3$, dans ce cas $F'_2\cap G'_3$ est soit vide, soit une sous-face stricte de $F'_2$, donc $\left<(F'_2\cap G'_3)\cap\{x\}\right>$ est une sous-face stricte de $F_2$. Puisque c'est aussi une sous-face stricte de $G_3$ alors $F_2\cap G_3$ est incluse dans la frontière de $F_2$ et $G_3$ et par conséquent $\open{F_2}\cap\open{G_3}=\emptyset$.
\end{itemize}

\itemx24remarquons que si $u$ et $v$ sont tous deux différents de $x$, puisque $F_2\subset\overline{\kappa_x}$ et que $\kappa_{u,v}\cap\kappa_x=\emptyset$ alors $F_2\cap [u,v]=\emptyset$ puisque $[u,v]\subset\kappa_{u,v}$. Supposons ensuite que $u=x$ par exemple, et raisonnons par l'absurde en supposant que $F_2\cap [u,v]\neq\{u\}$. Dans ce cas $F_2$ est la suspension de $F'_2$ par rapport à $u$ avec $F'_2\subset\overline{\kappa_u}$. En outre $\kappa_u$ étant étoilé par rapport à $u$ et $v\notin\kappa_u$ on trouve que forcément $[u,v]\cap F'_2\neq\emptyset$, ce qui contredit l'hypothèse que les arêtes du graphe ne rencontrent pas les polyèdres de $S$.

\itemx33puisque $F_3$ (respectivement $G_3$) est de seconde espèce, par définition de la suspension d'un complexe, $F'_3$ (respectivement $G'_3$) est incluse dans l'intersection de deux polyèdres $\delta_1$ et $\delta_2$ de $S$ (respectivement $\delta_3$ et $\delta_4$) tels que $p(\delta_1)=a$ et $p(\delta_2)=b$ (respectivement $p(\delta_3)=c$ et $p(\delta_4)=d$). Remarquons alors que si $F'_3$ (respectivement $G'_3$) est une sous-face de $G'_3$ (respectivement $F'_3$) alors parmi tous les $\delta_i$ pouvant convenir on peut choisir $\delta_1=\delta_3$ et $\delta_2=\delta_4$ et donc $\{a,b\}=\{c,d\}$. Dès lors, $F_3$ (respectivement $G_3$) est une sous-face de $G_3$ (respectivement $F_3$) et donc par conséquent $\open{F_3}\cap\open{G_3}=\emptyset$ ou $F_3=G_3$. On supposera donc que $F'_3\cap G'_3$ est soit vide, soit une sous-face stricte de $F'_3$ et de $G'_3$. Quitte à permuter les notations il reste là encore trois cas à traiter séparément:
\begin{itemize}
\item si $\{a,b\}\cap\{c,d\}=\emptyset$ alors par hypothèse $\kappa_{a,b}\cap\kappa_{c,d}=\emptyset$, de plus comme vu dans les cas précédents, puisque $F_3\cap\partial\kappa_{a,b}\subset F'_3$ et $G_3\cap\partial\kappa_{c,d}\subset G'_3$ alors $F_3\cap G_3=F_3\cap\overline{\kappa_{a,b}}\cap G_3\cap\overline{\kappa_{c,d}}=F_3\cap\partial\kappa_{a,b}\cap G_3\cap\partial\kappa_{c,d}\subset F'_3\cap G'_3$;
\item si $a=c$ et $b\neq d$, puisque par hypothèse $F'_3\subset\overline{\kappa_b}$ et $G'_3\subset\overline{\kappa_d}$ avec $\kappa_b\cap\kappa_d=\emptyset$, $b\in\kappa_b$ et $d\in\kappa_d$ alors déjà $\left<F'_3\cup\{b\}\right>\cap\left<G'_3\cup\{d\}\right>=F'_3\cap G'_3$. De plus $\kappa_{a,b}\cap\kappa_d=\emptyset$ et $\kappa_{a,d}\cap\kappa_b=\emptyset$ avec $F_3\subset\overline{\kappa_{a,b}}$ et $G_3\subset\overline{\kappa_{a,d}}$ donc $F_3\cap G_3=\left<F'_3\cup\{a,b\}\right>\cap\left<G'_3\cup\{a,d\}\right>=\left<F'_3\cup\{a\}\right>\cap\left<G'_3\cup\{a\}\right>=\left<(F'_3\cap G'_3)\cup \{a\}\right>$. Or puisque $b\notin\kappa_a$ et $d\notin\kappa_a$ alors $\left<(F'_3\cap G'_3)\cup\{a\}\right>$ est une sous-face stricte de $F_3$ et $G_3$;
\item si $a=c$ et $b=d$ alors puisque par hypothèse $[a,b]$ est en occlusion simple par rapport à $\alpha_3\cup\beta_3$, $F_3\cap G_3=\left<(F'_3\cap G'_3)\cup [a,b]\right>$. Or $F'_3\cap G'_3$ est soit vide, soit une sous-face stricte de $F'_3$ ou $G'_3$ donc $F_3\cap G_3$ est une sous-face stricte de $F_3$ et $G_3$.
\end{itemize}
Dans les trois cas, $F_3\cap G_3$ est une sous-face stricte de $F_3$ et de $G_3$ et par conséquent $\open{F_3}\cap\open{G_3}=\emptyset$.

\itemx34quitte à échanger les notations il y a trois cas à envisager:
\begin{itemize}
\item si $u=a$ et $v=b$ alors $[u,v]$ est une sous-face de $F_3$ et donc $\open{F_3}\cap\open{G_4}=\emptyset$;
\item si $u=a$ et $v\neq b$, puisque $F_3\subset\overline{\kappa_{u,b}}$ alors par hypothèse $v\notin F_3$. En écrivant que $F_3=\left<F'_3\cup\{u,b\}\right>=\left<\left<F'_3\cup\{b\}\right>\cup\{u\}\right>$, raisonnons par l'absurde et supposons que $]u,v[\cap F_3\neq\emptyset$. Puisque $F_3$ est la suspension de $\left<F'_3\cup\{b\}\right>$ par rapport à $u$ et que $v\notin F_3$, le segment $]u,v[$ intersecte alors $\left<F'_3\cup\{b\}\right>$ en au moins un point. Par hypothèse, $\kappa_{u,v}\cap\kappa_b=\emptyset$ avec $\kappa_{u,v}$ et $\kappa_b$ ouverts, donc on a aussi que $\kappa_{u,v}\cap\overline{\kappa_b}=\emptyset$, en outre $\left<F'_3\cup\{b\}\right>\subset\overline{\kappa_b}$ car $\kappa_b$ est étoilé par rapport à $b$. En résumé, en supposant que $]u,v[\cap F_3\neq\emptyset$ on obtient quatre comparaisons ensemblistes contradictoires:
\begin{equation}
\begin{aligned}
]u,v[&\subset\kappa_{u,v}&\kappa_{u,v}\cap\overline{\kappa_b}&=\emptyset\\
\left<F'_3\cup\{b\}\right>&\subset\overline{\kappa_b}&]u,v[\cap\left<F'_3\cup\{b\}\right>&\neq\emptyset.
\end{aligned}
\end{equation}
Par conséquent, $]u,v[\cap F_3$ est nécessairement vide; 
\item si $u\neq a$ et $v\neq b$ alors puisque $[u,v]\cap\overline{\kappa_{a,b}}=\emptyset$ (car $[u,v]$ est inclus dans l'ouvert $\kappa_{u,v}$ et $\kappa_{u,v}\cap\kappa_{a,b}=\emptyset$) et $F_3\subset\overline{\kappa_{a,b}}$ (car $\kappa_a$ et $\kappa_b$ étoilés par rapport respectivement à $a$ et $b$) on obtient que $[u,v]\cap F_3=\emptyset$.
\end{itemize}
Dans les trois cas, on trouve bien que $\open{F_3}\cap\open{G_4}=\emptyset$.

\itemx44par définition d'un graphe ses arêtes sont d'intérieurs disjoints donc $\open{F_4}\cap\open{G_4}=\emptyset$.
\end{description}

Ceci termine la démonstration que l'ensemble de polyèdres $\S(S,G,p)$ est un complexe $k+1$-dimensionnel, puisqu'on a démontré que ses sous-faces sont d'intérieurs disjoints deux à deux.
\end{proof}

Dans ce qui suit, on va devoir suspendre un complexe qui forme la frontière d'un ouvert de façon à remplir de polyèdres $n$-dimensionnel tout le volume de cet ouvert. La notion d'occlusion totale va permettre de caractériser une situation où le remplissage se fait automatiquement.

\begin{definition}[Occlusion totale]\label{definitionocclusiontotale}
On dira qu'un graphe linéaire $G$ est en occlusion totale par rapport à un complexe $n-1$-dimensionnel $S$ avec le choix $p$ si $p$ est surjectif et s'il existe un ouvert borné $O$ et une famille d'ouverts deux à deux disjoints $(\kappa_x)_{x\in T}$ vérifiant les trois conditions suivantes:
\begin{enumerate}
\item $\U(S)=\partial O$ et $\overline{O}=\bigcup_{x\in T}\overline{\kappa_x}$;
\item $\forall x\in T\colon\kappa_x$ est étoilé par rapport à $x$ et $\U(p^{-1}(x))\subset\overline{\kappa_x}$;
\item $\forall (x,y)\in T^2\colon\{x,y\}\notin A\Rightarrow x=y\text{ ou }\overline{\kappa_x}\cap\overline{\kappa_y}=\emptyset$.
\end{enumerate}
\end{definition}

Le lemme suivant justifie qu'en situation d'occlusion totale, la suspension va remplir tout l'espace vide à l'intérieur de l'ouvert $O$.

\begin{lemma}[Remplissage en occlusion totale]\label{lemmaocclusiontotale}
Si le graphe linéaire $G$ est en occlusion totale par rapport au complexe $S$ avec le choix $p$, et si $p$ est un choix adapté à la suspension alors $S'=\S(S,G,p)$ est un complexe qui vérifie
\begin{equation}
\overline{O}\subset\U(S').
\end{equation}

Si $S$ vérifie aussi les hypothèses du lemme~\ref{lemmaocclusionsimple} avec les mêmes ouverts $\kappa_x$ (qui, rappelons-le, établissent une condition suffisante pour un choix adapté) l'inclusion inverse est elle aussi vérifiée:
\begin{equation}
\overline{O}=\U(S').
\end{equation}
\end{lemma}

\begin{proof}
D'abord on sait que $S'$ est un complexe, si de plus $G$ est réduit à un seul sommet le lemme est évident. Supposons donc à présent que ce n'est pas le cas.

Par définition, un graphe linéaire non réduit à un seul sommet est tel que chaque sommet a un ou deux voisins distincts. Dans le second cas, en notant $x$ un sommet de $G$ et $y$ et $z$ ses deux voisins, puisque les $\left(\kappa_u\right)_{u\in T}$ sont disjoints deux à deux alors $\partial\kappa_x=(\partial O\cap\overline{\kappa_x}\setminus(\overline{\kappa_y}\cup\overline{\kappa_z}))\cup(\partial\kappa_y\cap\overline{\kappa_x})\cup(\partial\kappa_z\cap\overline{\kappa_x})$. De plus par hypothèse on a aussi que $\overline{\kappa_y}\cap\overline{\kappa_z}=\emptyset$ car $y$ et $z$ ne peuvent être voisins (sinon $(y,x,z,y)$ formerait un $3$-cycle) donc $\partial\kappa_y\cap\partial\kappa_z=\emptyset$. Cette union est par conséquent une union d'ensembles disjoints deux à deux:
\begin{equation}\label{equationocclusionsimpleA}
\partial\kappa_x=(\partial O\cap\overline{\kappa_x}\setminus(\overline{\kappa_y}\cup\overline{\kappa_z}))\sqcup(\partial\kappa_y\cap\overline{\kappa_x})\sqcup(\partial\kappa_z\cap\overline{\kappa_x}).
\end{equation}
De même si $x$ a un seul voisin $y$ alors
\begin{equation}
\partial\kappa_x=(\partial O\cap\overline{\kappa_x}\setminus\overline{\kappa_y})\sqcup(\partial\kappa_y\cap\overline{\kappa_x}).
\end{equation}

À présent, si on choisit $a\in O$ on peut alors trouver $x\in T$ tel que $a\in\overline{\kappa_x}$ (en effet, $ O\subset\overline{\bigcup_{x\in T}\kappa_x}$). Il nous faut considérer l'ordre du sommet $x$: supposons d'abord que $\O(x)=2$ et notons $y$ et $z$ ses deux voisins distincts.

Si $x=a$ alors puisque $p$ est surjectif on peut trouver $\delta\in S$ tel que $p(\delta)=x$, de plus par définition $\left<\delta\cup\{x\}\right>\in\S(S,G,p)$ donc $a\in\U(S')$.

Si $x\neq a$ alors puisque $x,a\in\overline{\kappa_x}$ et que $\kappa_x$ est un ouvert borné, on peut donc trouver $b\in\partial\kappa_x$ tel que $a\in[b,x[$. Et puisque $\kappa_x$ est étoilé par rapport à $x$ alors on a aussi que $]b,x]\subset\kappa_x$. D'après~\eqref{equationocclusionsimpleA} il nous faut considérer trois cas possibles:
\begin{itemize}
\item si $b\in\partial O\cap\overline{\kappa_x}\setminus(\overline{\kappa_y}\cup\overline{\kappa_z})$ alors puisque $\U(S)=\partial O$ on peut trouver $\delta\in S$ tel que $b\in\delta$. Remarquons aussi que puisque $b\notin\overline{\kappa_y}\cup\overline{\kappa_z}$ alors nécessairement $p(\delta)=x$ (en effet si $p(\delta)=y$ par exemple alors par hypothèse on aurait que $\delta\subset\overline{\kappa_y}$ ce qui est impossible car $b\in\delta$). Donc $\left<\delta\cup\{x\}\right>\in\S(S,G,p)$, or $[b,x]\subset\left<\delta\cup\{x\}\right>$ par conséquent $a\in\left<\delta\cup\{x\}\right>\in S'$;
\item si $b\in\partial\kappa_y\cap\overline{\kappa_x}$ alors notons $c$ l'intersection de la demi-droite d'origine $x$ et de direction $\overrightarrow{yx}$ avec $\partial\kappa_x$. Cette intersection existe car $\kappa_x$ est un ouvert borné et $x\in\kappa_x$, et est unique car $\kappa_x$ est étoilé par rapport à $x$. Puisque $\kappa_y$ est étoilé par rapport à $y$ alors nécessairement $c\notin\partial\kappa_y$ (sinon on devrait avoir $]c,x[\subset\kappa_y$, ce qui n'est pas le cas car $]c,x[\subset\kappa_x$).

À présent si on note $H$ un $2$-plan affine contenant les trois points $x$, $y$ et $a$, puisque $c\in H$ il vient
\begin{equation}
H\cap\partial\kappa_x\neq H\cap\partial\kappa_y\cap\overline{\kappa_x}.
\end{equation}
En faisant l'intersection des membres de l'égalité~\eqref{equationocclusionsimpleA} avec le $2$-plan $H$, on obtient
\begin{equation}
H\cap\partial\kappa_x=O'\sqcup F_1\sqcup F_2
\end{equation}
avec
\begin{equation}
\begin{aligned}
O'&=H\cap\partial O\cap\overline{\kappa_x}\setminus(\overline{\kappa_y}\cup\overline{\kappa_z})\\
F_1&=H\cap\partial\kappa_y\cap\overline{\kappa_x}\\
F_2&=H\cap\partial\kappa_z\cap\overline{\kappa_x}.
\end{aligned}
\end{equation}

En résumé, en se plaçant dans $H$ on a montré que le fermé $H\cap\partial\kappa_x$ est l'union disjointe de l'ouvert $O'$ et des deux fermés non vides $F_1$ et $F_2$, avec $F_1\neq H\cap\partial\kappa_x$. Il est facile de constater en plus que $H\cap\partial\kappa_x$ est connexe, puisque c'est la frontière de l'ouvert $H\cap\kappa_x$ de $H$, qui est borné et étoilé par rapport à $x$. Par conséquent l'intersection de $\overline{O'}$ et de $F_1$ n'est pas vide, en d'autres termes:
\begin{equation}
D=H\cap\partial\kappa_x\cap\partial O\cap\partial\kappa_y\neq\emptyset.
\end{equation}
Notons $d$ un élément de $D$, et démontrons à présent qu'il est possible de trouver une face $F$ de $S$ telle que $d\in F$ et $F\subset\overline{\kappa_x}\cap\overline{\kappa_y}$.

Soient $S_x$ et $S_y$ les sous-complexes de $S$ formés des polyèdres respectivement inclus dans $\overline{\kappa_x}$ et $\overline{\kappa_y}$:
\begin{align}
S_x&=\left\{\delta\in S\colon\delta\subset\overline{\kappa_x}\right\}&S_y&=\left\{\delta\in S\colon\delta\subset\overline{\kappa_y}\right\}.
\end{align}
Puisque par hypothèse $\partial O=\U(S)$ il vient
\begin{equation}
\partial O\cap\overline{\kappa_x}\cap\overline{\kappa_y}=\bigcup_{\delta\in S_x,\delta'\in S_y}\delta\cap\delta'.
\end{equation}
En outre, $S_x\cap S_y=\emptyset$ donc par définition d'un complexe, $\delta\cap\delta'$ est soit vide, soit une sous-face commune de $\delta$ et $\delta'$. On remarque aussi que $S_x$ et $S_y$ ne sont pas vides, car ils contiennent chacun au moins un polyèdre qui contient $d$ (puisque $d\in\partial O\cap\overline{\kappa_x}\cap\overline{\kappa_y}$).

Considérons par ailleurs pour $\epsilon>0$, la boule ouverte $B_\epsilon$ de centre $d$ et de rayon $\epsilon$, et notons
\begin{equation}
U_\epsilon=B_\epsilon\cap\partial O.
\end{equation}
Par hypothèse $\overline{O}=\bigcup_t\overline{\kappa_t}$ donc $\partial O\subset\bigcup_t\overline{\kappa_t}$, en outre $d$ n'appartient à aucun des fermés $\overline{\kappa_t}$ pour $t\notin\{x,y\}$, donc il existe $\epsilon_0>0$ tel que pour tout $\epsilon<\epsilon_0$ 
\begin{equation}
U_\epsilon\subset\overline{\kappa_x}\cup\overline{\kappa_y}.
\end{equation}
À présent considérons les deux cas possibles:
\begin{itemize}
\item si $d$ n'est sommet d'aucun polyèdre de $S$ alors il existe $\epsilon_1>0$ tel que $U_\epsilon$ ne contient aucun sommet des polyèdres de $S$ pour $\epsilon\leq\epsilon_1$. Dans ce cas, puisque $d$ ne peut être que sur la frontière de tout polyèdre qui le contient (car il est dans $\overline{\kappa_x}\cap\overline{\kappa_y}$), $U_\epsilon$ est donc l'intersection de $B_\epsilon$ et de l'union d'une famille finie $A_1,\ldots,A_m$ de demi-hyperplans affines contenant $d$
\begin{equation}
U_\epsilon=B_\epsilon\cap\bigcup_iA_i.
\end{equation}
Chacun des $A_i\cap B_\epsilon$ est un morceau d'un polyèdre $n-1$-dimensionnel $\delta_i$ inclus dans $\overline{\kappa_x}$ ou $\overline{\kappa_y}$, et il est clair que pour tout $i$, $\partial\kappa_x\cap\partial\kappa_y\cap\partial O$ contient $\partial A_i\cap B_\epsilon$ qui contient lui-même $d$, $\partial A_i$ désignant ici la frontière du demi-hyperplan $A_i$ prise relativement à l'hyperplan entier. Si on considère par exemple le sous-espace affine $\partial A_1$ (de dimension $n-2$), $F=\delta_1\cap\partial A_1$ est une face de $\delta_1$, et puisque $S$ est un complexe alors $F$ est dans $\partial\kappa_x\cap\partial\kappa_y\cap\partial O$ (puisqu'on vient de voir qu'une boule $n-2$ dimensionnelle incluse dans $F$ y est déjà, si la face entière n'y était pas cela impliquerait qu'un sommet d'un autre polyèdre de $S$ serait dans $F$ privé de ses propres sommets, ce qui contredirait le fait que $S$ est un complexe);
\item si $d$ est sommet d'un polyèdre de $S$ alors il existe $\epsilon_2>0$ tel que $U_\epsilon$ ne contient aucun autre sommet de $S$ pour $\epsilon\leq\epsilon_2$. Par un raisonnement analogue sur un point $e\in\partial\kappa_x\cap\partial\kappa_y\cap\partial O\cap B_\epsilon$ distinct de $d$, il est possible de trouver une face incluse dans $\partial\kappa_x\cap\partial\kappa_y\cap\partial O$ qui contient $d$ et $e$.
\end{itemize}

Dans les deux cas, on a donc trouvé une face $F$ telle que
\begin{equation}
d\in F\subset\overline{\kappa_x}\cap\overline{\kappa_y}\cap\partial O.
\end{equation}
Il nous reste encore à vérifier que le polyèdre $\left<F\cup\{x,y\}\right>$ est de dimension $n$. Puisque $\kappa_x$ et $\kappa_y$ sont respectivement étoilés par rapport à $x$ et $y$, et $F\subset\overline{\kappa_x}\cap\overline{\kappa_y}$ alors le sous-espace affine minimal $H'$ de dimension $n-2$ contenant $F$ ne contient pas $x$ et $y$, et sa direction n'est pas parallèle à celle de la droite $(x,y)$, donc $\left<F,\{x\}\right>$ est de dimension $n-1$ et l'hyperplan affine le contenant ne contient pas $y$. Dès lors
\begin{equation}
\left<(\delta_1\cap\delta_2)\cup\{x,y\}\right>\in\S(S,G,p)
\end{equation}
avec par construction, $a\in\left<\{x,y,d\}\right>\subset\left<\{x,y\}\cup (\delta_1\cap\delta_2)\right>$, c'est à dire que $a\in\U(S')$;
\item si $b\in\partial\kappa_z\cap\overline{\kappa_x}$ alors il suffit de refaire le raisonnement précédant en permutant $y$ et $z$.
\end{itemize}

Si $\O(x)=1$ alors en notant $y$ le voisin de $x$ il est possible de refaire un raisonnement semblable en remplaçant $\kappa_z$ par $\emptyset$.

Pour conclure, dans tous les cas on a bien démontré que $a\in\U(S')$, c'est à dire que $\overline{O}\subset\U(S')$. L'inclusion inverse est clairement vérifiée dès que les $(\kappa_x)_{x\in T}$ vérifient les hypothèses du lemme~\ref{lemmaocclusionsimple}, puisque dans ce cas tous les nouveaux polyèdres créés lors de la suspension restent inclus dans $\bigcup_{x\in T}\overline{\kappa_x}$:
\begin{itemize}
\item les polyèdres de première espèce obtenus par suspension par rapport à un sommet $x$ sont inclus dans $\overline{\kappa_x}$;
\item ceux de seconde espèce obtenus par suspension par rapport à deux sommets voisins $x$ et $y$ sont inclus dans $\overline{\kappa_x}\cup\overline{\kappa_y}$.
\end{itemize}
\end{proof}

À présent donnons une définition pour décrire la situation où un complexe combine les hypothèses des lemmes~\ref{lemmaocclusionsimple} et~\ref{lemmaocclusiontotale} par rapport à un ouvert en forme de <<~tube~>> parcouru par un graphe linéaire, en utilisant des parties $(\kappa_x)_{x\in T}$ obtenues par découpage du tube par des hyperplans perpendiculaires aux arêtes du graphe et pouvant être placés à $\epsilon$ près.

Supposons que $n\geq2$, considérons un complexe $n-1$-dimensionnel $S$, un ouvert $O$ tel que $\partial O=\U(S)$ et un graphe linéaire $G=(T,A)$ tel que $T\subset O$. Pour deux sommets voisins $x$ et $y$ de $G$ et $r\in\left]-\frac{\dist(x,y)}{2},\frac{\dist(x,y)}{2}\right[$ on notera $H(x,y,r)$ l'hyperplan affine perpendiculaire à et passant par le segment $[x,y]$ à distance $r+\frac{\dist(x,y)}{2}$ de $x$, $H^+(x,y,r)$ le demi-espace affine de frontière $H(x,y,r)$ qui contient $x$ et $H^-(x,y,r)$ l'autre demi-espace correspondant. Soit $f$ une application de $T^2$ dans $\mathbb{R}$, anticommutative (c'est à dire telle que $\forall (x,y)\in T^2\colon f(x,y)=-f(y,x)$) et vérifiant
\begin{equation}
\forall \{x,y\}\in A\colon\vert f(x,y)\vert<\frac{\dist(x,y)}{2}.
\end{equation}
Supposons que $\#T>1$, pour un sommet $x\in T$ selon les cas:
\begin{itemize}
\item si $\O(x)=0$ (dans ce cas $G$ est réduit à un seul sommet par hypothèse de linéarité) on posera $\kappa_x(f)=O$;
\item si $\O(x)=1$ alors soit $y$ le voisin de $x$, on notera $\kappa_x(f)$ l'adhérence de la composante connexe de $H^+(x,y,f(x,y))\cap O$ qui contient $x$;
\item si $\O(x)=2$ alors soient $y$ et $z$ les deux voisins distincts de $x$, on notera $\kappa_x(f)$ l'adhérence de la composante connexe de $H^+(x,y,f(x,y))\cap H^+(x,z,f(x,z))\cap O$ qui contient $x$.
\end{itemize}
Et pour finir on définit encore la famille de polyèdres $n-1$-dimensionnels $S(f)$ par
\begin{equation}
\delta\in S(f)\Longleftrightarrow\exists(\delta',x)\in S\times T\colon\delta=\delta'\cap\kappa_x(f)\text{ et }\dim\delta=n-1.
\end{equation}

\begin{definition}[Suspension tubulaire]\label{definitionsuspensiontubulaire}
Soient $S$ un complexe $n-1$-dimen\-sion\-nel, $O$ un ouvert tel que $\partial O=\U(S)$ et $G=(T,A)$ un graphe linéaire tel que $T\subset O$. Pour $\epsilon>0$ donné, on dira que $S$ est $\epsilon$-tubulaire par rapport à $G$ et $O$ si $G$ est réduit à un seul sommet, ou s'il existe deux applications $f$ et $g\colon T^2\rightarrow\mathbb{R}$ anticommutatives et vérifiant:
\begin{enumerate}
\item $\forall \{x,y\}\in A$:
\begin{equation}
\begin{aligned}
\left(f(x,y),g(x,y)\right)&\in\left]-\frac{\dist(x,y)}{2},\frac{\dist(x,y)}{2}\right[^2\\
\left\vert f(x,y)-g(x,y)\right\vert&\geq\epsilon;
\end{aligned}
\end{equation}
\item pour toute application anticommutative $h\colon T^2\rightarrow\mathbb{R}$ telle que $\min(f,g)\leq h\leq\max(f,g)$, $S(h)$ est un complexe vérifiant les hypothèses des lemmes~\ref{lemmaocclusionsimple} et~\ref{lemmaocclusiontotale} par rapport au graphe $G$ et $O$ en utilisant les ouverts $\kappa_x(h)$ et le choix $p$ de suspension défini par $p(\delta)=x\Leftrightarrow\delta\subset\kappa_x(h)$.
\end{enumerate}
\end{definition}

Lorsque les hypothèses de la définition sont vérifiées, $\partial O$ admet un hyperplan tangent $\H^{n-1}$-presque partout (en tant qu'union finie de polyèdres $n-1$-dimensionnels). On notera $\partial^\bot O$ le sous-ensemble de $\partial O$ où un tel hyperplan existe, et $\overrightarrow{n}(z)$ un vecteur unitaire normal à cet hyperplan lorsque $z\in\partial^\bot O$. Lorsque le graphe $G$ est réduit à un seul sommet $x$ on notera
\begin{equation}
\begin{aligned}
\alpha_-=\beta_-&=\dist(x,\partial O)&\gamma&=\left(\H^{n-1}(\partial O)\right)^\frac{1}{n-1}\\
\alpha_+=\beta_+&=\sup_{t\in\partial O}\dist(x,t)&\eta&=\inf_{t\in\partial^\bot O}\frac{\vert<\overrightarrow{n}(t),t-x>\vert}{\Vert t-x\Vert}.
\end{aligned}
\end{equation}
Lorsque $\#T>1$ on notera $\mathfrak{F}$ l'ensemble des applications anticommutatives comprises entre $f$ et $g$:
\begin{equation}
\mathfrak{F}=\left\{h\in\mathbb{R}^{T^2}\colon h\text{ anticommutative et }\min(f,g)\leq h\leq\max(f,g)\right\}
\end{equation}
et avec ces notations on définit encore les quantités suivantes:
\begin{equation}\label{equationsuspensiontubulaireA}
\begin{aligned}
\alpha_-&=\min_{\{x,y\}\in A}\dist(x,y)&\beta_-&=\inf_{\substack{h\in\mathfrak{F}\\\{x,y\}\in A\\z\in[x,y]}}\dist(z,\partial O\cap(\overline{\kappa_x(h)}\cup\overline{\kappa_y(h)}))\\
\alpha_+&=\max_{\{x,y\}\in A}\dist(x,y)&\beta_+&=\sup_{\substack{h\in\mathfrak{F}\\\{x,y\}\in A\\z\in[x,y]}}\sup_{t\in\partial O\cap(\overline{\kappa_x(h)}\cup\overline{\kappa_y(h)})}\dist(z,t)\\
&&\gamma&=\left(\sup_{\substack{h\in\mathfrak{F}\\x\in T}}\H^{n-1}(\partial O\cap\overline{\kappa_x(h)})\right)^\frac{1}{n-1}\\
&&\eta&=\inf_{\substack{h\in\mathfrak{F}\\\{x,y\}\in A\\z\in[x,y]\\t\in\partial^\bot O\cap(\overline{\kappa_x(h)}\cup\overline{\kappa_y(h)})}}\frac{\vert<\overrightarrow{n}(t),t-z>\vert}{\Vert t-z\Vert}.
\end{aligned}
\end{equation}

Le lemme suivant permet d'évaluer les constantes de forme optimales des polyèdres d'une suspension de complexe $\epsilon$-tubulaire. Comme dans le cas de la propriété~\ref{propertysuspension}, il serait possible de le donner sous une forme plus précise mais celle-ci nous suffira pour la suite.

\begin{lemma}[Rotondité d'une suspension tubulaire]\label{lemmasuspensiontubulaire}
Pour toute partie compacte $K\subset]0,+\infty[^9$ il existe des constantes $c_1$ et $c_2$ strictement positives telles que pour tout complexe $\epsilon$-tubulaire $S$ par rapport à un ouvert $O$ et un graphe $G$, si on peut trouver deux constantes $\rho_+$ et $\rho_-$ telles que
\begin{equation}\label{equationsuspensiontubulaireB}
\rho_+>\overline{\R}(S')>\underline{\R}(S')>\rho_->0
\end{equation}
et
\begin{equation}\label{equationsuspensiontubulaireC}
\left(n,\eta,\frac{\rho_+}{\rho_-},\frac{\epsilon}{\rho_+},\frac{\alpha_-}{\rho_+},\frac{\alpha_+}{\rho_+},\frac{\beta_-}{\rho_+},\frac{\beta_+}{\rho_+},\frac{\gamma}{\rho_-}\right)\in K
\end{equation}
alors on peut trouver $h\in\mathfrak{F}$ tel que $S'=\S(S(h),G,p)$ est un complexe $n$-dimensionnel vérifiant
\begin{equation}
\overline{\R}(S')\leq c_1\overline{\R}(S)\text{, }\underline{\R}(S')\geq c_2\underline{\R}(S)\text{ et }\U(S')=\overline{O}.
\end{equation}
\end{lemma}

\begin{proof}
On peut commencer par remarquer que le lemme est évident si $\#T=1$. En effet, dans ce cas il suffit de considérer $\delta\in\F(S)$, en notant $o$ un orthocentre de $\delta$ et $x\in T$ il vient
\begin{equation}
\begin{gathered}
\eta\beta_-\leq\dist(x,\affine(\delta))\leq\beta_+\\
\beta_-\leq\dist(x,o)\leq\beta_+\\
\rho_-\leq\overline{R}(\delta)\leq\rho_+
\end{gathered}
\end{equation}
d'où on tire
\begin{equation}
\begin{gathered}
\eta\frac{\beta_-}{\beta_+}\leq\frac{\dist(x,\affine(\delta))}{\dist(x,o)}\leq1\\
\frac{\beta_-}{\rho_+}\leq\frac{\dist(x,o)}{\overline{R}(\delta)}\leq\frac{\beta_+}{\rho_-}.
\end{gathered}
\end{equation}
En appliquant la propriété~\ref{propertysuspension} on obtient immédiatement les constantes désirées, on supposera donc pour la suite que $\#T>1$.

Soit $\delta\in\F_k(S)$ une sous-face de dimension $k\leq n-1$ d'un polyèdre de $S$, considérons les deux cas:
\begin{itemize}
\item si $k=1$ alors $\#\F_{k-1}(\delta)=2$;
\item si $k>1$, par hypothèse $\delta$ est contenu dans une boule de centre $c_\delta$ et de rayon au plus $\rho_+$. Notons respectivement $U_k$ et $V_k$ la surface de la sphère unité et le volume de la boule unité en dimension $k$. En rappelant que $\delta$ est convexe et contenu dans une boule de rayon $\rho_+$ il vient
\begin{equation}\label{equationsuspensiontubulaireD}
\H^{k-1}(\partial\delta)\leq U_k\rho_+^{k-1}
\end{equation}
et puisque les faces de $\delta$ sont disjointes et forment sa frontière on a aussi
\begin{equation}
\sum_{\alpha\in\F_{k-1}(\delta)}\H^{k-1}(\alpha)=\H^{k-1}(\partial\delta).
\end{equation}
De la même façon, les faces de $\delta$ contiennent une boule $k-1$-dimensionnelle de rayon au moins $\rho_-$ et donc
\begin{equation}\label{equationsuspensiontubulaireE}
\#\F_{k-1}(\delta)V_{k-1}\rho_-^{k-1}\leq\H^{k-1}(\partial\delta).
\end{equation}
On peut alors tirer de~\eqref{equationsuspensiontubulaireD} et~\eqref{equationsuspensiontubulaireE}:
\begin{equation}
\#\F_{k-1}(\delta)\leq\frac{U_k}{V_{k-1}}\left(\frac{\rho_+}{\rho_-}\right)^{k-1}.
\end{equation}
\end{itemize}
Par conséquent, on peut donc trouver $M>0$ ne dépendant que de $K$ tel que pour $0\leq k\leq n-1$ et pour tout sous-complexe $S'\subset S$
\begin{equation}\label{equationtubulaireA}
\#\F_k(S')\leq M\#\F_{n-1}(S').
\end{equation}

Puisqu'on a supposé que $G$ n'est pas réduit à un seul sommet alors $A\neq\emptyset$, on peut donc considérer $\{x,y\}\in A$ et poser
\begin{equation}
u=\min(f(x,y),g(x,y))\text{ et }v=\max(f(x,y),g(x,y)).
\end{equation}
Fixons $k\in\{0,\ldots,n-1\}$ et considérons le sous-complexe $S_k(x,y)$ des sous-faces de dimension $k$ de $S$ pour lesquelles il existe $w\in]u,v[$ tel que l'hyperplan affine $H(x,y,w)$ intersecte leur intérieur sans les contenir:
\begin{equation}
S_k(x,y)=\{\delta\in\F_k(S)\colon\exists w\in]u,v[,\open{\delta}\cap H(x,y,w)\neq\emptyset\text{ et }\delta\not\subset H(x,y,w)\}.
\end{equation}

Par ailleurs en prenant $k=n-1$ on sait que
\begin{equation}
\U(S_{n-1}(x,y))\subset\partial O\cap\bigcup_{u\leq w\leq v}H(x,y,w)
\end{equation}
sinon il y aurait des polyèdres de $S$ qui pourraient être découpés plusieurs fois par les hyperplans $H(x,y,w)$ pour $\{x,y\}\in A$, ce que les hypothèses du lemme interdisent. On en déduit que les $S_{n-1}(x,y)$ (pour $\{x,y\}\in A$) sont des sous-complexes disjoints deux à deux de $S$ qui vérifient
\begin{equation}
\U(S_{n-1}(x,y))\subset(\overline{\kappa_x(f)}\cup\overline{\kappa_x(g)})\cap\partial O
\end{equation}
et on en tire
\begin{equation}
\H^{n-1}(\U(S_{n-1}(x,y)))\leq\H^{n-1}((\kappa_x(f)\cup\kappa_x(g))\cap\partial O)\leq2\gamma^{n-1}.
\end{equation}
En outre, tous les polyèdres de $S_{n-1}(x,y)$ contiennent une boule $n-1$-dimen\-sionnelle de rayon au moins $\underline{\R}(S)>\rho_-$ donc
\begin{equation}
V_{n-1}\#S_{n-1}(x,y)\rho_-^{n-1}\leq\H^{n-1}(\U(S_{n-1}(x,y)))
\end{equation}
d'où on tire
\begin{equation}
\#S_{n-1}(x,y)\leq\frac{2\gamma^{n-1}}{\rho_-^{n-1}V_{n-1}}
\end{equation}
et d'après~\eqref{equationtubulaireA}
\begin{equation}\label{equationtubulaireB}
\sum_{k<n}\#S_k(x,y)\leq\frac{2nM}{V_{n-1}}\left(\frac{\gamma}{\rho_-}\right)^{n-1}\leq N
\end{equation}
avec $N$ entier qui là encore ne dépend que du compact $K$.

À présent considérons un polyèdre $\delta\in S_k(x,y)$ et $B$ une boule $k-1$-dimensionnelle de $\affine(\delta)$ inscrite dans $\delta$ centrée sur un orthocentre $c$ et de rayon $\underline{R}(\delta)$ (par hypothèse, $\underline{R}(\delta)>\rho_-$). Notons encore:
\begin{flalign}
u'&=\inf\{w\in[u,v]\colon H(x,y,w)\cap\open{\delta}\neq\emptyset\}&\delta^+(w)&=H^+(x,y,w)\cap\delta\\
v'&=\sup\{w\in[u,v]\colon H(x,y,w)\cap\open{\delta}\neq\emptyset\}&\delta^-(w)&=H^-(x,y,w)\cap\delta\\
S(\delta,w)&=
\begin{cases}
\{\delta^+(w),\delta^-(w)\}&\text{si $w\in]u',v'[$}\\
\{\delta\}&\text{si $w\notin]u',v'[$.}
\end{cases}
\end{flalign}
Par construction $H(x,y,u')\cap\open{B}=H(x,y,v')\cap\open{B}=\emptyset$ (car $B\subset\delta$), $\delta^+(w)$ et $\delta^-(w)$ sont deux polyèdres $k+1$-dimensionnels pour $u'<w<v'$, $S(\delta,w)$ est un complexe $k$-dimensionnel et il existe $w_0\in]u',v'[$ tel que $c\in H(x,y,w_0)$. Soient $x'\in\delta\cap H(x,y,u')$ et $y'\in\delta\cap H(x,y,v')$ et posons $B'=\S(B,x')\cup\S(B,y')$. Déjà $B'\subset\delta$ par convexité, en outre $B\cap  H^+(x,y,w_0)$ et $B\cap  H^-(x,y,w_0)$ contiennent chacun une boule de rayon $\underline{R}(\delta)/2$. En posant
\begin{align}
B^+(w)=&B'\cap H^+(x,y,w)&\text{et}&&B^-(w)=&B'\cap H^-(x,y,w)
\end{align}
et en considérant une homothétie de centre $x'$ il vient
\begin{equation}
\forall w\in]u',w_0]\colon\underline{R}(B^+(w))\geq\frac{w-u'}{w_0-u'}\cdot\frac{\underline{R}(\delta)}{2}
\end{equation}
et de manière symétrique
\begin{equation}
\forall w\in[w_0,v'[\colon\underline{R}(B^-(w))\geq\frac{w-v'}{w_0-v'}\cdot\frac{\underline{R}(\delta)}{2}.
\end{equation}
Réciproquement,
\begin{equation}
\forall w\in]w_0,v'[\colon\underline{R}(B^+(w))\geq\frac{\underline{R}(\delta)}{2}
\end{equation}
et
\begin{equation}
\forall w\in]u',w_0[\colon\underline{R}(B^-(w))\geq\frac{\underline{R}(\delta)}{2}.
\end{equation}
Par ailleurs puisque par hypothèse $\underline{R}(\delta)\geq\rho_-$ et $\overline{R}(\delta)\leq\rho_+$, alors
\begin{equation}
\min(w_0-u',v'-w_0)\geq\frac{\rho_-}{2\rho_+}(v'-u')
\end{equation}
et en posant
\begin{equation}
\psi_\delta(w)=
\begin{cases}
1&\text{si $w\leq u'$}\\
\fracd{\min(\vert w-u'\vert,\vert w-v'\vert)}{v'-u'}\cdot\fracd{\rho_-}{4\rho_+}&\text{si $u'<w<v'$}\\
1&\text{si $w\geq u'$}
\end{cases}
\end{equation}
il vient
\begin{equation}
\min_{\alpha\in S(\delta,w)}\underline{R}(\alpha)\geq\underline{R}(\delta).
\end{equation}

Avec nos notations, on peut donc écrire
\begin{equation}\label{equationtubulaireC}
\underline{\R}(S(x,y,w))\geq\min_{\substack{1\leq k\leq n-1\\\delta\in S_k(x,y)}}\psi_\delta(w)\underline{\R}(S)
\end{equation}
où $S(x,y,w)$ désigne le complexe obtenu par découpage des polyèdres de $S_{n-1}(x,y)$ par l'hyperplan $H(x,y,w)$. Pour $a\in\mathbb{R}$ et $b\in]0,+\infty]$ écrivons encore
\begin{equation}
\phi_{a,b}(w)=
\begin{cases}
1&\text{si $w\leq a$}\\
\fracd{w-a}{b/2}&\text{si $w\in]a,a+b/2]$}\\
\fracd{a+b-w}{b/2}&\text{si $w\in]a+b/2,a+b[$}\\
1&\text{si $w\geq a+b$}
\end{cases}
\end{equation}
et remarquons que pour $\delta\in S_k(x,y)$ on a
\begin{equation}
\psi_\delta(w)\geq\frac{\rho_-}{4\rho_+}\phi_{u',v'-u'}(w)
\end{equation}
avec $v'-u'\leq2\rho_+$. Par ailleurs il est facile de constater que 
\begin{equation}\label{equationtubulaireD}
\forall(a,b,\sigma)\in\mathbb{R}\times]0,2\rho_+]\times[0,1[\colon\H^1(\left\{w\in\mathbb{R}\colon\phi_{a,b}(w)\leq\sigma\right\})\leq b\sigma\leq2\rho_+\sigma.
\end{equation}
D'après~\eqref{equationtubulaireB} il vient
\begin{equation}
\sup_{w\in[u,v]}\min_{\substack{1\leq k\leq n-1\\\delta\in S_k(x,y)}}\psi_\delta(w)\underline{\R}(S)\geq\frac{\rho_-}{4\rho_+}\min_{\substack{(a_1,\ldots,a_N)\in\mathbb{R}^N\\(b_1,\ldots,b_N)\in]0,2\rho_+]^N}}\sup_{w\in[0,\epsilon]}\phi_{a,b}(w)
\end{equation}
et en utilisant~\eqref{equationtubulaireD}
\begin{multline}\label{equationtubulaireE}
\forall\sigma\in[0,1[,\forall(a_1,\ldots,a_N)\in\mathbb{R}^N,\forall(b_1,\ldots,b_N)\in]0,2\rho_+]^N\colon\\\H^1\left(\left\{w\in\mathbb{R}\colon\phi_{a,b}\leq\sigma\right\}\right)\leq2N\rho_+\sigma.
\end{multline}
En prenant $\sigma=\frac{\epsilon}{4N\rho_+}$ dans~\eqref{equationtubulaireE} et d'après~\eqref{equationtubulaireC} on peut donc trouver $w_{x,y}\in[u,v]$ tel que
\begin{equation}
\underline{\R}(S(x,y,w_{x,y}))\geq\frac{\epsilon\rho_-}{16N\rho_+^2}\underline{\R}(S)=A\underline{\R}(S)
\end{equation}
avec là encore $A>0$ qui ne dépend que de $K$. En posant $h(x,y)=w_{x,y}$ et en recommençant avec tous les couples de voisins du graphe, il est donc possible de trouver $h\in\mathfrak{F}$ telle que:
\begin{equation}\label{equationtubulaireF}
\underline{\R}(S(h))\geq A\underline{\R}(S)\text{ et }\overline{\R}(S(h))\leq\overline{\R}(S).
\end{equation}

Il nous reste encore à établir des relations uniformes sur les régularités après suspension. Pour cela soit $\delta\in\F(S')$ une sous-face de $S'$ et rappelons que comme vu dans la démonstration du lemme~\ref{lemmaocclusionsimple}, $\delta$ peut être de quatre sortes:
\begin{itemize}
\item si $\delta\in\F(S(h))$ alors par définition
\begin{equation}\label{equationtubulaireG}
\overline{R}(\delta)\leq\overline{\R}(S(h))\text{ et }\underline{R}(\delta)\geq\underline{\R}(S(h));
\end{equation}
\item si $\delta=[u,v]$ où $u$ et $v$ sont deux sommets voisins du graphe, alors 
\begin{equation}\label{equationtubulaireH}
\overline{R}(\delta)=\underline{R}(\delta)\in\left[\alpha_-,\alpha_+\right]\subset\left[\frac{\alpha_-}{\rho_+}\underline{\R}(S),\frac{\alpha_-}{\rho_-}\overline{\R}(S)\right];
\end{equation}
\item si $\delta=\S(F,u)$ où $F\in\F(S(h))$ est une sous-face de $S(h)$ et $u\in T$ un sommet du graphe, notons $c$ un orthocentre de $F$. On a
\begin{equation}
\frac{\dist(u,\affine(\delta))}{\dist(u,c)}\in\left[\eta\frac{\beta_-}{\beta_+},1\right]\text{ et }\frac{\dist(u,c)}{\overline{R}(F)}\in\left[\frac{\beta_-}{\rho_+},\frac{\beta_+}{\rho_-}\right]
\end{equation}
donc la propriété~\ref{propertysuspension} nous donne des constantes $a$ et $b$ ne dépendant que de $K$ telles que
\begin{equation}\label{equationtubulaireI}
\overline{R}(\delta)\leq a\overline{R}(F)\leq a\overline{\R}(S(h))\text{ et }\underline{R}(\delta)\geq b\underline{R}(F)\geq b\underline{\R}(S(h));
\end{equation}
\item si $\delta=\left<F\cup\{u,v\}\right>=\S(\S(F,u),v)$ où $F$ est une sous-face de $S(h)$ et $\{u,v\}\in A$ deux sommets voisins du graphe, notons $F'=\S(F,u)$ et $c'$ un orthocentre de $F'$. Déjà, comme vu pour le cas précédent on a
\begin{equation}\label{equationtubulaireJ}
\overline{R}(F')\leq a\overline{\R}(S(h))\text{ et }\underline{R}(F')\geq b\underline{\R}(S(h))
\end{equation}
d'où on tire d'après~\eqref{equationtubulaireF}
\begin{equation}\label{equationtubulaireK}
\overline{R}(F')\geq\underline{R}(F')\geq b\underline{\R}(S(h))\geq bA\rho_-\text{ et }\overline{R}(F')\leq a\rho_+.
\end{equation}
Par définition de la suspension tubulaire $F\subset H(u,v,h(u,v))$, $F'\subset\kappa_u$, $\open{F'}\cap\kappa_v=\emptyset$ et $\dim\affine(F)\leq n-2$, donc 
\begin{equation}
\left(\dist(v,\affine(F')),\dist(v,c')\right)\in\left[\dist(v,F'),\beta_+\right]^2.
\end{equation}
Par ailleurs 
\begin{equation}
\dist(v,F')\in\left[\frac{\min(\alpha_-,\beta_-)}{2},\max(\alpha_+,\beta_+)\right]
\end{equation}
et donc d'après~\eqref{equationtubulaireK}
\begin{equation}
\begin{split}
\frac{\dist(v,\affine(F'))}{\dist(v,c')}&\in\left[\frac{\min(\alpha_-,\beta_-)}{2\max(\alpha_+,\beta_+)},1\right]\\
\frac{\dist(v,c')}{\overline{R}(F')}&\in\left[\frac{\min(\alpha_-,\beta_-)}{2a\rho_+},\frac{\max(\alpha_+,\beta_+)}{bA\rho_-}\right]
\end{split}
\end{equation}
d'où on tire de nouveau d'après la propriété~\ref{propertysuspension} et l'inégalité~\eqref{equationtubulaireJ}
\begin{equation}\label{equationtubulaireL}
\overline{R}(\delta)\leq a'\overline{R}(F')\leq aa'\overline{\R}(S(h))\text{ et }\underline{R}(\delta)\geq b'\underline{R}(F')\geq bb'\underline{\R}(S(h))
\end{equation}
avec là encore $a'$ et $b'$ qui ne dépendent que de $K$.
\end{itemize}

En combinant l'inégalité~\eqref{equationtubulaireF} avec celles obtenues pour chacun des quatre cas (\eqref{equationtubulaireG},~\eqref{equationtubulaireH},~\eqref{equationtubulaireI} et~\eqref{equationtubulaireL}) on obtient les constantes $c_1$ et $c_2$ désirées, ce qui termine la démonstration du lemme.
\end{proof}

\section{Raccordement de complexes dyadiques}\label{sectionC}

On s'intéresse à présent à des complexes composés de cubes dyadiques de $\mathbb{R}^n$. Un cube dyadique de pas $r>0$ s'écrit comme le pavé $[0,r]^n$ dans une base orthonormale adaptée. Un famille de tels cubes disposés sur un pavage constitue naturellement un complexe.

\begin{definition}[Complexes dyadiques]\label{definitioncomplexedyadique}
On appellera complexe dyadique $n$-dimensionnel de pas $r>0$ toute famille $S$ de cubes dyadiques qui peut s'écrire
\begin{equation}
S=\left\{rz+[0,r]^n\colon z\in Z\right\}
\end{equation}
dans une base orthonormale adaptée.

On dira qu'un complexe $T$ est dyadique en surface s'il existe un complexe dyadique $S$ tel que
\begin{equation}
\F_\partial(T)=\F_\partial(S).
\end{equation}
\end{definition}

Comme précédemment on généralise la définition à des complexes dyadiques de dimension $k\leq n$ en se plaçant dans un sous-espace affine. Lorsqu'un complexe $T$ est dyadique en surface il est clair qu'il existe un seul complexe dyadique $S$ avec les mêmes faces sur la frontière. On généralise donc naturellement à $T$ les attributs de $S$ propres aux complexes dyadiques (par exemple les propriétés de groupement de la définition~\ref{definitiongroupement}). 

Les principales définitions ayant été données, on peut à présent énoncer le théorème principal.

\begin{theorem}[Fusion]\label{theoremfusion}
Il existe trois constantes strictement positives $\rho$, $c_1$ et $c_2$ ne dépendant que de $n$ telles que pour tout compact $K\subset\mathbb{R}^n$, pour tout ouvert $O\subset K$ et pour tous complexes $n$-dimensionnels $S_1$ et $S_2$ dyadiques unitaires en surface vérifiant
\begin{align}
\U(S_1)&=K\setminus O&\U(S_2)&\subset O&\min_{(x,y)\in\U(S_1)\times\U(S_2)}\Vert x-y\Vert&>\rho
\end{align}
on peut construire $S_3$ tel que $S'=S_1\cup S_2\cup S_3$ est un complexe $n$-dimensionnel vérifiant
\begin{align}
\U(S')&=K&\overline{\R}(S')&\leq c_1\overline{\R}(S_1\cup S_2)&\underline{\R}(S')&\geq c_2\underline{\R}(S_1\cup S_2).
\end{align}
\end{theorem}

Dans ce qui suit on notera $\phi$ une isométrie affine qui fait passer d'une base de $S_1$ à une base de $S_2$. La démonstration va consister à combler l'espace compris entre $S_1$ et $S_2$ par des couches successives de polyèdres de façon à former un seul complexe plus grand, en raisonnant par récurrence sur la dimension $n$. Auparavant, on va donner un lemme préliminaire.

\subsection{Subdivision}

Définissons une propriété caractérisant les complexes dyadiques qui sont des subdivisions de complexes dyadiques plus grands.

\begin{definition}[Groupement]\label{definitiongroupement}
Soit $S$ un complexe dyadique de pas $r$ et $p$ un entier non nul. On dira que $S$ est $\underbrace{p\times p\cdots\times p}_{k\text{ fois}}$-groupé s'il existe un complexe dyadique $T$ de pas $pr$ tel que $\U(S)=\U(T)$.

Dans ce cas on dira que $S$ est la $p\times p\cdots\times p$-subdivision de $T$, ou que $T$ est le $p\times p\cdots\times p$-groupement de $S$.
\end{definition}

On se réserve aussi le droit de parler de $p_1\times p_2\ldots\times p_k$-groupements par la suite (où les $p_i$ sont des entiers non nuls pas forcément tous égaux) lorsque le complexe $S$ s'y prête. Pour simplifier les notations, on remplacera aussi parfois $p\times p\times\ldots\times p$ par $p^k$.

\'Etant donnés un entier $m>1$ et un complexe dyadique en surface, le lemme suivant permet de construire une couche de polyèdres sur la frontière de façon à le rendre $m^k$-groupé en surface, tout en gardant un contrôle sur les constantes de forme obtenues.

\begin{lemma}[Subdivision]\label{lemmasubdivision}
Pour tout entier $p\geq1$ il existe deux constantes $c_1$ et $c_2$ telles que pour tout complexe $n$-dimensionnel $S$ dyadique en surface de pas $r$ il est possible de construire un complexe $S'\supset S$ dyadique en surface de pas $r/p$, $p^n$-groupé dans la même base et vérifiant
\begin{gather}
\max_{x\in\U(S')}\dist(x,\U(S))\leq r\sqrt{n}\label{equationsubdivisionA}\\
\begin{aligned}
\overline{\R}(S')&\leq c_2\overline{\R}(S)&\underline{\R}(S')&\geq c_1\underline{\R}(S)\label{equationsubdivisionB}.
\end{aligned}
\end{gather}
\end{lemma}

\begin{proof}
Considérons la famille $T$ des cubes dyadiques de pas $r$ dans une base adaptée à $S$ qui ont au moins un point commun avec $\partial\U(S)$ et qui sont d'intérieur disjoint avec tous les polyèdres de $S$. Il est clair que $S\cup T$ est un complexe dyadique en surface de pas $r$ qui vérifie l'inégalité~\eqref{equationsubdivisionA}, avec en outre
\begin{equation}
\partial\U(S\cup T)\cap\partial\U(S)=\emptyset.
\end{equation}

Soit $\delta\in T$ et notons $c_\delta$ son centre. L'ensemble $\F_{n-1}(\delta)\cap\F_\partial(T)$ contient des cubes dyadiques $n-1$-dimensionnels de pas $r$, qu'il est possible de subdiviser chacun naturellement en $p^{n-1}$ cubes de même dimension et de pas $r/p$. On notera $E_\delta$ le complexe obtenu après subdivision des cubes de $\F_{n-1}(\delta)\cap\F_\partial(T)$, et
\begin{equation}
F_\delta=\F_{n-1}(\delta)\setminus\F_\partial(T).
\end{equation}

Par construction, $E_\delta\cup F_\delta$ est un complexe $n-1$-dimensionnel vérifiant
\begin{equation}
\partial\delta=\U(E_\delta\cup F_\delta)
\end{equation}
et en faisant la suspension des cubes de $E_\delta\cup F_\delta$ par rapport à $c_\delta$ on se retrouve dans le cas particulier d'une suspension tubulaire par rapport à un graphe réduit à un seul sommet, avec des constantes $\beta_-$, $\beta_+$, $\gamma$ et $\eta$ qui ne dépendent que de $p$ et $n$. Dès lors en posant
\begin{equation}
S_\delta=\left\{\S(\alpha,c_\delta)\colon\alpha\in E_\delta\cup F_\delta\right\}
\end{equation}
on sait d'après le lemme~\ref{lemmasuspensiontubulaire} que $S_\delta$ est un complexe $n$-dimensionnel vérifiant $\U(S_\delta)=\delta$ et
\begin{align}\label{equationsubdivisionC}
\overline{\R}(S_\delta)&\leq c_1\overline{\R}(E_\delta\cup F_\delta)&\underline{\R}(S_\delta)&\geq c_2\underline{\R}(E_\delta\cup F_\delta)
\end{align}
avec $c_1$ et $c_2$ qui dépendent uniquement de $n$ et $p$.

Posons finalement
\begin{equation}
S'=S\cup\bigcup_{\delta\in T}S_\delta.
\end{equation}
Par construction $S'$ est un complexe $n$-dimensionnel dyadique en surface de pas $r/p$, $p^n$-groupé, qui vérifie l'inégalité~\eqref{equationsubdivisionA}, ainsi que~\eqref{equationsubdivisionB} d'après~\eqref{equationsubdivisionC}.
\end{proof}

Dans la démonstration du théorème qui va suivre, pour $p>1$ et $\rho'>0$ donnés, en supposant que $\rho>3\sqrt{n}+\rho'$ et en appliquant le lemme~\ref{lemmasubdivision} à $S_1$ et $S_2$ (quitte à enlever ensuite les cubes rajoutés à $S_2$ qui ne sont pas dans $\overline{O}$) on peut aussi supposer que tous deux sont $p\times p$-groupés (pour alléger les notations on supposera qu'ils sont toujours unitaires, et on les notera encore $S_1$ et $S_2$) et les hypothèses du théorème resteront vérifiées avec $\rho'$. En particulier si on prend $p>3\rho'+3\sqrt{n}$, en notant
\begin{equation}
\begin{aligned}
A&=\left\{x\in O\colon\dist(x,\U(S_1))>\frac{p}{3}\right\}\\
O'&=O\setminus A
\end{aligned}
\end{equation}
alors $\U(S_1)\subset O'\subset O$ et toute composante connexe de $O'$ contient une unique composante connexe de $\U(S_1)$. Par ailleurs, en ajoutant à $S_1$ tous les cubes dyadiques unitaires possibles dans une même base situés à distance au moins $\rho'$ de $\U(S_2)$ et inclus dans $\overline{O}$ on a $\U(S_1)\supset A$ (car tout cube dyadique unitaire qui intersecte $\partial A\cap O$ est à distance au moins $p/3-\sqrt{n}>\rho'$ de $\U(S_1)$). Par conséquent, quitte à subdiviser préalablement $S_1$ et $S_2$, à ajouter des cubes à $S_1$ et à travailler séparément dans les composantes connexes de $O'$ on va supposer en plus que $\U(S_2)$ est connexe, ce qui implique en particulier que $\partial\U(S_2)$ l'est.

\subsection{Fusion en dimension $2$}

Traitons d'abord le cas de la fusion de deux complexes dyadiques bidimensionnels. On verra que la construction qu'on va expliciter va être encore réutilisable dans la suite de la démonstration en dimension plus grande.

\begin{lemma}\label{lemmafusion2}
Le théorème~\ref{theoremfusion} est vrai lorsque $n=2$.
\end{lemma}

\begin{proof}
Supposons donc que $n=2$, que les hypothèses du théorème de fusion sont vérifiées. Dans ces conditions, les portions de $S_1$ et $S_2$ qui nous intéressent sont deux graphes au sens de la définition~\ref{definitiongraphe}, les arêtes étant les faces (ou segments) qui sont respectivement dans la frontière de $O$ et de $\U(S_1)$. On va supposer qu'on dispose d'une distance $\rho$ assez grande entre les deux complexes pour faire toutes les constructions dont on va avoir besoin dans la démonstration qui suit, on verra à la fin que $\rho$ peut être borné indépendamment des complexes considérés.

Remarquons qu'une isométrie $\phi$ qui fait passer d'une base de $S_1$ à l'une de $S_2$ peut être décomposée en une rotation linéaire $r_\theta$ d'angle $\theta$ suivie d'une translation $\tau$:
\begin{equation}
\phi=\tau\circ r_\theta.
\end{equation}
En faisant éventuellement une permutation et\slash ou des inversions des vecteurs de la base de l'un des complexes, on peut aussi supposer que l'angle $\theta$ est compris entre $-\frac{\pi}{2}$ et $0$ modulo $\frac{\pi}{2}$. Pour tout entier $p>0$, la rotation $r_\theta$ peut alors être décomposée en $p$ rotations d'angle $\theta'=\frac{\theta}{p}\in\left[-\frac{\pi}{2p},0\right]$:
\begin{equation}
r_\theta=\left(r_{\theta/p}\right)^p=\left(r_{\theta'}\right)^p.
\end{equation}
Enfin en posant $\theta''=\frac{\pi}{4}+\frac{\theta'}{2}$ on constate que $\theta'=2\theta''$ modulo $\frac{\pi}{2}$ et là encore à des permutations\slash inversions près des vecteurs des bases considérées, $r_\theta$ peut s'écrire comme le produit de $2p$ rotations d'angle $\theta''\in\left[\frac{\pi}{4}-\frac{\pi}{4p},\frac{\pi}{4}\right]$. Si l'on construit successivement $2p-1$ complexes dyadiques en couronne autour de $S_2$ (en supposant $\rho$ assez grand) dont une base est l'image de celle du complexe précédant par une rotation d'angle $\theta''$, le problème de fusionner $S_1$ avec $S_2$ se ramène à faire fusionner successivement deux à deux ces complexes.

De la même façon en posant $\theta'''=\theta''-\frac{\pi}{4p}$, la composée de $2p$ rotations d'angle $\theta'''$ est une rotation d'angle égal à $\theta$ modulo $\frac{\pi}{2}$, avec 
\begin{equation}
\theta'''\in\left[\frac{\pi}{4}-\frac{\pi}{2p},\frac{\pi}{4}-\frac{\pi}{4p}\right].
\end{equation}
On supposera donc à partir de maintenant pour simplifier la démonstration que 
\begin{equation}
\theta\in\left[\theta_{\min},\theta_{\max}\right]\subset\left]0,\frac{\pi}{4}\right[
\end{equation}
où $\theta_{\min}$ peut être choisi arbitrairement proche de $\frac{\pi}{4}$, au prix d'un nombre de transitions $p$ assez grand qui ne dépend pas des complexes à faire fusionner.

À présent considérons la subdivision $4\times4$-groupée du complexe $S_2$ obtenue en utilisant le lemme~\ref{lemmasubdivision} avec $p=4$ (on appellera aussi $S_2$ le complexe obtenu, et on supposera encore qu'il est unitaire pour ne pas surcharger les notations). Chaque couple de faces distinctes de $\F_\partial(S_2)$ ayant un point en commun peuvent former un angle plat, aigu ou obtus. Dans ce qui suit, pour deux points $a=(a_1,a_2)$ et $b=(b_1,b_2)$ de $\mathbb{R}^2$ on notera $\dist_{\max}(a,b)$ la distance dite habituellement <<~du maximum~>> qui les sépare, définie par
\begin{equation}
\dist_{\max}(a,b)=\max\left(\vert a_1-b_1\vert,\vert a_2-b_2\vert\right).
\end{equation}
Considérons l'ensemble des points de coordonnées $\mathbb{Z}^2+\left(\frac{1}{2},\frac{1}{2}\right)$ dans une base de $S_2$ (les centres des cubes de $S_2$ en font par exemple partie), et notons $T$ le sous-ensemble de ces points qui sont à l'extérieur de $\U(S_2)$ et à distance $\frac{1}{2}$ pour la distance $\dist_{\max}$:
\begin{equation}
T=\left\{z\in\mathbb{Z}^2+\left(\frac{1}{2},\frac{1}{2}\right)\colon\dist_{\max}(z,\U(S_2))=\frac{1}{2}\right\}.
\end{equation}
Cet ensemble de sommets peut être muni naturellement d'une structure de graphe, en posant que les voisins d'un point de $T$ sont les points de $T$ à distance minimale:
\begin{multline}
A=\left\{\{x,y\}\subset T\colon x\in T,y\in T\setminus\{x\}\text{ et }\right.\\\left.\forall z\in T\setminus\{x\},\Vert z-x\Vert\geq\Vert y-x\Vert\right\}
\end{multline}
Notons $G$ le graphe formé du couple $(T,A)$ et vérifions rapidement que $G$ est linéaire cyclique. D'abord tout point de $T$ a au moins deux voisins car si $x\in T$ alors il existe $x'\in\F_0(S_2)$ tel que $x'\in\U(S_2)$ et $\dist_{\max}(x',x)=\frac{1}{2}$. On a vu que $x'$ était le sommet d'un angle plat, aigu ou obtus, et on peut trouver $y'$ et $z'$ distincts dans $\F_0(S_2)\cap\partial\U(S_2)$ tels que $\dist(x',y')=\dist(x',z')=1$. puisque l'on a subdivisé quatre fois alors nécessairement seul l'un des trois sommets $(x',y',z')$ de la frontière du complexe peut être le sommet d'un angle non plat. Par conséquent on peut alors trouver $y$ et $z$ distincts et distincts de $x$ dans $\mathbb{Z}^2+\left(\frac{1}{2},\frac{1}{2}\right)$ tels que $\dist_{\max}(y',y)=\frac{1}{2}$, $\dist_{\max}(z',z)=\frac{1}{2}$ et $\dist(x,y)=\dist(x,z)=1$. Par conséquent $x$ a au moins deux voisins: $y$ et $z$.

Raisonnons à présent par l'absurde et supposons que $x$ a au moins trois voisins distincts $a$, $b$ et $c$. Par construction $a$, $b$ et $c$ sont à distance $1$ de $x$ donc forment les sommets d'un triangle rectangle isocèle dont $x$ est le milieu de l'hypoténuse et par exemple $a$ est l'angle droit. Considérons alors l'union $U$ des boules fermées de centres $x$, $a$, $b$ et $c$ et de rayon $\frac{1}{2}$ pour la distance $\dist_{\max}$ (il s'agit en fait de quatre cubes dyadiques dans une base de $S_2$). Puisque l'on a subdivisé $S_2$ quatre fois alors $\U(S_2)\cap U$ est forcément compris d'un seul côté de la droite $(b,c)$, ce qui implique soit que $\dist_{\max}(a,\U(S_2))=\frac{3}{2}$, soit que $\dist_{\max}(a,\U(S_2))=0$. Dans tous les cas, on a une contradiction avec la définition des points de $T$. 

On vient donc de montrer que tout sommet du graphe $G$ a exactement deux voisins. On a en outre supposé que $\partial\U(S_2)$ est connexe, par conséquent $G$ est cyclique. Et puisque $G$ a au minimum huit sommets (ce serait le cas si $S_2$ n'était composé que d'un seul cube) alors il ne contient aucun $3$-cycle. Par conséquent $G$ est bien linéaire cyclique.

Soient $u$ et $v$ deux réels tels que $1<u<v<\frac{3}{2}$. On définit une bande en forme de couronne autour de $S_2$ à distance comprise entre $u$ et $v$:
\begin{equation}
\K_{u,v}=\{x\in\mathbb{R}^n\colon u<\dist(x,\U(S_2))<v\}.
\end{equation}
La frontière de cette bande a une composante connexe <<~extérieure~>> située à distance $v$ de $\U(S_2)$:
\begin{equation}
\K_v=\{x\in\mathbb{R}^n\colon\dist(x,\U(S_2))=v\}.
\end{equation}
Lorsque $v$ est fixé, $\forall u<v$ on peut trouver un entier $p$ tel que $\frac{\sqrt{2}}{p}<v-u$, donc en subdivisant $p$ fois $S_1$, si $\delta\in S_1$ et $\delta\cap\K_v\neq\emptyset$ alors $\dist(\delta,\U(S_2))>u$. On va alors compléter $S_1$ en lui rajoutant tous les cubes inclus dans $\overline{O}$, de pas $\frac{1}{p}$ dans la même base et qui contiennent un point à distance au moins $v$ de $\U(S_2)$:
\begin{equation}
S'_1=S_1\cup\left\{\delta\colon\delta\subset\overline{O}\text{ et }\sup_{x\in\delta}\dist(x,\U(S_2))\geq v\right\}.
\end{equation}
On remarquera que par construction, tous les polyèdres de $S'_1$ sont à distance au moins $u$ de $\U(S_2)$, donc que $\partial\U(S'_1)$ a une composante connexe incluse dans $\K_{u,v}$. Là encore, pour simplifier les notations on notera encore $S_1$ à la place de $S'_1$, puisque $S_1$ vérifie aussi les hypothèses du théorème avec une autre constante $\rho$.

Notons $T'$ l'ensemble des sommets du graphe $G$ qui sont alignés avec leurs deux voisins. Remarquons que pour toute arête $\{x,y\}\in A$ du graphe $G$, $x\in T'$ ou $y\in T'$ car on avait subdivisé $S_2$ quatre fois. Pour tout sommet $x\in T'$ soit $(d)$ la droite perpendiculaire à celle qu'il forme avec ses deux voisins, on notera $\C_x$ le cône de sommet $x$ et dont les points forment un angle (non orienté) compris entre $0$ et $\theta_{\min}$ avec la droite $(d)$:
\begin{equation}
\C_x=\left\{z\in\mathbb{R}^n\colon\widehat{\left(x,y\right),\left(d\right)}\in\left[0,\theta_{\min}\right]\right\}.
\end{equation}
Nous allons à présent discuter des différentes configurations possibles des sommets du graphe.

\subsubsection*{Cas d'un angle plat}

\begin{center}\fbox{\includegraphics[width=0.9\linewidth]{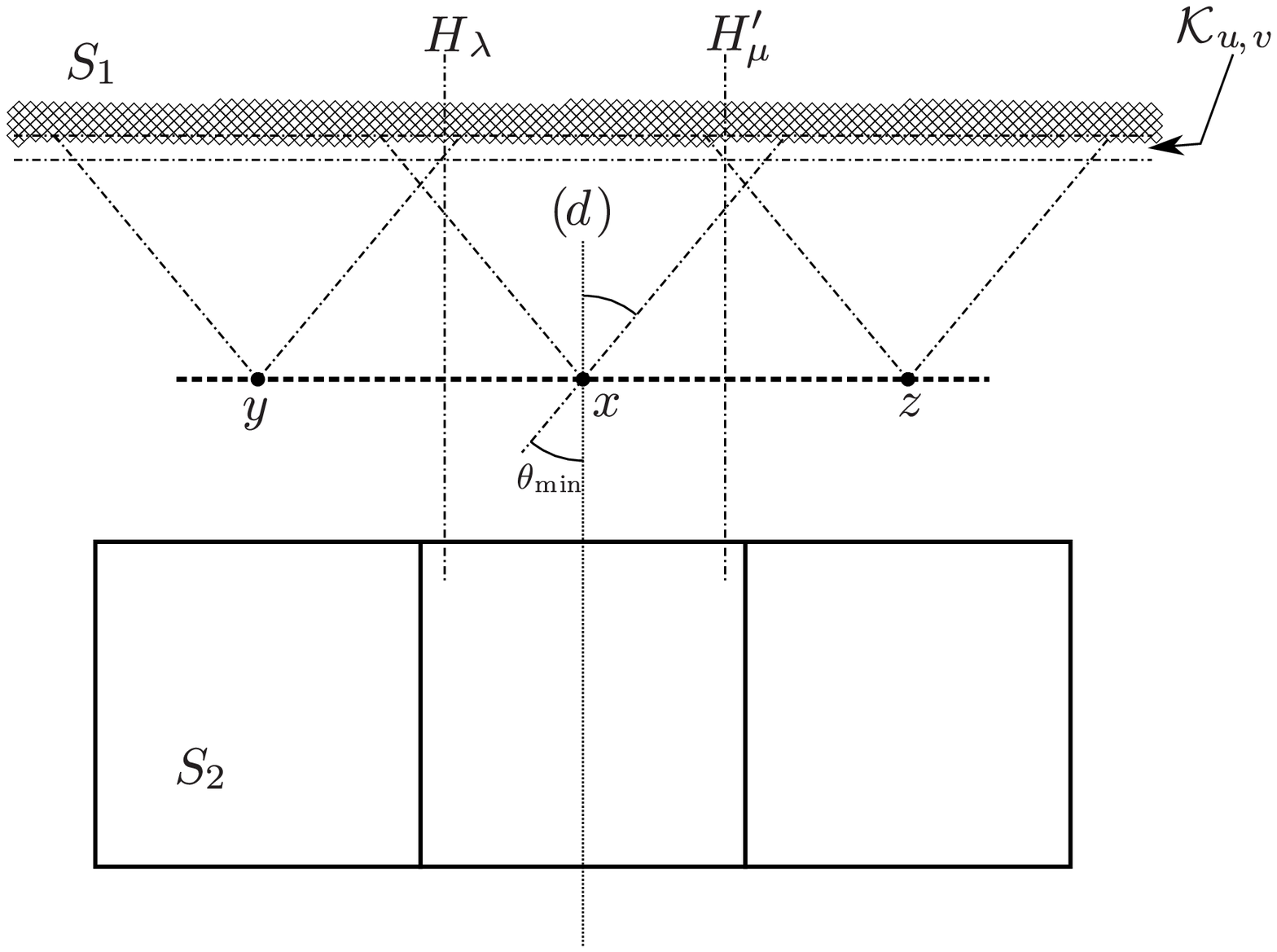}}\end{center}

Soient $x$ et $y$ deux sommets voisins du graphe tous deux alignés avec leurs voisins ($x,y\in T'$).
Remarquons d'abord que puisque l'on a supposé que $\theta\in[\theta_{\min},\theta_{\max}]$ alors $\partial S_1\cap\K_{u,v}\cap\C_x$ est en occlusion simple par rapport à $x$ (respectivement $\partial S_1\cap\K_{u,v}\cap\C_y$ par rapport à $y$). En outre, $\partial S_1\cap\K_{u,v}\cap\C_x\cap\C_y$ est en occlusion simple par rapport au segment $[x,y]$: en effet, remarquons que les faces de $S_1$ qui sont dans cet ensemble forment un angle compris entre $\theta_{\min}$ et $\frac{\pi}{4}$ avec la normale à la droite $(x,y)$, donc coupent la droite $(x,y)$ en un point extérieur au segment $[x,y]$ car $u>1$.

En choisissant $\theta_{\min}$ assez proche de $\frac{\pi}{4}$ il est possible d'obtenir que $\tan\theta_{\min}>\frac{1}{2}$, et en prenant
\begin{equation}\label{equationfusion2A}
v>\frac{1}{2}+\frac{1}{2\tan\theta_{\min}}
\end{equation}
on obtient que $\K_v\cap\C_x\cap\C_y$ est un segment non vide parallèle à $[x,y]$, de longueur $\nu>0$ qui ne dépend que de $v$ et $\theta_{\min}$. On notera $H_\lambda$ une droite perpendiculaire à ce segment, à distance $\lambda$ de l'une de ses extrémité (pour $\lambda\in[0,\nu])$.

\subsubsection*{Cas d'un angle aigu}

\begin{center}\fbox{\includegraphics[width=0.9\linewidth]{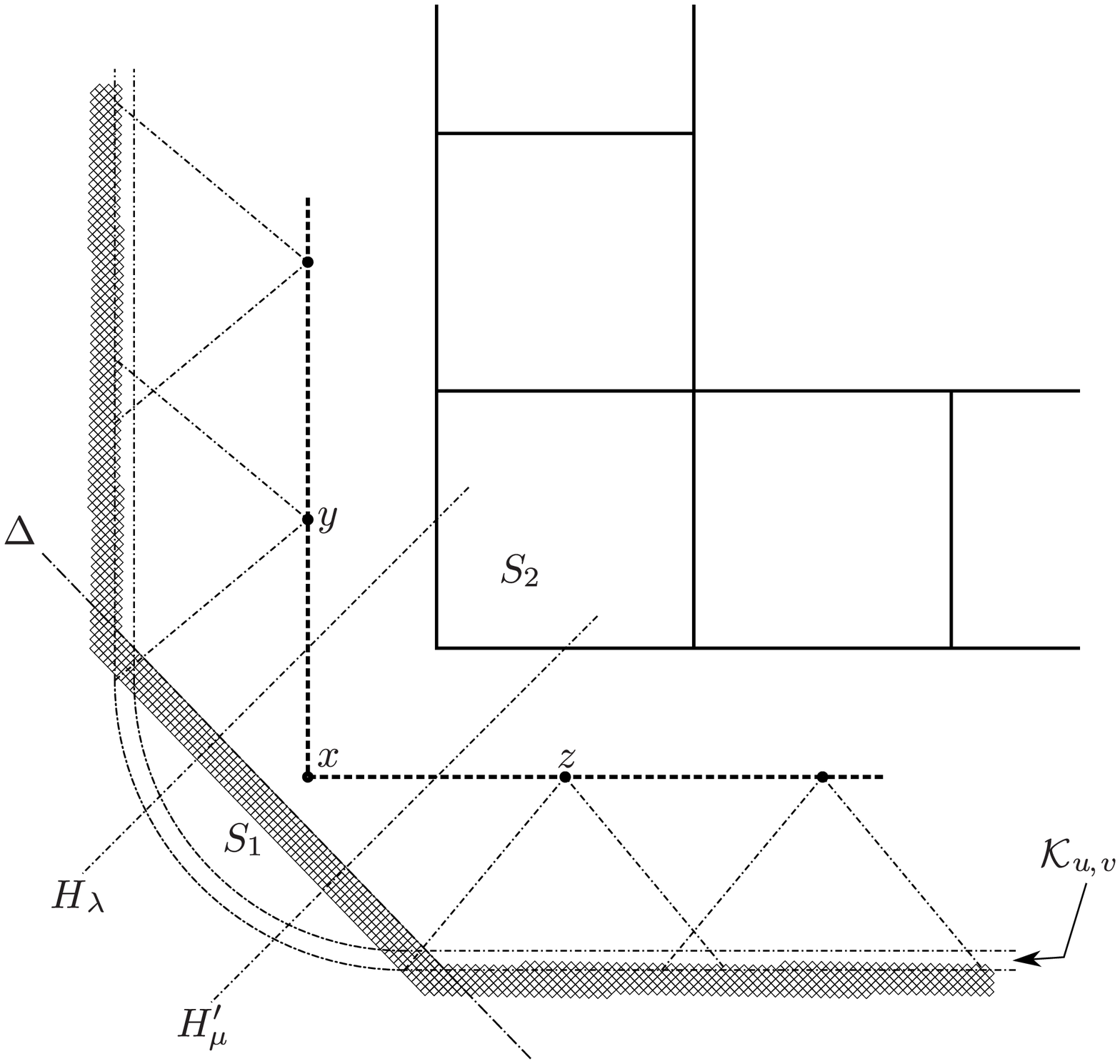}}\end{center}

Traitons à présent le cas des angles aigus: soit $x$ le sommet de l'angle, et $y$ et $z$ ses deux voisins. Si $\theta_{\min}$ est suffisamment proche de $\frac{\pi}{4}$, il est possible de trouver une droite $\Delta$ parallèle à l'un des vecteurs de la base canonique de $S_1$ et telle qu'elle passe par des sommets de la frontière de $S_1$, situés respectivement dans $\K_{u,v}\cap\C_y$ et $\K_{u,v}\cap\C_z$. Si on suppose qu'on a assez subdivisé $S_1$ (il suffit que son pas soit inférieur à $\frac{1}{100}$ par exemple) on peut aussi imposer que la droite $\Delta$ ne coupe pas les segments $[x,y]$ et $[x,z]$ et même, on peut trouver une constante $C>0$ ne dépendant que du pas de subdivision et de $\theta_{\min}$ telle qu'on puisse imposer $\dist(x,\Delta)>C$.

Notons $\Delta^+$ le demi-plan délimité par $\Delta$ qui contient $x$, $\Delta^-$ l'autre demi-plan. On rajoute alors à $S_1$ tous les cubes qui sont dans $\Delta^-\cap\{t\in\mathbb{R}^n\colon\dist(t,\U(S_2))\leq v\}$. On prend cette fois-ci des droites $H_\lambda$ et $H'_\mu$ qui sont perpendiculaires à $\Delta$ et à distance respective $\lambda$ et $\mu$ de $x$, du côté respectivement de $y$ et de $z$, avec par exemple $\lambda,\mu\in[\frac{1}{3},\frac{2}{3}]$.

\subsubsection*{Cas d'un angle obtus}

\begin{center}\fbox{\includegraphics[width=0.9\linewidth]{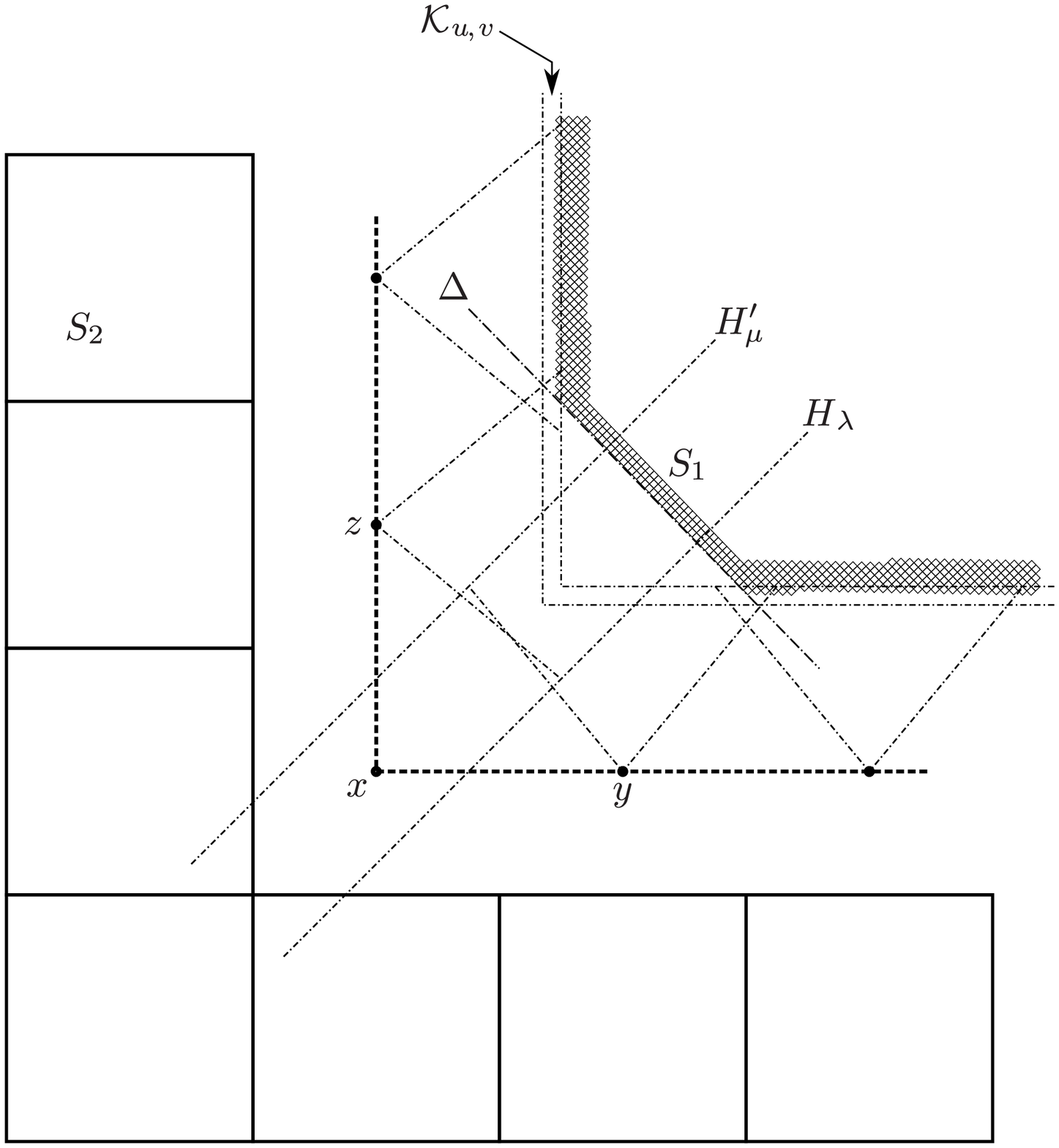}}\end{center}

Traitons à présent le cas des angles obtus: soit $x$ le sommet de l'angle, et $y$ et $z$ ses deux voisins. On va procéder quasiment de la même façon que pour les angles aigus. Notons $y'$ le voisin de $y$ qui n'est pas $x$, $z'$ celui de $z$ qui n'est pas $x$. Si $\theta_{\min}$ est suffisamment proche de $\frac{\pi}{4}$, il est possible de trouver une droite $\Delta$ parallèle à l'un des vecteurs de la base canonique de $S_1$ et telle qu'elle passe par des sommets de la frontière de $S_1$, situés respectivement dans $\K_{u,v}\cap\C_y\cap\C_{y'}$ et $\K_{u,v}\cap\C_z\cap\C_{z'}$. Si on suppose qu'on a assez subdivisé $S_1$ (il suffit que son pas soit inférieur à $\frac{1}{100}$ par exemple) on peut aussi imposer que la droite $\Delta$ ne coupe pas les segments $[x,y]$ et $[x,z]$ et même, on peut trouver une constante $C>0$ ne dépendant que du pas de subdivision et de $\theta_{\min}$ telle qu'on puisse imposer $\dist(x,\Delta)>C$.

Notons $\Delta^+$ le demi-plan délimité par $\Delta$ qui contient $x$, $\Delta^-$ l'autre demi-plan. On retire alors à $S_1$ tous les cubes qui sont dans $\Delta^-\cap\{t\in\mathbb{R}^n\colon\dist(t,\U(S_2))\geq u\}$ (cette opération est possible sans modification des complexes initiaux si on avait supposé $\rho>2$ par exemple, puisque dans ce cas on ne retire que ce qu'on avait surajouté). On prend des droites $H_\lambda$ et $H'_\mu$ qui sont perpendiculaires à $\Delta$ et à distance respective $\lambda$ et $\mu$ de $x$, du côté respectivement de $y$ et de $z$, pour $\lambda,\mu\in[\frac{1}{3},\frac{2}{3}]$ par exemple.

\subsubsection*{Conclusion}

Notons $\epsilon=\min\left(\nu,\frac{1}{3}\right)$ et $\Sigma$ le complexe composé des faces de $S_1$ et $S_2$ incluses dans leurs frontières en <<~vis-à-vis~>>:
\begin{equation}
\Sigma=\left\{F\in\F_1(S_1\cup S_2)\colon F\subset\partial\U(S_1)\text{ ou }F\subset\partial\U(S_2)\cap\overline{O}\right\}.
\end{equation}
Il est clair que $\Sigma$ est $\epsilon$-tubulaire par rapport à $G$ et $O$ en utilisant les droites $H_\lambda$ et $H_\mu$ qu'on a construites pour le découpage. Considérons les hypothèses du lemme~\ref{lemmasuspensiontubulaire}, par construction on obtient facilement les bornes suivantes:
\begin{align}
\alpha_+&=\alpha_-=1&\beta_-&>\frac{u}{2}&\beta_+&<2v&\gamma&<100&\eta&\in[\cos\theta_{\max},\cos\theta_{\min}].
\end{align}
Puisque comme on l'a vu, $\theta_{\min}$, $\theta_{\max}$, $u$ et $v$ peuvent être choisis indépendamment des complexes à fusionner il est clair qu'on peut trouver des constantes $\rho_+$, $\rho_-$ et un compact $K$ qui ne dépendent pas de $S_1$ et $S_2$ tels que~\eqref{equationsuspensiontubulaireA} et~\eqref{equationsuspensiontubulaireB} soient vérifiées, ce qui nous donne les constantes $c_1$ et $c_2$ recherchées.
\end{proof}

Cette démonstration va aussi nous servir à prouver le lemme~\ref{lemmacanalisations2}, utile dans la suite de la démonstration du théorème de fusion en dimension plus grande. Auparavant donnons quelques définitions.

\begin{definition}[Graphe canonique]\label{definitiongraphecanonique}
Soit $S$ un complexe dyadique de pas $r$, pour tout cube $\delta\in S$ notons $c_\delta$ son centre. Le graphe canonique $G=(T,A)$ associé à $S$ est défini naturellement par
\begin{align}
T=&\{c_\delta\colon\delta\in S\}&A=&\left\{\{x,y\}\subset T\colon\Vert x-y\Vert=r\right\}
\end{align}
et on dira qu'un complexe dyadique est connexe, cyclique ou linéaire si son graphe canonique l'est. 
\end{definition}

Plaçons-nous temporairement dans $\mathbb{R}^2$, pour $(x,y)\in\mathbb{R}^2$ et $r>0$ on notera
\begin{equation}
\begin{aligned}
\Delta(x,y,r)&=[x,x+r]\times[y,y+r]\\
\V(\Delta(x,y,r))&=\{\Delta(x+u,y+v,r)\colon(u,v)\in\{-r,0,r\}^2\}\\
\V^*(\Delta(x,y,r))&=\{\Delta(x+u,y+v,r)\colon(u,v)\in\{-r,0,r\}^2\text{ et }\vert u+v\vert=r\}.
\end{aligned}
\end{equation}
En remarquant que ces définitions ne dépendent pas du choix de la base de cube dyadique considérée, on les généralise à des complexes dyadiques bidimensionnels dans $\mathbb{R}^n$ ($n\geq2$) en se plaçant dans le $2$-plan affine correspondant. On dira que $\V(\delta)$ est l'ensemble des cubes voisins de $\delta$, $\V^*(\delta)$ celui des cubes tangents à $\delta$. On notera aussi
\begin{equation}
\extrem(S)=\left\{\delta\in S\colon\#(\V^*(\delta)\cap S)=1\right\}.
\end{equation}

Si on considère un complexe bidimensionnel $S$ et son graphe canonique associé $G=(T,A)$ on a par définition les équivalences
\begin{equation}
\forall (\alpha,\beta)\in S^2\colon\{c_\alpha,c_\beta\}\in A\Leftrightarrow\alpha\in\V^*(\beta)\Leftrightarrow\beta\in\V^*(\alpha)
\end{equation}
où $c_\alpha$ et $c_\beta$ désignent les centres de $\alpha$ et $\beta$. En remarquant que le graphe canonique d'un complexe dyadique n'a jamais de $3$-cycles on peut encore écrire:
\begin{align}
S\text{ linéaire}&\Leftrightarrow S\text{ connexe et }\forall\delta\in S\colon \#(\V^*(\delta)\cap S)\leq2\\
S\text{ linéaire cyclique}&\Leftrightarrow S\text{ connexe et }\forall\delta\in S\colon \#(\V^*(\delta)\cap S)=2\\
S\text{ linéaire non cyclique}&\Rightarrow\#\extrem(S)=2\\
S\text{ linéaire cyclique}&\Rightarrow\extrem(S)=\emptyset.
\end{align}

De façon à simplifier les énoncés par la suite, introduisons une notion pour désigner des complexes dyadiques bidimensionnels composés de sous-complexes linéaires en forme de sillons qui ne peuvent se rencontrer qu'en leurs extrémités.

\begin{definition}[Complexes sulciformes]\label{definitionsulciforme}
Soit $S$ un complexe dyadique bidimensionnel. On dira que $S$ est un complexe quasi-sulciforme s'il existe une partition $(S_i)_{i\in I}$ de $S$ en sous-complexes linéaires n'ayant pas de faces communes sur leur frontière hors de leurs extrémités:
\begin{equation}
\forall (i,j)\in I^2\colon i\neq j\Rightarrow\F_\partial(S_i\setminus\extrem(S_i))\cap\F_\partial(S_j\setminus\extrem(S_j))=\emptyset
\end{equation}
et dans ce cas on dira que les $S_i$ sont des sillons de $S$.

Si de plus les $S_i$ n'ont aucune face commune (extrémités comprises), autrement dit si
\begin{equation}
\forall (i,j)\in I^2\colon i\neq j\Rightarrow\F_\partial(S_i)\cap\F_\partial(S_j)=\emptyset
\end{equation}
alors on dira que $S$ est sulciforme. Si en outre tous les sillons de $S$ sont cycliques, on dira que $S$ est sulciforme en cycles.
\end{definition}

On pourra remarquer que dans le cas d'un complexe sulciforme, les sillons sont définis de manière unique. Dans le cas d'un complexe $4\times4$-groupé, le lemme suivant indique la possibilité d'y creuser des sillons (le complexe appelé $T$), de sorte que les cubes qui restent dans le complémentaire et forment leur <<~talus~>> (les complexes appelés $U_1$ et $U_2$) composent eux aussi des sillons, éventuellement deux fois plus larges, mais qui ne peuvent se rencontrer qu'au niveau de l'une de leurs extrémités.

\begin{lemma}[Laboureur]\label{lemmalaboureur}
Soit $S$ un complexe dyadique bidimensionnel $4\times4$-groupé. On peut construire une partition de $S$ par trois sous-complexes disjoints $T$, $U_1$ et $U_2$ tels que:
\begin{enumerate}
\item $S=T\sqcup U_1\sqcup U_2$;
\item $T$ est sulciforme en cycles et $\F_\partial(S)\subset\F_\partial(T)$, c'est à dire que les faces de la frontière de $S$ sont des faces de la frontière de $T$;
\item $U_1$ est sulciforme;
\item $U_2$ est la $2\times2$-subdivision d'un complexe quasi-sulciforme $U'_2$;
\item $\forall\delta\in U'_2\setminus \extrem(U'_2)\colon\U(U_1)\cap\delta=\emptyset$, c'est à dire que seules les extrémités des sillons de $U'_2$ peuvent toucher les cubes de $U_1$.
\end{enumerate}
\end{lemma}

\begin{figure}
\begin{center}\includegraphics[width=0.9\textwidth]{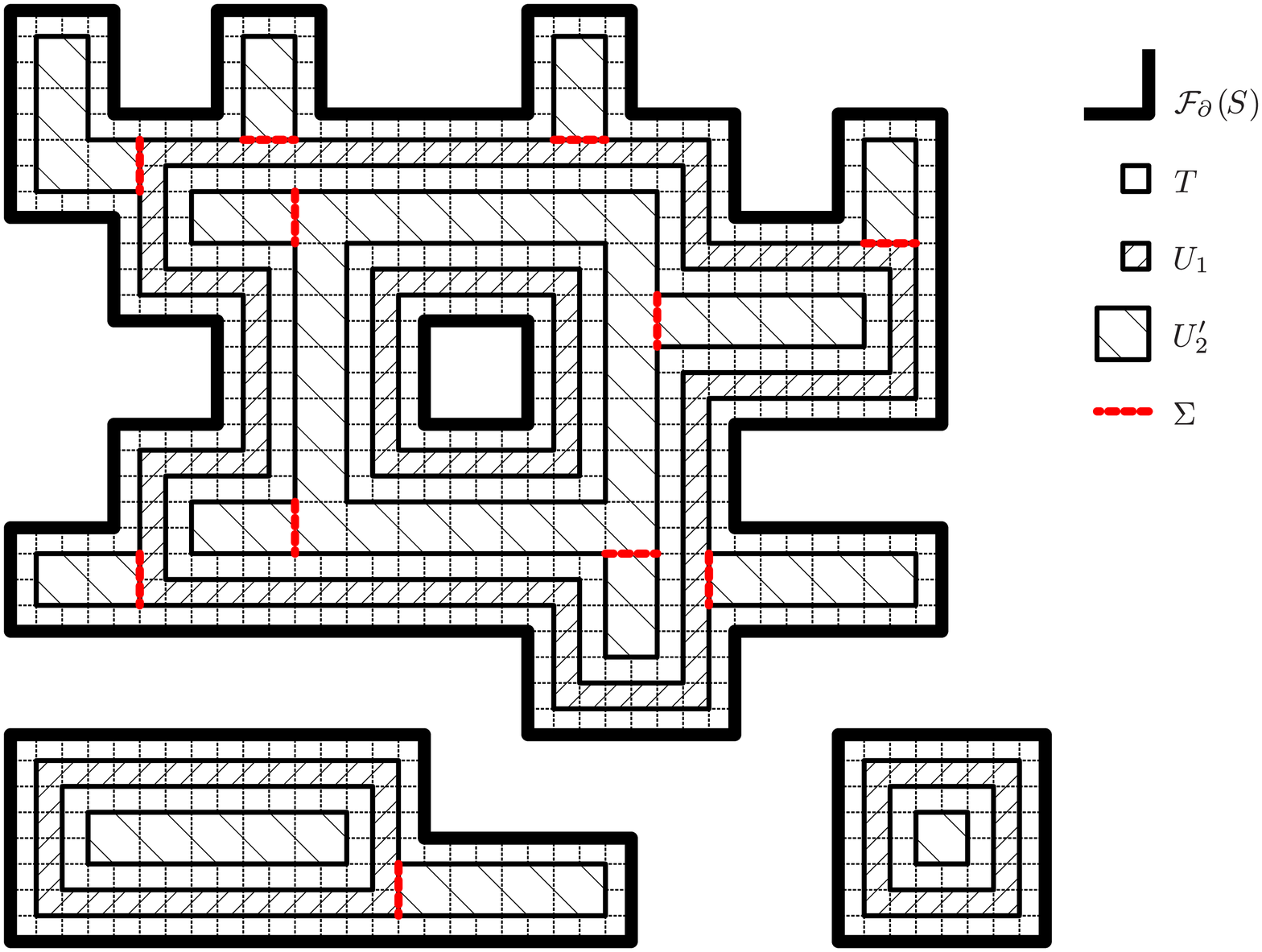}\end{center}
\caption{Lemme~\ref{lemmalaboureur} du laboureur appliqué à un complexe $4\times4$-groupé}
\label{figurelaboureurA}
\end{figure}

\begin{proof}
On se place dans $\mathbb{R}^2$ et on suppose que l'on dispose d'un complexe dyadique unitaire $S$ $4\times4$-groupé dans la base canonique de $\mathbb{R}^2$.

Notons $\mathbb{D}^2$ l'ensemble des cubes dyadiques bidimensionnels unitaires dans la même base, les sous-ensembles finis de $\mathbb{D}^2$ forment l'ensemble des complexes dyadiques unitaires par rapport à cette base. On notera $\mathfrak{D}$ l'application qui enlève à un complexe dyadique les cubes situés sur sa frontière:
\begin{equation}
\mathfrak{D}\colon
\begin{cases}
\P(\mathbb{D}^2)\rightarrow\P(\mathbb{D}^2)\\
T\mapsto\mathfrak{D}(T)=T\setminus\left\{\delta\in T\colon\#(\V(\delta)\cap T)<9\right\}.
\end{cases}
\end{equation}
On définit alors la suite $(S_n)_{n\in\mathbb{N}}$ de complexes:
\begin{equation}
\begin{cases}
S_0=S\\
S_{n+1}=\mathfrak{D}(S_n)
\end{cases}
\end{equation}
Pour tout complexe $T$, $\#\mathfrak{D}(T)<\#T$ donc il existe un entier $N$ tel que $\forall n\geq N\colon S_n=\emptyset$.

On va chercher dans un premier temps à montrer que 
\begin{equation}
\forall T\in\P(\mathbb{D}^2)\colon T\text{ est un complexe $4\times4$-groupé}\Rightarrow\mathfrak{D}^2(T)\text{ est $4\times4$-groupé.}
\end{equation}
Remarquons d'abord que puisque $\mathfrak{D}(T)$ est défini à partir du $1$-voisinage de chacun des cubes de $T$, si $T\in\P(\mathbb{D}^2)$, alors pour tous $p,q\in\mathbb{N}$, la restriction de $\mathfrak{D}(T)$ au carré $[p,q]^2$ (c'est à dire le sous-complexe
\begin{equation}
\mathfrak{D}(T)\vert_{[p,q]^2}=\left\{\delta\in\mathfrak{D}(T)\colon\open{\delta}\cap[p,q]^2\neq\emptyset\right\}
\end{equation}
voir la définition~\ref{definitionrestriction}) est connue à partir de la restriction de $T$ à $[p-1,q+1]^2$:
\begin{equation}
\mathfrak{D}(T)\vert_{[p,q]^2}=\mathfrak{D}\left(T\vert_{[p-1,q+1]^2}\right).
\end{equation}
Par conséquent si l'on connaît $T\vert_{[0,12]^2}$ alors on peut calculer $\mathfrak{D}^2(T)\vert_{[2,10]^2}$. Dès lors, si l'on dispose d'une méthode pour énumérer tous les $\mathfrak{D}^2(T\vert_{[0,12]^2})$ pour tous les complexes $T$ unitaires $4\times4$-groupés possibles (par rapport à l'origine $(0,0)$) tels que $\U(T)\subset[0,12]^2$, et vérifier que les sous-complexes obtenus restreints à $[2,10]^2$ sont tous eux-mêmes $4\times4$-groupés (par rapport à l'origine $(2,2)$), ceci suffira à prouver que $\mathfrak{D}^2$ conserve la propriété de $4\times4$-groupement.

Dans un second temps, définissons pour $n\geq0$:
\begin{align}
T_n&=S_n\setminus S_{n+1}&T&=\bigcup_{n\in\mathbb{N}}T_{2n}&U&=\bigcup_{n\in\mathbb{N}}T_{2n+1}.
\end{align}
Si l'on prouve que pour tout complexe $S$ $4\times4$-groupé par rapport à l'origine $(0,0)$, les cubes de la différence ensembliste $S\setminus\mathfrak{D}(S)$ ont chacun deux voisins tangents --- ceci peut se faire autour du carré central $[4,8]^2$, puisque $S$ est $4\times4$-groupé --- autrement dit si
\begin{equation}
\forall\delta\in(S\setminus\mathfrak{D}(S))\vert_{[4,8]^2}\colon\#\left(\V^*(\delta\cap(S\setminus\mathfrak{D}(S))\vert_{[3,9]^2}\right)=2
\end{equation}
alors on aura aussi démontré que $T$ est sulciforme en cycles.

Il suffira pour conclure d'extraire de $U$ deux sous-complexes $U_1$ et $U_2$ qui justifient le lemme. Pour cela, on donnera aussi un algorithme qui permet de dresser la liste de toutes les configurations possibles de $(\mathfrak{D}(S)\setminus\mathfrak{D}^2(S))\vert_{[3,9]^2}$ à une isométrie laissant le carré $[3,9]^2$ invariant. On discutera ensuite des différents cas trouvés, et de la façon d'extraire des sous-complexes $U_1$ et $U_2$ vérifiant les propriétés annoncées, ce découpage pouvant être déterminé localement par pavage.

Les algorithmes et le programme en C, ainsi que la fin de la preuve sont donnés en annexe.
\end{proof}

Pour alléger les énoncés à venir, on va généraliser la notion de sillons à des complexes non dyadiques en définissant la notion de canalisation.

\begin{definition}[Canalisations]\label{definitioncanalisation}
Soient $O$ un ouvert borné de $\mathbb{R}^n$, $S$ un complexe $n-1$-dimensionnel, $\G$ une famille de graphes linéaires et $\epsilon>0$. On dira que le couple $(S,\G)$ est une $\epsilon$-canalisation de $O$ si $\partial O=\U(S)$ et si pour toute composante connexe $\Omega$ de $O$ on peut trouver un sous-complexe $\Sigma$ de $S$ et un graphe $G=(T,A)$ de $\G$ vérifiant:
\begin{align}
T&\subset\Omega&\partial\Omega&=\U(\Sigma)&\text{$\Sigma$ est $\epsilon$-tubulaire par rapport à $\Omega$ et $G$}.
\end{align}
\end{definition}

On peut par exemple vérifier sans problème que les faces de la frontière d'un complexe sulciforme et les graphes canoniques de ses sillons forment une canalisation. Il nous reste encore à établir un dernier lemme qui va clore notre étude des complexes dyadiques bidimensionnels. Il sera utilisé dans le cadre de la démonstration du théorème de fusion en dimension $n>2$. Pour deux complexes bidimensionnels dyadiques unitaires $S_1$ et $S_2$ de $\mathbb{R}^2$, on notera $\theta$ un angle de rotation pris modulo $\pi/4$ qui fait passer d'une base de $S_1$ à une base de $S_2$.

\begin{lemma}[Canalisations complémentaires en dimension $2$]\label{lemmacanalisations2}
Il existe quatre constantes $\theta_{\min}$, $\theta_{\max}$, $p$, $\epsilon$ et un compact $K\subset]0,+\infty[$ tels que
\begin{align}
0&<\theta_{\min}<\theta_{\max}<\frac{\pi}{4}&p&\in\mathbb{N}\setminus\{0\}&\epsilon&>0
\end{align}
et pour tous complexes bidimensionnels dyadiques unitaires de $\mathbb{R}^2$ tels que $S_2$ est $28\times28$-groupé, $\theta\in[\theta_{\min},\theta_{\max}]$ et
\begin{equation}\label{equationcanalisations2A}
\min_{(x,y)\in\U(S_1)\times\U(S_2)}\Vert x-y\Vert>8
\end{equation}
alors il est possible de construire
\begin{itemize}
\item un ouvert $O$ tel que $\U(S_2)\subset O$ et 
\begin{equation}\label{equationcanalisations2B}
\inf_{(x,y)\in\U(S_1)\times O}\Vert x-y\Vert>\sqrt{2},
\end{equation}
\item un complexe $T_1$ dyadique de pas $\frac{1}{p}$ dans la même base que $S_1$ (on notera $O_1=\open{\U(T_1)}$),
\item un sous-complexe $T_2$ $7\times7$-groupé de $S_2$ (on notera $O_2=\open{\U(T_2)}$) tel que
\begin{equation}\label{equationcanalisations2C}
\F_\partial(S_2)\subset\F_\partial(T_2),
\end{equation}
\item un complexe unidimensionnel $\Sigma$ <<~placé à $\epsilon$-près~>> (voir plus bas),
\item trois familles $\G_1$, $\G_2$ et $\G_3$ de graphes linéaires
\end{itemize}
vérifiant les propriétés suivantes:
\begin{enumerate}
\item $\overline{O_1}\cap\overline{O_2}=\emptyset$ et $O_1\cup O_2\subset O$;
\item $\G_1$ et un sous-complexe de $\F_\partial(T_1)$ forment une canalisation de $O\setminus\overline{O_1}$;
\item $\G_2$ et un sous-complexe de $\Sigma\cup\F_\partial(T_2)$ forment une canalisation de $\open{\U(S_2)}\setminus(\overline{O_2}\cup\U(\Sigma))$;
\item $\G_3$ et un sous-complexe de $\F_\partial(T_1\cup T_2)$ forment une canalisation de $O\setminus\overline{(O_1\cup O_2)}$;
\item on peut trouver des constantes $\rho_+$ et $\rho_-$ telles que pour toutes les suspensions tubulaires des canalisations mentionnées ci-dessus, les relations~\eqref{equationsuspensiontubulaireA} et~\eqref{equationsuspensiontubulaireB} du lemme~\ref{lemmasuspensiontubulaire} sont vérifiées.
\end{enumerate}
\end{lemma}

Par <<~placé à $\epsilon$ près~>> on entend que pour tout segment $[x,y]\in\Sigma$, l'une de ses extrémités peut être déplacée à l'intérieur d'une boule de rayon $\epsilon$ sans que cela change quoi que ce soit aux cinq points exprimés plus haut.

\begin{figure}
\begin{center}\includegraphics[width=0.8\textwidth]{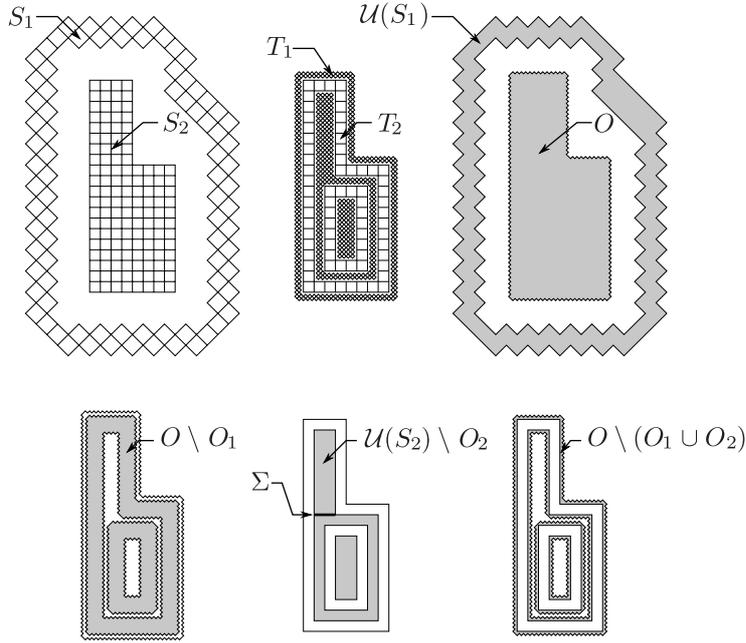}\end{center}
\caption{Canalisations complémentaires en dimension $2$}
\label{figurecanalisations2A}
\end{figure}

\begin{proof}
Appliquons d'abord le lemme~\ref{lemmalaboureur} au $7\times7$-groupement de $S_2$ (qui est lui-même $4\times4$-groupé car $S_2$ est $28\times28$-groupé par hypothèse): on obtient les trois sous-complexes $T$, $U_1$ et $U_2$ annoncés dans ce lemme. On notera $\G$ la famille des graphes canoniques associés aux sillons de $T$ (ce sont des graphes linéaires cycliques, car $T$ est sulciforme en cycles), et on définit $T_2$ comme la $7\times7$-subdivision de $T$ ($T_2$ est bien un sous-complexe de $S_2$ car $T$ est un sous-complexe du $7\times7$-groupement de $S_2$), et l'ouvert $O_2$ par $\overline{O_2}=\U(T_2)$.

Considérons, pour $r\in]0,7/3[$ et $p\in\mathbb{N}\setminus\{0\}$ tel que $\frac{7}{p}<\frac{7}{6\sqrt{2}}$ le complexe $T_1(r,p)$ de cubes dyadiques de pas $\frac{7}{p}$ dans la même base que $S_1$ et les ouverts $O$ et $O_1$ définis par:
\begin{equation}
\begin{aligned}
T_1(r,p)&=\left\{\delta\colon\forall x\in\delta,r<\dist(x,\U(T_2))<6\right\}\label{equationcanalisations2D}\\
\overline{O}&=\U\left(\left\{\delta\colon\forall x\in\delta,\dist(x,\U(S_2))<6\right\}\right)\\
\overline{O_1}&=\U(T_1(r,p)).
\end{aligned}
\end{equation}
Soit $\Omega$ une composante connexe de $O\setminus\overline{(O_1\cup O_2)}$: d'après les conditions sur $r$ et $p$, sa frontière a exactement deux composantes connexes, l'une portée par des faces de la frontière de $T_1(r,p)$, l'autre par des faces de la frontière de $T$ (on rappelle que $T$ est le $7\times7$-groupement de $T_2$). Notons respectivement $t_1$ et $t_2$ les complexes composés de ces faces: $t_1\subset\F_\partial(T_1(r,p))$ et $t_2\subset\F_\partial(T)$.

Utilisons à présent ce qu'on a vu dans la démonstration du lemme~\ref{lemmafusion2} en dimension $2$: il est possible de choisir $\theta_{\min}$, $\theta_{\max}$ et $v\in]1,2[$ qui vérifient~\eqref{equationfusion2A} de façon à ce que tout complexe dyadique unitaire situé à distance au moins $v$ de $\U(S_2)$ et de même base que $S_1$, suffisamment subdivisé, puisse être complété pour former avec $t_2$ et les cubes de $t_1$ en <<~vis-à-vis~>> une suspension tubulaire par rapport au graphe utilisé dans la démonstration. On peut prendre $r=\frac{13}{6}$, ce qui nous assure:
\begin{align}
r&>2>v&r&<\frac{7}{3}.
\end{align}

On choisira $p$ suffisamment grand (à la fois pour que $\frac{7}{p}<\frac{7}{6\sqrt{2}}$, et pour pouvoir faire la même construction que dans la démonstration du lemme~\ref{lemmafusion2}), et on notera $T_1$ le complexe obtenu à partir de $T_1(r,p)$ pour les valeurs correspondantes, auquel on a éventuellement ajouté les cubes supplémentaires nécessaires à la construction du tube. Ceci termine donc la démonstration du quatrième point, en notant $\G_3$ la famille des graphes utilisés pour les différentes composantes connexes $\Omega$ possibles.

Considérons maintenant les graphes de $\G$ et subdivisons-les naturellement sept fois (en rajoutant cinq sommets uniformément répartis sur chaque arête). Les graphes obtenus sont toujours linéaires cycliques, notons $\G_2$ la famille qu'ils composent. Soit $\Omega$ une composante connexe de $O\setminus \overline{O_1}$: par construction sa frontière est portée par un sous-complexe $s$ de $\F_\partial(T_1)$, et il existe un graphe $g\in\G_1$ dont les sommets sont dans $\Omega$. En fait, ce sous-complexe est formé par l'union des sous-complexes qui interviennent dans les suspensions tubulaires des composantes connexes de $O\setminus\overline{(O_2\cup O_1(r,p))}$ (ceux qu'on a appelés $t_1$ dans le cas précédant), et les arêtes du graphe $g$ sont situés à distance $7$ d'une arête parallèle des graphes du cas précédant, à l'intérieur de $\U(S_1)$. Il est alors facile de constater que puisque l'on a pu réaliser une suspension tubulaire avec $t_1$ dans le cas précédant, on peut encore le faire avec $s$. Ceci termine alors la démonstration du deuxième point. En outre, par construction d'après~\eqref{equationcanalisations2A} et~\eqref{equationcanalisations2D} on a
\begin{equation}
\min_{\substack{x\in\U(T_1)\\y\in\U(S_1)}}\Vert x-y\Vert>\min_{\substack{x\in\U(S_1)\\y\in\U(S_2)}}\Vert x-y\Vert-\min_{\substack{x\in\U(T_1)\\y\in\U(S_2)}}\Vert x-y\Vert>8-6=2>\sqrt{2}
\end{equation}
donc~\eqref{equationcanalisations2B} est vérifiée.

Revenons à présent à la construction qui intervient dans la démonstration du lemme~\ref{lemmalaboureur} (qui figure en annexe): on a découpé le complexe $T\setminus\bigcup_{t\in\T_1}t$ avec des segments de longueur $2$ longeant un $2\times2$ cube $\delta$ de $U'_2$, extrémal ou isolé ($\#\V^*(\delta)\leq1$) pour séparer les sillons de $U_1$ et $U_2$. On notera $\Sigma$ le complexe de dimension $1$ qui contient tous ces segments: puisque de façon évidente, le complexe formé des faces de la frontière d'un complexe dyadique linéaire est en situation de suspension tubulaire par rapport à son graphe canonique, le troisième point est démontré en prenant pour $\G_2$ la famille des graphes canoniques de $U_1$ et $U'_2$. Il est de plus tout à fait possible de découper les segments formant $\Sigma$ (par exemple en leur milieu), et il est clair qu'on peut donner $\epsilon>0$ ne dépendant pas de $S_1$ et $S_2$ tel qu'on puisse déplacer ces points à l'intérieur de boules de rayon $\epsilon$ sans que cela ne change l'existence de ces suspensions tubulaires.

Tous les autres points du lemmes sont vérifiés par construction (l'existence de $\rho_-$, $\rho_+$ et d'un compact $K$ vérifiant~\eqref{equationtubulaireA} et~\eqref{equationtubulaireB} qui ne dépendent pas de $S_1$ et $S_2$ ayant été démontrée avec le théorème de fusion en dimension $2$ pour $\Sigma$ fixé). De plus, en prenant $\epsilon$ assez petit, on peut trouver $K$ suffisamment grand qui ne dépend pas du choix de $\Sigma$ à $\epsilon$ près.
\end{proof}

\subsection{Fusion en dimension quelconque}

Notre objectif est de démontrer le théorème de fusion de complexes dyadiques par récurrence sur la dimension $n$ dans le cas d'une rotation planaire: en découpant les deux complexes à faire fusionner en <<~tranches~>> d'épaisseur $1$ et parallèles au plan de la rotation, on fera la fusion sur les portions planaires des frontières de la tranche en utilisant l'hypothèse de récurrence pour construire des <<~couvercles~>>, puis on va exhiber des graphes linéaires (en dimension $2$) à l'intérieur des tranches pour faire une suspension tubulaire remplissant les <<~boîtes~>> ainsi formée, et ainsi remplir tout l'espace entre les deux complexes. Le problème est qu'en dimension $n>2$, contrairement au cas précédent en dimension $2$, les faces de la frontière des deux complexes en vis-à-vis peuvent être très éloignées les unes des autres (en particulier si le complexe central a de larges morceaux de frontière formés de cubes alignés parallèlement au plan de la rotation). On n'a alors plus de borne sur les régularités des suspensions tubulaires, car la distance des points et des arêtes du graphe au tube n'est plus majorée indépendamment des complexes à faire fusionner.

Pour contourner le problème on se propose de généraliser le lemme~\ref{lemmacanalisations2} pour montrer qu'il est possible de creuser des canalisations complémentaires dans les deux complexes, de façon à ce qu'en les <<~encastrant~>> l'un dans l'autre, on dispose d'une borne supérieure à la distance séparant les faces en vis-à-vis et perpendiculaires au plan de rotation. Avant de commencer définissons la notion de restriction d'un complexe à un sous-ensemble, qui va nous être commode pour énoncer le lemme à venir.

\begin{definition}[Restriction d'un complexe]\label{definitionrestriction}
Soient $S$ un complexe $n$-dimen\-sion\-nel et $A$ une sous-partie de $\mathbb{R}^n$, on définit la restriction de $S$ à $A$ comme l'ensemble des intersections avec $A$ des polyèdres de $S$ dont l'intérieur est non disjoint de $A$:
\begin{equation}
S\vert_A=\left\{\delta\cap A\colon\delta\in S\text{ et }\open{\delta}\cap A\neq\emptyset\right\}.
\end{equation}
\end{definition}

Il est facile de vérifier que lorsque $A$ est une intersection finie de demi-espaces affines (par exemple si $A$ est un sous-espace affine, comme dans ce qui va suivre) alors $S\vert_A$ est encore un complexe lorsqu'elle est non vide.

Réciproquement on va utiliser des suspensions tubulaires bidimensionnelles pour construire des suspensions tubulaires $n$-dimensionnelles par rapport à un produit cartésien du tube bidimensionnel, flanqué de <<~couvercles~>> orthogonaux. Le lemme suivant permet d'évaluer les régularités extrêmes obtenues lors de cette opération.

On suppose que $n>2$, que $S$ est un complexe de dimension $1$ de $\mathbb{R}^2$, tubulaire par rapport à une graphe $G=(T,A)$ et un ouvert $O$. Pour $r>0$ on notera $S'$ le complexe obtenu par produit cartésien des cubes de $S$ par $[-r,r]^{n-2}$, $O'$ le produit cartésien de $O$ par $]-r,r[^{n-2}$, et $G'=(T',A')$ le plongement de $G$ dans $\mathbb{R}^n$:
\begin{equation}
\begin{aligned}
S'&=\left\{\delta\times[-r,r]^{n-2}\colon\delta\in S\right\}\\
O'&=O\times]-r,r[^{n-2}\\
T'&=\left\{(x,y,0,\ldots,0)\colon(x,y)\in T\right\}\\
A'&=\left\{\{(x,y,0,\ldots,0),(x',y',0,\ldots,0)\}\colon\{(x,y),(x',y')\}\in A\right\}.
\end{aligned}
\end{equation}
Pour $k\in\{0,\ldots,n-3\}$ et $r>0$ on notera
\begin{equation}
\begin{aligned}
O_{2k+1}&=\overline{O}\times[-r,r]^{k}\times\{-r\}\times[-r,r]^{n-3-k}\\
O_{2k+2}&=\overline{O}\times[-r,r]^{k}\times\{r\}\times[-r,r]^{n-3-k}
\end{aligned}
\end{equation}
et on supposera aussi qu'il existe une famille de $2(n-2)$ complexes $n-1$-di\-mensionnels $S_1,\ldots,S_{2(n-2)}$ de $\mathbb{R}^n$ tels que
\begin{equation}
\forall k\in\{1,\ldots,2(n-2)\}\colon\U(S_k)=O_k\text{ et }\F_\partial(S_k)=S'\vert_{O_k}.
\end{equation}

\begin{lemma}[Produit cartésien d'une suspension tubulaire]\label{lemmaproduitcartesien}
Pour tout compact $K\subset]0,+\infty[^9$ et pour tout $r>0$ il existe un compact $K'\subset]0,+\infty[^9$ tel que pour tout complexe unidimensionnel $\epsilon$-tubulaire $S$ par rapport à un ouvert $O\subset\mathbb{R}^2$ et un graphe $G$, si on peut trouver deux constantes $\rho_+$ et $\rho_-$ telles que les relations~\eqref{equationsuspensiontubulaireB} et~\eqref{equationsuspensiontubulaireC} du lemme~\ref{lemmasuspensiontubulaire} sont vérifiées avec en plus
\begin{equation}
\forall k\in\{1,\ldots,2(n-2)\}\colon\rho_+>\overline{\R}(S_k)>\underline{\R}(S_k)>\rho_-
\end{equation}
alors la famille de polyèdres $n-1$-dimensionnels
\begin{equation}
S''=S'\cup\bigcup_{k}S_k
\end{equation}
est un complexe $\epsilon$-tubulaire par rapport à $O'$ et $G'$ et il est possible de trouver $\rho'_-$ et $\rho'_+$ qui vérifient~\eqref{equationsuspensiontubulaireB} et~\eqref{equationsuspensiontubulaireC} avec le compact $K'$.
\end{lemma}

\begin{proof}
Remarquons d'abord que les polyèdres des $S_k$ sont tous dans des hyperplans parallèles aux arêtes du graphe. Dès lors les hypothèses de suspension tubulaire sont automatiquement vérifiées pour $S''$ puisque les hyperplans du découpage tubulaire $n$-dimensionnel sont le produit cartésien de ceux de la suspension bidimensionnelle par $\mathbb{R}^{n-2}$, donc perpendiculaires aux hyperplans contenant les polyèdres des $S_k$.

En notant $\alpha'_\pm$, $\beta'_\pm$, $\gamma'$ et $\eta'$ les constantes relatives à $O'$, $G'$ et $S''$ qui sont définies dans~\eqref{equationsuspensiontubulaireA} on peut déjà donner immédiatement les bornes suivantes par construction:
\begin{equation}
\begin{aligned}
\alpha'_-&=\alpha_-&\beta'_-&\geq\min(\beta_-,r)\\
\alpha'_+&=\alpha_+&\beta'_+&\leq\beta_++r\sqrt{n-2}\\
&&\eta'&\in\left[\frac{\beta_-}{\beta_-+r\sqrt{n-2}}\cdot\eta,1\right].
\end{aligned}
\end{equation}

En outre la mesure $n-1$-dimensionnelle des polyèdres de $S'$ est obtenue à partir de la mesure unidimensionnelle de ceux de $S$ multipliée par $(2r)^{n-2}$. Si l'on considère un sommet $x\in T$ et le sommet $x'$ correspondant de $T'$, par définition tous les polyèdres du découpage tubulaire des $S_k$ qui vont être mis en correspondance avec $x'$ par le choix de suspension tubulaire sont inclus dans le cylindre $C=B(x,s)\times [-r,r]^{n-2}$ où $B(x,s)$ désigne une boule $2$-dimensionnelle centrée en $x$ et de rayon $s=2\max(\alpha_+,\beta_+)$. Par ailleurs, pour tout $k\in\{1,\ldots,2(n-2)\}$ on a
\begin{equation}
\H^{n-1}(O_k\cap C)\leq us^2(2r)^{n-3}
\end{equation}
où $u$ désigne le volume de la boule unitaire en dimension $2$. On en tire
\begin{equation}
\gamma'\in\left[((2r)^{n-2}\gamma)^{\frac{1}{n-1}},((2r\gamma+us^2)(2r)^{n-3})^{\frac{1}{n-1}}\right].
\end{equation}
Dès lors, en prenant 
\begin{align}
\rho'_-&=\min(\rho_-,r)&\rho'_+&=\rho_++r\sqrt{n-2}
\end{align}
l'inégalité~\eqref{equationsuspensiontubulaireB} est vérifiée pour $S''$, ce qui termine la démonstration du lemme.
\end{proof}

Poursuivons à présent notre construction de complexes <<~imbriqués~>> en généralisant le lemme~\ref{lemmacanalisations2} des canalisations complémentaires en dimension $n\geq3$. L'idée ici est de considérer des complexes dyadiques dans des bases ayant subi une rotation planaire l'une par rapport à l'autre, parallèlement aux deux premiers vecteurs de ces bases. Dans ce cas, le lemme suivant indique qu'il est possible de construire des canalisations à chaque <<~étage~>> parallèlement au $2$-plan de la rotation, tout en restant suffisamment loin du complexe <<~extérieur~>> pour se laisser la place de construire des raccords.

\begin{figure}
\begin{center}\includegraphics[width=0.9\textwidth]{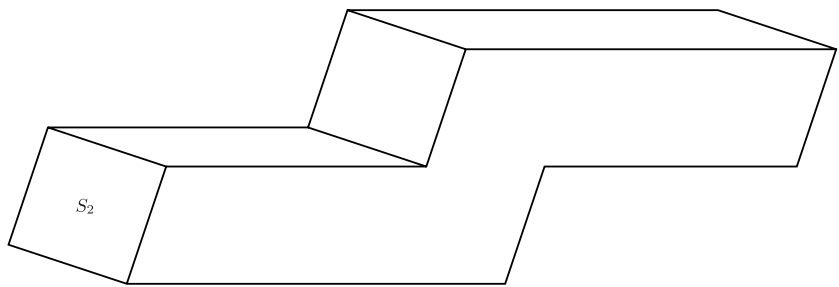}\end{center}
\caption{Un exemple de complexe $S_1$ en dimension $3$}
\label{figurecanalisationsnA}
\end{figure}

\begin{figure}
\begin{center}\includegraphics[width=0.9\textwidth]{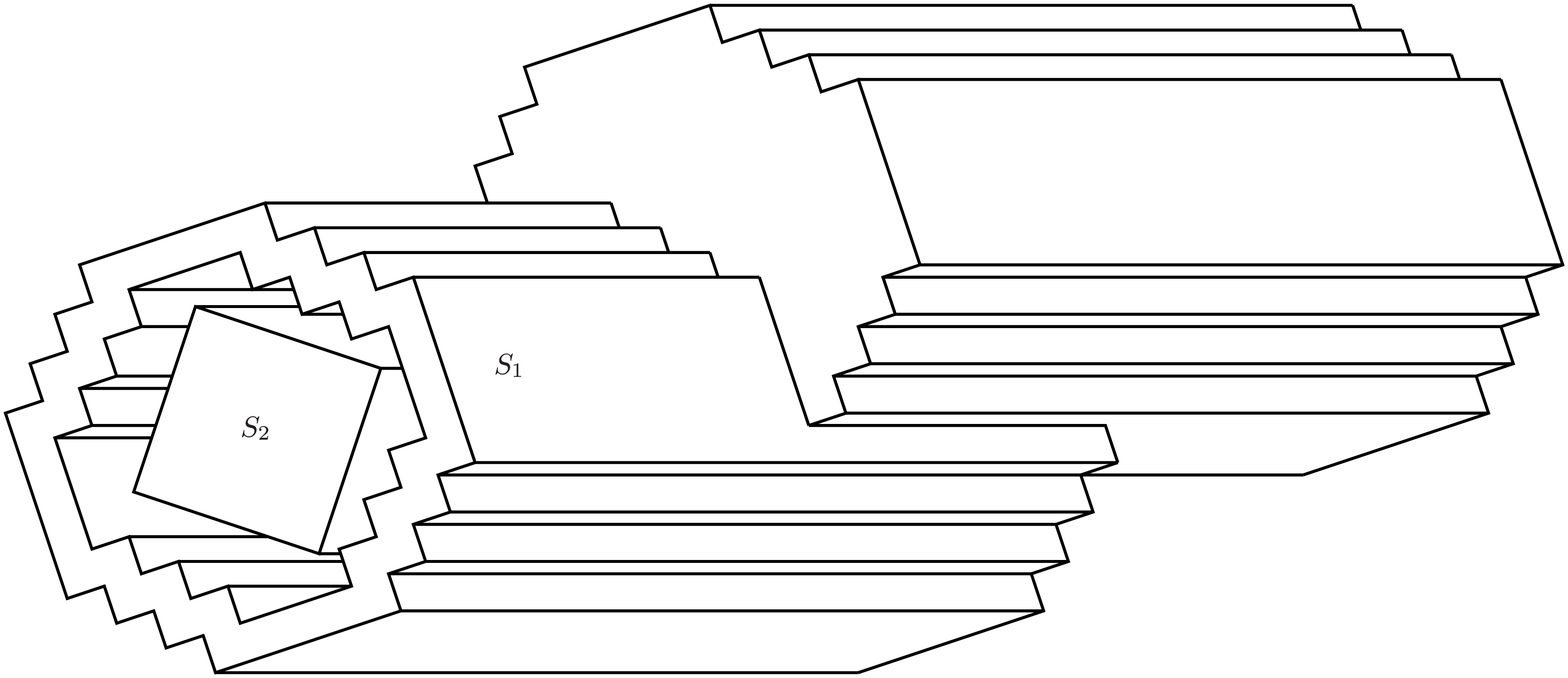}\end{center}
\caption{Le complexe $S_2$ représenté autour du complexe de la figure précédente}
\label{figurecanalisationsnB}
\end{figure}

\begin{figure}
\begin{center}\includegraphics[width=0.9\textwidth]{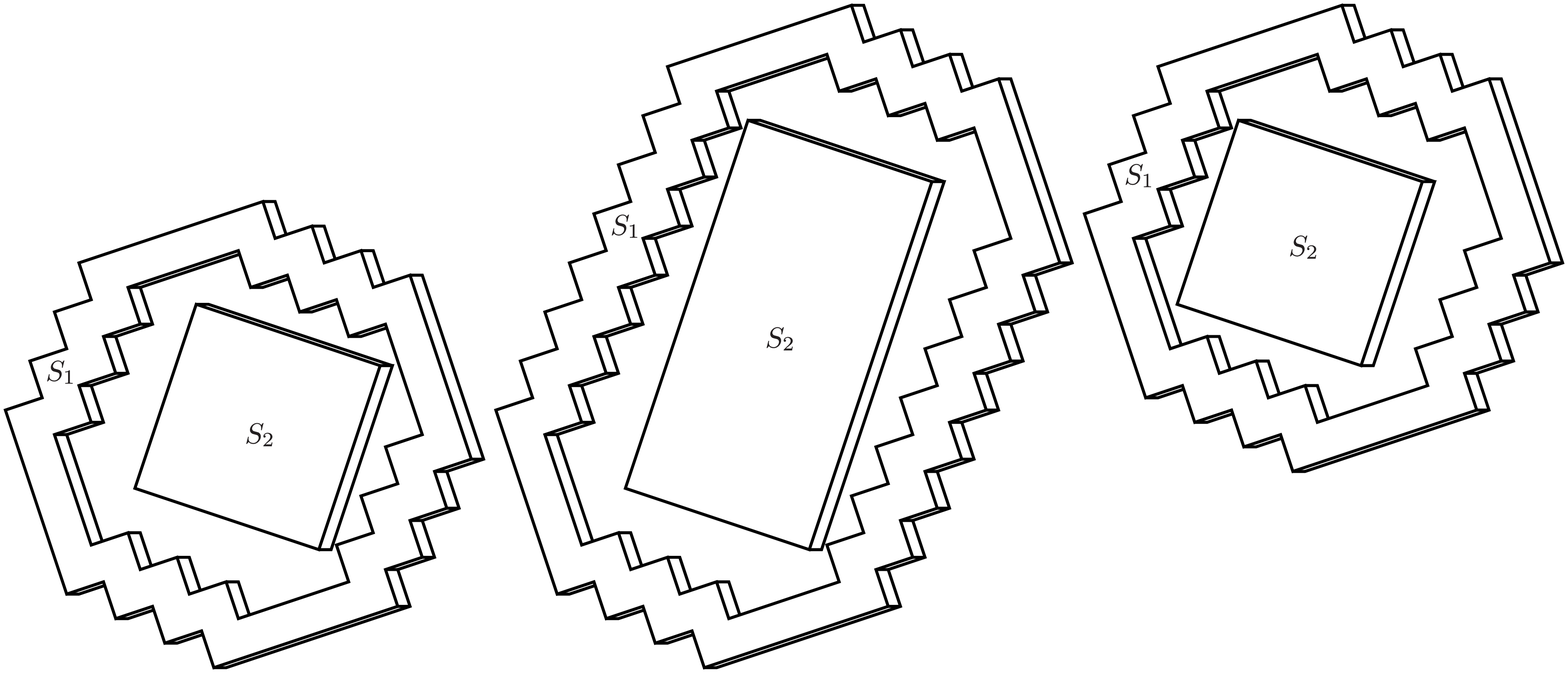}\end{center}
\caption{Trois <<~étages~>> dans lesquels on va creuser des canalisations}
\label{figurecanalisationsnC}
\end{figure}

\begin{figure}
\begin{center}\includegraphics[width=0.9\textwidth]{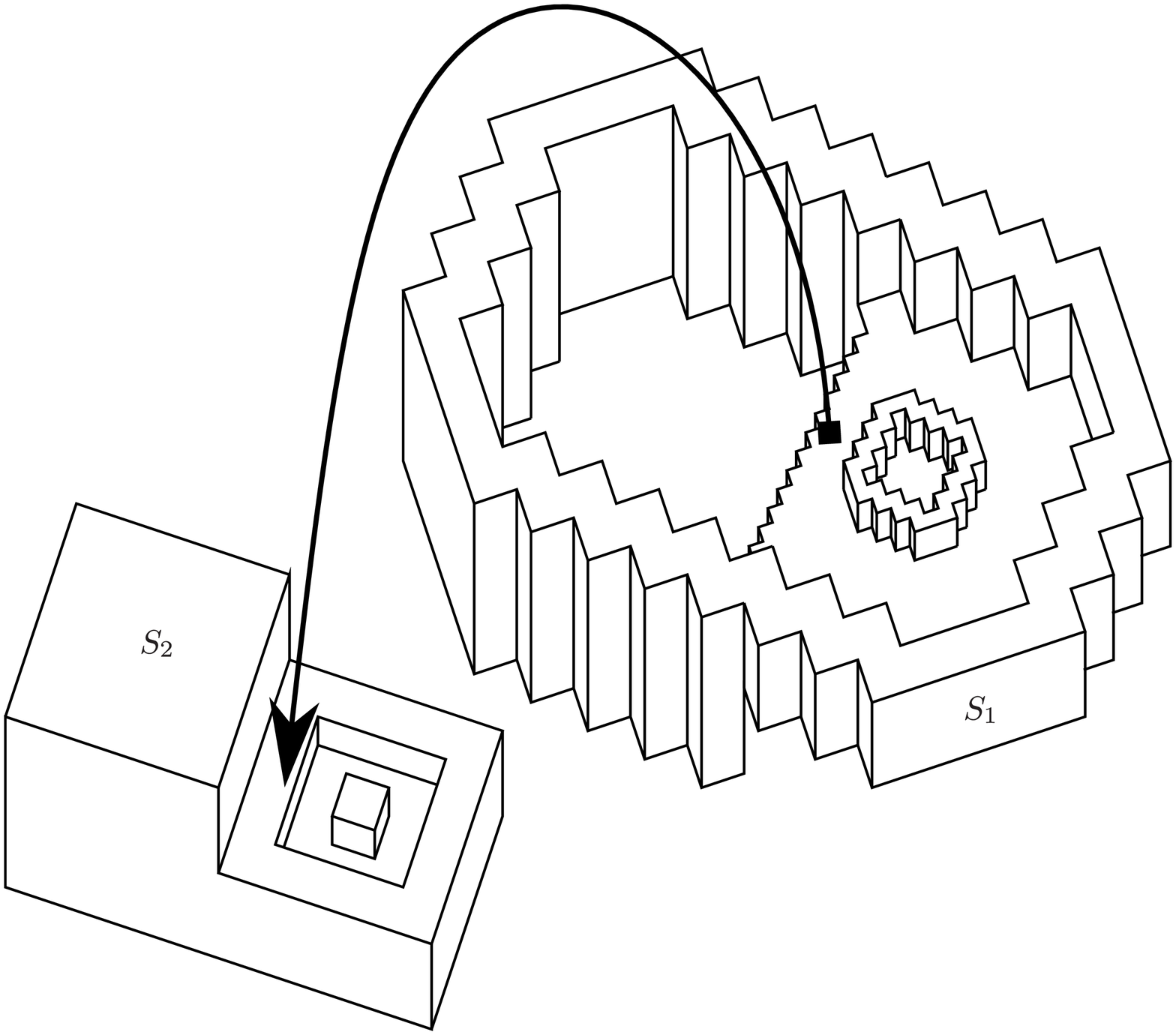}\end{center}
\caption{Exemples de canalisations complémentaires en dimension $3$, qui peuvent s'imbriquer l'une dans l'autre}
\label{figurecanalisationsnD}
\end{figure}

On va donc supposer que $n>2$, que $S_1$ et $S_2$ sont deux complexes dyadiques $n$-dimensionnels unitaires tels qu'une base canonique de $S_1$ soit l'image de celle de $S_2$ par une rotation affine d'angle $\theta\in\left[0,\frac{\pi}{4}\right[$ dans le $2$-plan $\mathbb{R}^2\times\{0\}^{n-2}$. Pour simplifier on supposera qu'une base canonique de $S_2$ est celle de $\mathbb{R}^n$, et pour $z=(z_3,\ldots,z_n)\in\mathbb{Z}^{n-2}$ on définit le $2$-plan affine
\begin{equation}
H_z=z+\mathbb{R}^2\times \left\{\frac{1}{2}\right\}^{n-2}.
\end{equation}

\begin{lemma}[Canalisations complémentaires en dimension quelconque]\label{lemmacanalisationsn}
Il existe quatre constantes $\theta_{\min}$, $\theta_{\max}$, $p$, $\epsilon$ et un compact $K\subset]0,+\infty[$ ne dépendant que de $n$ tels que
\begin{align}
0&<\theta_{\min}<\theta_{\max}<\frac{\pi}{4}&p&\in\mathbb{N}\setminus\{0\}&\epsilon&>0
\end{align}
et pour tous complexes dyadiques $n$-dimensionnels unitaires de $\mathbb{R}^n$ tels que $S_2$ est $28^n$-groupé, $\theta\in[\theta_{\min},\theta_{\max}]$ et
\begin{equation}
\min_{(x,y)\in\U(S_1)\times\U(S_2)}\Vert x-y\Vert>8
\end{equation}
alors il est possible de construire
\begin{itemize}
\item un ouvert $O$ tel que $\U(S_2)\subset O$ et 
\begin{equation}\label{equationcanalisationsnA}
\inf_{x\in O}\dist(x,\U(S_1))>\sqrt{2},
\end{equation}
\item un complexe $T_1$ dyadique de pas $\frac{1}{p}$, $1\times1\times p^{n-2}$-groupé dans la même base que $S_1$ (on notera $O_1=\open{\U(T_1)}$),
\item un sous-complexe $T_2$ de $S_2$ (on notera $O_2=\open{\U(T_2)}$) tel que
\begin{equation}\label{equationcanalisationsnB}
\max_{x\in\U(S_2)}\dist(x,\U(T_2))\leq2\sqrt{n},
\end{equation}
\item un complexe unidimensionnel $\Sigma$ <<~placé à $\epsilon$-près~>> (voir plus bas),
\item une famille $\G$ de graphes linéaires
\end{itemize}
vérifiant les propriétés suivantes:
\begin{enumerate}
\item $\overline{O_1}\cap\overline{O_2}=\emptyset$ et $O_1\cup O_2\subset O$;
\item pour tout $z\in\mathbb{Z}^{n-2}$ tel que $O\cap H_z\neq\emptyset$, il existe un sous-complexe de $(\F_\partial(T_1\cup T_2)\cup\Sigma)\vert_{H_z}$ et une sous-famille de graphes de $\G$ formant une canalisation de $(O\cap H_z)\setminus (O_1\cup O_2\cup\U(\Sigma))$ (en se plaçant dans le $2$-plan $H_z$);
\item on peut trouver des constantes $\rho_+$ et $\rho_-$ telles que pour toutes les suspensions tubulaires des canalisations mentionnées ci-dessus, les relations~\eqref{equationsuspensiontubulaireA} et~\eqref{equationsuspensiontubulaireB} du lemme~\ref{lemmasuspensiontubulaire} sont vérifiées.
\end{enumerate}
\end{lemma}

Par <<~placé à $\epsilon$ près~>> on entend que pour tout segment $[x,y]\in\Sigma$, l'une de ses extrémités peut être déplacée à l'intérieur d'une boule $2$-dimensionnelle parallèle à $H_z$ et de rayon $\epsilon$ sans que cela change quoi que ce soit aux trois points exprimés plus haut.

\begin{proof}
Prenons $\theta_{\min}$ et $\theta_{\max}$ égales aux constantes du lemme~\ref{lemmacanalisations2} et soient $S_1$ et $S_2$ deux complexes dyadiques vérifiant les hypothèses du lemme~\ref{lemmacanalisationsn}. On va construire $T_1$ et $T_2$ étage par étage, en découpant $S_1$ et $S_2$ par des tranches $I_z$ (pour $z\in(28\mathbb{Z})^{n-2}$) d'épaisseur $28$ autour de $H_z$:
\begin{equation}
I_z=\mathbb{R}^2\times (z+[0,28]^{n-2}).
\end{equation}
Fixons $z\in(28\mathbb{Z})^{n-2}$ tel que $S_2\vert_{H_z}\neq\emptyset$ et notons $U_1,\ldots,U_{2(n-2)}$ les morceaux de frontière plane de $I_z$ de la forme (pour $k\in\{1,\ldots,n-2\}$):
\begin{equation}
\begin{aligned}
U_{2k-1}&=\mathbb{R}^2\times[0,28]^k\times\{0\}\times[0,28]^{n-3-k}\\
U_{2k}&=\mathbb{R}^2\times[0,28]^k\times\{28\}\times[0,28]^{n-3-k}.
\end{aligned}
\end{equation}
Appelons aussi $U^*_1,\ldots,U^*_{2(n-2)}$ les <<~milieux~>> de ces morceaux (il s'agit de deux $2$-plans parallèles an plan de la rotation):
\begin{equation}
\begin{aligned}
U^*_{2k-1}&=\mathbb{R}^2\times\{14\}^k\times\{0\}\times\{14\}^{n-3-k}\\
U^*_{2k}&=\mathbb{R}^2\times\{14\}^k\times\{28\}\times\{14\}^{n-3-k}.
\end{aligned}
\end{equation}
Constatons que pour $1\leq j\leq2(n-2)$ les restrictions de $\F_\partial(S_1)$ et de $\F_\partial(S_2)$ à $U^*_j$ sont deux complexes dyadiques bidimensionnels unitaires, et notons-les respectivement $S_1(j)$ et $S_2(j)$:
\begin{align}
S_1(j)&=(\F_\partial(S_1))\vert_{U^*_j}&S_2(j)&=(\F_\partial(S_2))\vert_{U^*_j}.
\end{align}
En observant que les $U^*_j$ sont des $2$-plan affines parallèles, on va considérer les projection respectives $S_1^*(j)$ et $S_2^*(j)$ de $S_1(j)$ et $S_2(j)$ sur le $2$-plan vectoriel parallèle à $U^*_j$, qu'on identifie à $\mathbb{R}^2$ pour ne pas alourdir les notations. Considérons la famille de complexes bidimensionnels de $\mathbb{R}^2$:
\begin{equation}
\mathfrak{S}_2=\left\{t=\bigcap_{l\in K}S_2^*(l)\setminus\bigcup_{l\notin K}S_2^*(l)\colon K\subsetneqq\{1,\ldots,2n-4\}\text{ et }t\neq\emptyset\right\}.
\end{equation}
Par construction les complexes de $\mathfrak{S}_2$ sont dyadiques unitaires, bidimensionnels dans la base canonique de $\mathbb{R}^2$, $28\times28$-groupés, disjoints deux à deux et pour tout $j\in\{1,\ldots,2(n-2)\}$, $S_2^*(j)$ est un complexe formé par une union disjointe d'un ensemble $\mathfrak{S}_2(j)$ de certains complexes de $\mathfrak{S}_2$:
\begin{equation}
\forall j\in\left\{1,\ldots,2(n-2)\right\}\colon\exists\mathfrak{S}_2(j)\subset\mathfrak{S}_2,S_2^*(j)=\bigsqcup_{t\in\mathfrak{S}_2(j)}t.
\end{equation}

\begin{figure}
\begin{center}\includegraphics[width=0.9\textwidth]{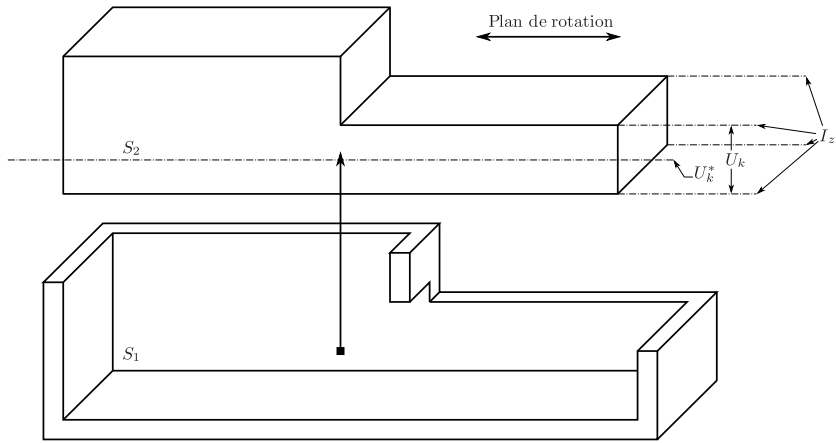}\end{center}
\caption{Un exemple en dimension plus grande}
\label{figurecanalisationsnE}
\end{figure}

\begin{figure}
\begin{center}\includegraphics[width=0.9\textwidth]{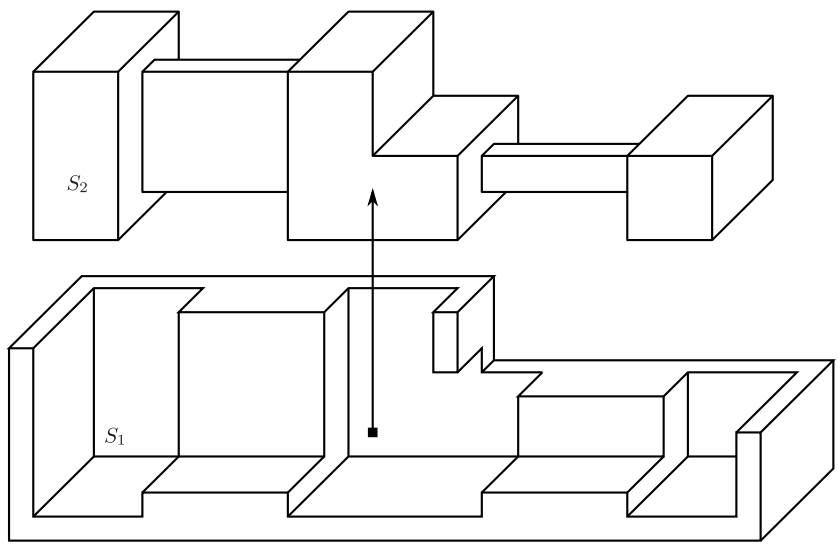}\end{center}
\caption{Canalisations complémentaires du cas précédent}
\label{figurecanalisationsnF}
\end{figure}

Pour $j$ fixé, soit $\Sigma_2\in\mathfrak{S}_2(j)$ et posons $\Sigma_1=S_1^*(j)$; appliquons alors le lemme~\ref{lemmacanalisations2} des canalisations complémentaires bidimensionnelles à $\Sigma_1$ et $\Sigma_2$ (respectivement notés $S_1$ et $S_2$ dans l'énoncé du lemme~\ref{lemmacanalisations2}): on obtient les trois complexes $\Theta_1(j,\Sigma_2)$, $\Theta_2(j,\Sigma_2)$ et $\Sigma(j,\Sigma_2)$ (respectivement notés $T_1$, $T_2$ et $\Sigma$ dans l'énoncé du lemme), et les trois familles de graphes $\G_1(j,\Sigma_2)$, $\G_2(j,\Sigma_2)$ et $\G(j,\Sigma_2)$. On va ôter de $\Theta_1(j,\Sigma_2)$ les cubes qui, dans la démonstration du lemme~\ref{lemmacanalisations2} forment une couronne autour de $\Sigma_2$, c'est à dire ne garder que ceux qui sont inclus dans l'union des cubes du complexe original $\Sigma_2$ (en d'autres termes, on ne garde que ceux qui peuvent <<~s'insérer~>> dans les trous qu'on a creusés dans $\Sigma_2$ dans la démonstration du lemme):
\begin{equation}
\Theta'_1(j,\Sigma_2)=\left\{\delta\in\Theta_1\colon\delta\subset\U(\Sigma_2)\right\}.
\end{equation}

Toujours d'après le lemme~\ref{lemmacanalisations2}, on sait que $\Theta_2(j,\Sigma_2)$ est sulciforme en cycles, et d'après~\eqref{equationcanalisations2C} que la construction effectuée préserve les cubes qui touchent la frontière de $\Theta_2$: 
\begin{equation}
\F_\partial(\Sigma_2)\subset\F_\partial(\Theta_2(j,\Sigma_2)).
\end{equation}
Si l'on considère un sillon cyclique de $\Theta_2(j,\Sigma_2)$ alors d'après le quatrième point du lemme, il existe un sous-complexe de $\F_\partial(\Theta_1(j,\Sigma_2))$ qui est située en <<vis-à-vis>>, à l'intérieur de $\Sigma_2$, et dont l'union des polyèdres (ici de dimension $1$) forme une courbe fermée de $\mathbb{R}^2$ qui est la frontière d'un ouvert borné. Si on note $\Omega(j,\Sigma_2)$ l'union de tous ces ouverts, alors $\Theta'_1(j,\Sigma_2)\subset\overline{\Omega(j,\Sigma_2)}\subset\U(\Sigma_2)$.

Dans ces conditions, nos complexes vérifient les propriétés suivantes:
\begin{enumerate}
\item $\Theta_2(j,\Sigma_2)$ est $7\times7$-groupé (d'après le lemme en dimension $2$);
\item $\G_1(j,\Sigma_2)$ et un sous-complexe de $\F_\partial(\Theta_1(j,\Sigma_2))$ forment une canalisation de $\Omega(j,\Sigma_2)\setminus\U(\Theta'_1(j,\Sigma_2))$;
\item $\G_2(j,\Sigma_2)$ et un sous-complexe de $\Sigma(j,\Sigma_2)\cup\F_\partial(\Theta_2(j,\Sigma_2))$ forment une canalisation de $\open{\U(\Sigma_2)}\setminus\U(\Theta_2(j,\Sigma_2)\cup\Sigma(j,\Sigma_2))$;
\item $\G_3(j,\Sigma_2)$ et un sous-complexe de $\F_\partial(\Theta'_1(j,\Sigma_2)\cup\Theta_2(j,\Sigma_2))$ forment une canalisation de $\open{\U(\Sigma_2)}\setminus\U(\Theta'_1(j,\Sigma_2)\cup\Theta_2(j,\Sigma_2))$.
\end{enumerate}
En se rappelant que $\mathfrak{S}_2$ est composée de complexes disjoints, on peut alors définir:
\begin{equation}
\begin{aligned}
\Theta_1(j)&=\bigsqcup_{\Sigma_2\in\mathfrak{S}_2(j)}\Theta'_1(j,\Sigma_2)&\G_1(j)&=\bigcup_{\Sigma_2\in\mathfrak{S}_2(j)}\G_1(j,\Sigma_2)\\
\Theta_2(j)&=\bigsqcup_{\Sigma_2\in\mathfrak{S}_2(j)}\Theta_2(j,\Sigma_2)&\G_2(j)&=\bigcup_{\Sigma_2\in\mathfrak{S}_2(j)}\G_2(j,\Sigma_2)\\
\Sigma(j)&=\bigsqcup_{\Sigma_2\in\mathfrak{S}_2(j)}\Sigma(j,\Sigma_2)&\G_3(j)&=\bigcup_{\Sigma_2\in\mathfrak{S}_2(j)}\G_3(j,\Sigma_2)\\
\Omega(j)&=\bigcup_{\Sigma_2\in\mathfrak{S}_2(j)}\Omega(j,\Sigma_2)
\end{aligned}
\end{equation}

\begin{figure}
\begin{center}\includegraphics[width=0.9\textwidth]{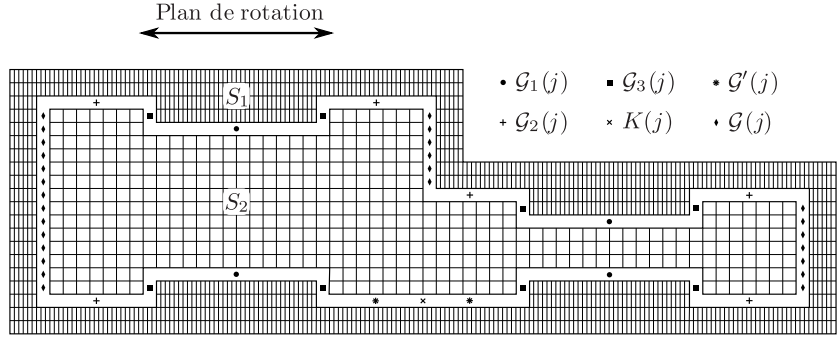}\end{center}
\caption{Les différents graphes utilisés pour les canalisations}
\label{figurecanalisationsnG}
\end{figure}

Les trois familles $\Theta_1(j)$, $\Theta_2(j)$ et $\Sigma(j)$ forment bien des complexes car chacun des éléments qui intervient dans les unions qui les définissent est inclus dans un unique $\U(\Sigma_2)\in\mathfrak{S}_2$. Ces ensembles vérifient en plus les propriétés suivantes:
\begin{enumerate}
\item $\Theta_2(j)$ est $7\times7$-groupé, $\Theta_2(j)\subset S^*_2(j)$, $\F_\partial(S^*_2(j))\subset\F_\partial(\Theta_2(j))$ et $\Omega(j)\subset\U(S^*_2(j))$;
\item $\G_1(j)$ et un sous-complexe de $\F_\partial(\Theta_1(j))$ forment une canalisation de $\Omega(j)\setminus\U(\Theta_1(j))$;
\item $\G_2(j)$ et un sous-complexe de $\Sigma(j)\cup\F_\partial(\Theta_2(j))$ forment une canalisation de $\open{\U(S^*_2(j))}\setminus\U(\Theta_2(j)\cup\Sigma(j))$;
\item $\G_3(j)$ et un sous-complexe de $\F_\partial(\Theta_1(j)\cup\Theta_2(j))$ forment une canalisation de $\open{\U(S^*_2(j))}\setminus\U(\Theta_1(j)\cup\Theta_2(j))$.
\end{enumerate}

Considérons à présent l'ensemble $K(j)$ des frontières communes des sillons des $\Theta_2(j,\Sigma_2)$, lorsque $\Sigma_2$ parcourt $\mathfrak{S}_2(j)$:
\begin{multline}
K(j)=\left\{s=\U(s_1)\cap\U(s_2)\colon\Sigma_2\in\mathfrak{S}_2(j),s\neq\emptyset,\right.\\\left.
\text{$s_1$ et $s_2$ sillons de $\Theta_2(j,\Sigma_2)$ et $s_1\neq s_2$}\right\}.
\end{multline}
$K(j)$ est composé de courbes formées de segments parallèles aux vecteurs de la base de $S^*_2$. Soient $p$ et $r$ les constantes qui interviennent dans la démonstration du lemme~\ref{lemmacanalisations2} en dimension $2$, et éventuellement en ajoutant $1$ à $p$ on s'assure qu'il est impair; définissons alors le complexe composé de cubes dyadiques de pas $\frac{1}{p}$ dans la même base que $S^*_1(j)$, situés à distance comprise entre $r$ et $6$ de $\U(S^*_2(j))$:
\begin{equation}
\Theta'_1(j)=\left\{\delta\colon r<\min_{x\in\delta}\dist(x,\U(S^*_2(j))<6\right\}.
\end{equation}
Puisque les courbes de $K(j)$ sont incluses dans $\U(S^*_2(j))$, elles ne rencontrent pas $\U(\Theta'_1(j))$. Considérons l'ensemble des extrémités de ces courbes (parmi celles qui ne sont pas fermées) qui ne font partie d'aucune autre courbe de $K(j)$, et construisons à partir de ces points l'ensemble des segments dont l'autre extrémité est le sommet le plus proche de $\U(S^*_2(j))$: appelons ce nouvel ensemble de segments $K'(j)$. Ces segments sont d'intérieurs disjoints deux à deux et des courbes de $K(j)$ (car les extrémités des segments formant les courbes de $K(j)$ sont au moins à distance $7$ les unes des autres, et situées à distance au plus $r+1<3$ d'un sommet de $S^*_2(j)$), dès lors l'ensemble des segments de $K'(j)$ et ceux formant les courbes de $K(j)$ (éventuellement en les découpant aux endroits où ils se rencontrent) forme un complexe unidimensionnel.

Notons $O(j)$ l'ouvert de $\mathbb{R}^2$ dont l'adhérence est l'union des cubes dyadiques de pas $\frac{1}{p}$ et à distance au plus $6$ de $\U(S^*_2(j))$:
\begin{equation}
\overline{O(j)}=\U\left(\left\{\delta\colon\min_{x\in\delta}\dist(x,\U(S^*_2(j)))<6\right\}\right).
\end{equation}
Par un argument similaire à celui utilisé dans la démonstration du lemme~\ref{lemmacanalisations2}, quitte à ajouter\slash enlever quelques cubes à $\Theta'_1(j)$, on peut trouver une famille $\G(j)$ de graphes linéaires formant avec le complexe unidimensionnel $K'(j)\cup\F_\partial(\Theta'_1(j))\cup\F_\partial(S^*_2(j))$ une canalisation de $O(j)\setminus\U(S^*_2(j))$. De plus, si l'on suppose que l'on a subdivisé de façon naturelle les segments de $K(j)$ en morceaux de longueur $1$ (on notera $K''(j)$ ce nouveau complexe) alors par le même argument, il est aussi possible de trouver une famille $\G'(j)$ de graphes formant avec un sous-complexe de $K'(j)\cup K''(j)\cup\F_\partial(\Theta'_1(j))$ une canalisation de $O(j)\setminus(\U(\Theta'_1(j))\cup\U(K'(j)\cup K''(j)))$. De plus, comme dans la démonstration du lemme~\ref{lemmacanalisations2} il est possible de déplacer l'une des extrémités des segments de $K'(j)$ et $K''(j)$ à l'intérieur d'une boule de rayon $\epsilon$ suffisamment petit pour que les régularités extrêmes des suspensions tubulaires restent comparables à celles de $S^*_1(j)$ et $S^*_2(j)$, avec des constantes uniformes.

Définissons encore:
\begin{equation}
\begin{aligned}
\Theta''_1(j)&=\Theta_1(j)\cup\Theta'_1(j)&\Sigma'(j)&=\Sigma(j)\cup K'(j)\cup K''(j)\\
\G'_1(j)&=\G_1(j)\cup\G'(j)&\G'_3(j)&=\G_3(j)\cup\G(j).
\end{aligned}
\end{equation}
On dispose à présent de deux complexes qui vérifient les propriétés suivantes:
\begin{enumerate}
\item $\G'_1(j)$ et un sous-complexe de $\F_\partial(\Theta''_1(j))\cup\Sigma'(j)$ forment une canalisation de $O(j)\setminus(\U(\Theta_1(j))\cup\U(\Sigma'(j)))$;
\item $\G'_3(j)$ et un sous-complexe de $\F_\partial(\Theta''_1(j)\cup\Theta_2(j))$ forment une canalisation de $\open{\U(S^*_2(j))}\setminus\U(\Theta''_1(j)\cup\Theta_2(j))$.
\end{enumerate}

Il est temps à présent d'utiliser toutes ces canalisations bidimensionnelles pour construire $T_1$ et $T_2$. On notera cette fois, pour $r_1<r_2$, $z\in(28\mathbb{Z})^{n-2}$ et $j\in\{1,\ldots,2(n-2)\}$:
\begin{equation}
V(r_1,r_2,z,j)=
\begin{cases}
z+[0,28]^k\times[r_1,r_2]\times[0,28]^{n-3-k}&\text{si $j=2k-1$}\\
z+[0,28]^k\times[28-r_2,28-r_1]\times[0,28]^{n-3-k}&\text{si $j=2k$}
\end{cases}
\end{equation}
c'est à dire, les $\mathbb{R}^2\times V(r_1,r_2,z,j)$ sont des <<~tranches~>> de $I_z$ à distance comprise entre $r_1$ et $r_2$ des portions planaires de sa frontière. Commençons par ôter de $S_2$ l'ensemble des cubes de sa frontière qui ne sont pas dans $\U(\Theta_2(j))\times V(0,2,z,j)$ (cela revient à enlever sur une épaisseur de deux cubes, parallèlement à la frontière de $I_z$, les cubes dont la projection sur $\mathbb{R}^2$ a été enlevée dans le lemme~\ref{lemmacanalisations2}) en posant:
\begin{equation}
t_2(z,j)=\left\{\delta\in S_2\vert_{\mathbb{R}^2\times V(0,2,z,j)}\colon\delta\subset\U(\Theta_2(j))\times V(0,2,z,j)\right\}.
\end{equation}
On définira aussi l'ouvert $O(z)$ par:
\begin{equation}
\overline{O(z)}=\bigcup_{j}\overline{O(j)}\times V(0,1,z,j).
\end{equation}

Soient $z'\in(28\mathbb{Z})^{n-2}$, et $j'\in\{1,\ldots,2(n-2)\}$ et notons $S$ et $S'$ les sous-complexes formés des cubes de $S_2$ ayant au moins un sommet commun avec ceux de $\F_\partial(S_2)$, restreints respectivement à $V(0,1,z,j)$ et $V(0,1,z',j')$. Dans ces conditions:
\begin{equation}
S\vert_{V(0,1,z',j')}=S'\vert_{V(0,1,z,j)}.
\end{equation}
Considérons alors les restrictions respectives de $t_2(z,j)$ à $V(0,1,z',j')$, et de $t_2(z',j')$ à $V(0,1,z,j)$: puisque d'après~\eqref{equationcanalisations2C} seules les restrictions des cubes ayant un sommet commun avec la frontière de $S_2$ qui n'avaient pas de sommet commun avec la frontière des $S^*_2(j)$ ont été enlevés à $S_2$ pour former les $t_2(z,j)$ quand on a utilisé le lemme~\ref{lemmacanalisations2} alors il vient
\begin{align}
S&\subset t_2(z,j)&S'&\subset t_2(z',j')
\end{align}
et donc
\begin{equation}\label{equationcanalisationsnC}
t_2(z,j)\vert_{V(0,1,j',z')\cap V(0,1,j,z)}=t_2(z',j')\vert_{V(0,1,j',z')\cap V(0,1,j,z)}.
\end{equation}
Puisque $t_2(z,j)$ et $t_2(z',j')$ sont $7\times7\times2^{n-2}$-groupés, cette relation est aussi vraie avec des couches d'épaisseur $2$, dès lors $t_2(z,j)\vert_{V(0,2,j,z)}=t_2(z,j)$ et $t_2(z',j')\vert_{V(0,2,j',z')}=t_2(z',j')$ et on trouve finalement:
\begin{equation}
t_2(z,j)\vert_{V(0,2,j',z')}=t_2(z',j')\vert_{V(0,2,j,z)}.
\end{equation}
En posant
\begin{equation}
T_2(z)=\left(\bigcup_{1\leq j\leq2(n-2)}t_2(z,j)\right)\sqcup S_2\vert_{\mathbb{R}^2\times[2,25]^{n-2}}
\end{equation}
d'après~\eqref{equationcanalisationsnC} $T_2(z)$ est donc un sous-complexe de $S_2$ tel que $T_2(z)\vert_{U^*_j}$ est égal à $\Theta_2(j)$ si on identifie $U^*_j$ avec $\mathbb{R}^2$. En outre $T_2(z)$ vérifie~\eqref{equationcanalisationsnB} par construction, puisqu'on a enlevé des polyèdres sur une couche d'épaisseur $2$ autour de la frontière de $S_2$.

Définissons par ailleurs $t_1(z,j)$ comme l'ensemble des cubes dyadiques de pas $\frac{1}{p}$ dans la même base que $S_1$ et qui sont inclus dans $\U(\Theta''_1(j))\times V(-1,1,j,z)$:
\begin{equation}
t_1(z,j)=\left\{\delta\colon\delta\subset\U(\Theta''_1(j))\times V(-1,1,j,z)\right\}.
\end{equation}
Soient $z'$ et $j'$ tels que $V(-1,1,j,z)=V(-1,1,j',z')$ et $I_z\neq I_{z'}$ (c'est à dire que $I_z$ et $I_{z'}$ ont pour frontière commune l'un des $U^*_j$ qui est égal à l'un des $U^*_{j'}$). On remarque déjà que $t_1(s,j)$ est $1\times1\times p^{n-2}$-groupé et que 
\begin{equation}
t_1(z,j)\vert_{V(-1,0,j,z)}=t_1(z',j')\vert_{V(-1,0,j',z')}
\end{equation}
puisque l'on a utilisé pour construire ces deux complexes respectivement les restrictions de $S_1$ et $S_2$ à $U^*_j$ et $U^*_{j'}$, qui sont égales.

Construisons maintenant $t'_1(z)$ comme l'ensemble des cubes dyadiques de pas $\frac{1}{p}$ dans la même base que $S_1$, inclus dans $\mathbb{R}^2\times[1,27]^{n-2}$ et à distance comprise entre $r$ et $6$ de $S_2$:
\begin{equation}
t'_1(z)=\left\{\delta\colon r<\min_{x\in\delta}\dist\left(x,\U\left(S_2\vert_{I_z}\right)\right)<6\right\}.
\end{equation}
De même on définit l'ouvert $O'(z)$ dont l'adhérence est composée des cubes situés à distance au plus $6$:
\begin{equation}
\overline{O'(z)}=\U\left(\left\{\delta\colon\min_{x\in\delta}\dist\left(x,\U\left(S_2\vert_{I_z}\right)\right)<6\right\}\right).
\end{equation}
Là encore, par un argument identique à celui utilisé dans la démonstration du lemme~\ref{lemmacanalisations2} on peut montrer qu'il est possible, en ajoutant\slash supprimant quelques cubes à $t'_1(z)$, de faire des suspensions tubulaires\footnote{Ici intervient le fait qu'on ait supposé $p$ impair, pour que $H_z$ coupe les cubes de $t'_1(z)$ en passant par leur centre, et pas le long de leur frontière.} avec les faces en vis-à-vis de $t'_1(z)\vert_{H_{z'}}$ et $S_2\vert_{H_{z'}}$ pour tout $\displaystyle z'\in z+\left\{1,\ldots,26\right\}^{n-2}$. En outre, il est possible de le faire de façon à ce que $t'_1(z)$ soit $1\times1\times p^{n-2}$-groupé (en ajoutant\slash supprimant des paquets de cubes $1\times1\times p^{n-2}$-groupés). On notera $\G(z)$ l'ensemble des graphes utilisés pour réaliser ces suspensions tubulaires. Pour finir on pose:
\begin{equation}
T_1(z)=\left(\bigsqcup_{1\leq j\leq2(n-2)}t_1(z,j)\right)\cup t'_1(z).
\end{equation}

Il ne nous reste plus qu'à montrer que $T_1(z)$ et $T_2(z)$ vérifient les propriétés annoncées dans le lemme à l'intérieur de $I_z$: soit $z'\in[0,27]^{n-2}$, et notons $\nu=\left(\frac{1}{2},\ldots,\frac{1}{2}\right)\in\mathbb{R}^{n-2}$. Deux cas sont possibles:
\begin{itemize}
\item si $z'+\nu\in[2,26]^{n-2}$ alors $T_2(z)\vert_{H_{z+z'}}$ est égal à $S_2(z)\vert_{H_{z'}}$, et $T_1(z)\vert_{H_{z+z'}}$ est égal à $T'_1(z)$ donc on peut utiliser la famille des graphes de $\G(z)$ mentionnée plus haut pour faire une canalisation de $O'(z)\cap H_{z+z'}\setminus\U(T_1(z)\cup T_2(z))$;
\item si $z+z'+\nu\notin[2,26]^{n-2}$ alors notons $z_{\min}$ et $z_{\max}$ respectivement les coordonnées minimale et maximale de $z'$. Là encore, considérons deux cas possibles:
\begin{itemize}
\item si $z_{\min}>1$ et $z_{\max}<26$ alors soit $j$ tel que $z+z'+\nu\in V(1,2,j,z)$: on a démontré que $\G_2(j)$ et un sous-complexe de $\Sigma(j)\cup\F_\partial(\Theta_2(j))$ forment une canalisation de $\open{\U(S^*_2(j))}\setminus\U(\Theta_2(j)\cup\Sigma(j))$. Dès lors en considérant que les graphes de $\G_2(j)$ sont des graphes du $2$-plan $H_{s+z'}$, et que le complexe unidimensionnel $\Sigma(j)$ est un complexe du $2$-plan $H_{z+z'}$ alors il vient que $\G_2(j)$ et un sous-complexe de $\Sigma(j)\cup\F_\partial(T_2(z)\vert_{H_{z+z'}})$ forment une canalisation de $\open{\U(S_2)}\cap H_{z+z'}\setminus\U(T_2(z)\cup\Sigma(j))$.

De plus, par le même argument que dans le cas précédant, $\G(z)$ et un sous-complexe de $\F_\partial((T_1(z)\cup S_2)\vert_{H_{z+z'}})$ forment une canalisation de $O'(z)\cap H_{z+z'}\setminus\U(T_1(z)\cup S_2)$. Rappelons que
\begin{equation}
\F_\partial\left(S_2\vert_{H_{z+z'}}\right)\subset\F_\partial\left(T_2(z)\vert_{H_{z+z'}}\right)
\end{equation}
donc $\G(z)\cup\G_2(j)$ et un sous-complexe de $\Sigma(j)\cup\F_\partial\left((T_1(z)\cup T_2(z))\vert_{H_{z+z'}}\right)$ forment une canalisation de $O'(z)\cap H_{z+z'}\setminus\U(T_1(z)\cup T_2(z))$. On notera respectivement $\G'(z)$ et $\Sigma(z)$ l'ensemble de tous les graphes et l'ensemble de tous les complexes unidimensionnels $\Sigma(j)$ utilisés dans ce cas pour faire les canalisations;
\item si $z_{\min}=0$ ou si $z_{\max}=26$ alors soit $j$ tel que $z+z'+\nu\in V(0,1,j,z)$: on peut refaire une démonstration analogue (cette fois en utilisant les graphes $\G'_1(j)$ et $\G'_3(j)$) pour exhiber des canalisations de $O(z)\cap H_{z+z'}\setminus\U(T_1(z)\cup T_2(z))$. Cette fois-ci on notera respectivement $\G''(z)$ et $\Sigma'(z)$ l'ensemble de tous les graphes et l'ensemble de tous les complexes unidimensionnels $\Sigma(j)$ utilisés pour faire ces canalisations.
\end{itemize}
\end{itemize}

Il est temps de conclure en posant:
\begin{equation}
\begin{aligned}
T_1&=\bigcup_{z\in(28\mathbb{Z})^{n-2}}T_1(z)&T_2&=\bigcup_{z\in(28\mathbb{Z})^{n-2}}T_2(z)\\
\Sigma&=\bigcup_{z\in(28\mathbb{Z})^{n-2}}\Sigma(z)\cup\Sigma'(z)&\overline{O}&=\bigcup_{z\in(28\mathbb{Z})^{n-2}}\overline{O(z)}\cup\overline{O'(z)}\\
\G&=\bigcup_{z\in(28\mathbb{Z})^{n-2}}\G(z)\cup\G'(z)\cup\G''(z)
\end{aligned}
\end{equation}
On vient de démontrer les points 1 et 2 du lemme, l'existence des constantes $\rho_+$ et $\rho_-$ ainsi que le compact $K$ mentionnés dans le point (3) ont quant à elles été obtenues en appliquant le lemme~\ref{lemmacanalisations2}.
\end{proof}

À présent, on dispose de tous les lemmes nécessaires pour démontrer que le théorème~\ref{theoremfusion} de fusion est une propriété inductive sur $n$ pour $n\geq2$.

\begin{lemma}[Fusion en dimension quelconque]\label{lemmafusionn}
Soit $n\geq3$. Si le théorème~\ref{theoremfusion} est vrai en dimension $n-1$, alors il est vrai en dimension $n$.
\end{lemma}

\begin{proof}
On suppose que $n>2$, que le théorème de fusion est vrai en dimension $n-1$ et que $S_1$ et $S_2$ vérifient les hypothèses du théorème de fusion en dimension $n$.

L'isométrie affine $\phi$ de changement de base entre les deux complexes peut être décomposée en $N=\frac{(n+1)(n+2)}{2}$ rotations affines successives dans des plans engendrés par des couples de vecteurs de la base canonique de $S_2$ par exemple. Il est donc possible, en ajoutant des couches de polyèdres successives autour de $S_2$ et en supposant $\rho$ suffisamment grand, de se ramener au cas où $\phi$ est une rotation d'angle $\theta$ dans le plan engendré par deux vecteurs d'une base de $S_2$. Il suffira de réaliser la fusion en faisant $N$ transitions pour démontrer le théorème. Et comme dans le cas de la dimension $2$ il est même possible, en insérant encore les étapes intermédiaires nécessaires, de supposer que
\begin{equation}
\theta\in\left[\theta_{\min},\theta_{\max}\right]
\end{equation}
avec $\theta_{\min}<\theta_{\max}$ deux constantes arbitraires prises dans $\left]0\frac{\pi}{4}\right[$. Bien évidemment, le nombre de transitions à effectuer ne dépendra que de $n$ et des constantes $\theta_{\min}$ et $\theta_{\max}$ choisies. Pour simplifier l'écriture de la démonstration, en notant $(u_1,\ldots,u_n)$ une base canonique de $S_2$ on supposera aussi, quitte à permuter ses vecteurs, que $\phi$ est une rotation dans le $2$-plan $\vect(u_1,u_2)$.

À présent supposons qu'on ait simplement subdivisé $S_2$ vingt-huit fois grâce au lemme~\ref{lemmasubdivision} de subdivision, et appliquons le lemme~\ref{lemmacanalisationsn} à $S_1$ et à $S_2$: on obtient les trois complexes $T_1$, $T_2$ et $\Sigma$ annoncés dans le lemme. D'après~\eqref{equationcanalisationsnB} on sait que seuls les cubes qui étaient à distance au plus $2\sqrt{n}$ de la frontière de $S_2$ n'apparaissent pas dans $T_2$, donc si l'on suppose que $\rho$ est suffisamment grand pour qu'on ait rajouté une couche d'épaisseur $2$ de cubes autour de $S_2$ après subdivision, alors on peut supposer que $T_2\subset S_2$ et que les complexes $T_1$, $T_2$ et $\Sigma$ vérifient toujours les propriétés annoncées. Pour simplifier encore, on notera toujours $S_2$ le complexe obtenu après y avoir découpé les canalisations du lemme, au lieu de $T_2$. Maintenant considérons le complexe $T_1$ : il est de pas $\frac{1}{p}$ donc il suffit de subdiviser $p$ fois $S_1$ pour le raccorder à $T_1$, la condition~\eqref{equationcanalisationsnA} nous assurant qu'on dispose de l'espace nécessaire pour le faire. Par commodité là encore on notera toujours $S_1$ le complexe obtenu après raccordement.

Il est temps d'utiliser notre hypothèse de récurrence. Considérons la famille des hyperplans affines (pour $1\leq k\leq2(n-2)$ et $z\in\mathbb{Z}^{n-2}$):
\begin{equation}
H_{k,z}=z+\mathbb{R}^2\times\mathbb{R}^{k}\times\{0\}\times\mathbb{R}^{n-3-k}
\end{equation}
et les restrictions $\F_{n-1}(S_1)\vert_{H_{k,z}}$ et $\F_{n-1}(S_2)\vert_{H_{k,z}}$: lorsqu'elles sont non vides, ce sont deux complexes $n-1$-dimensionnels qui vérifient les hypothèses du théorème de fusion à ceci près qu'ils ne sont pas forcément à distance suffisante l'un de l'autre. Cependant on peut tout à fait supposer que l'on avait subdivisé préalablement $S_1$ et $S_2$ suffisamment pour que ce soit le cas (d'un nombre de fois $q$ qui ne dépend pas des complexes considérés) et qu'on a jusqu'ici travaillé sur des $q^n$-groupements, de façon à ce que $\frac{q}{2p}>\rho'$ où $\rho'$ est cette distance minimale imposée par le théorème en dimension $n-1$ et $p$ la constante donnée par le lemme~\ref{lemmacanalisationsn}. En appliquant le théorème de fusion en dimension $n-1$ on peut donc remplir toutes les composantes connexes de $\overline{O\cap H_{k,z}\setminus\U(S_2)}$ de complexes $n-1$-dimensionnels, dont on notera $\Theta$ l'union.

Considérons maintenant le complexe unidimensionnel $\Sigma$, et notons $\Sigma'$ l'ensemble des produits cartésiens de ses segments par $\left[-\frac{1}{2},\frac{1}{2}\right]^{n-2}$: par une démonstration analogue à celle du lemme~\ref{lemmasuspensiontubulaire} on peut montrer que parmi les positions à $\epsilon$ près des sommets de $\Sigma$ qui peuvent être déplacés, il est possible d'en trouver telles que les régularités extrêmes des découpages des faces de $\Theta$ par les polyèdres de $\Sigma'$ dans les couches $H_{k,z}$ peuvent être bornées par des constantes multiplicatives ne dépendant pas de $S_1$ et $S_2$ par rapport aux régularités de $S_1$ et $S_2$.

À présent il ne reste plus qu'à utiliser les graphes de la famille $\G$ fournie par le lemme~\ref{lemmacanalisationsn}, puis le lemme~\ref{lemmaproduitcartesien} pour conclure: on peut faire des suspensions tubulaires dans toutes les composantes connexes de $O\setminus\U(\Theta\cup S_1\cup S_2\cup\Sigma')$ et de ce fait, remplir toute l'adhérence de $O$ de polyèdres, de façon à obtenir un complexe $n$-dimensionnel vérifiant toutes les conditions voulues, ce qui achève la démonstration du lemme~\ref{lemmafusionn} et par là celle du théorème~\ref{theoremfusion} par récurrence.
\end{proof}

\backmatter

\section{Preuve informatique du lemme du laboureur}\label{sectionD}

On va à présent donner les algorithmes d'énumération annoncés dans le début de la démonstration du lemme~\ref{lemmalaboureur}. Il y a en tout $2^9=512$ cas à traiter, qu'on pourrait ramener à $2^6=64$ à des changements de base près. Ce nombre encore trop élevé justifie le recours à un outil informatique au lieu d'un travail à la main fastidieux et guère intéressant.

Les algorithmes utilisés étant très simples, on se dispensera de donner une preuve détaillée de leur fonctionnement correct. 

\subsection{Méthode algorithmique}

Dans les algorithmes à venir, pour tout complexe $S$ dyadique unitaire, on va représenter $S\vert_{[0,12]^2}$ par une matrice d'entiers $M=(m_{i,j})_{(i,j)\in\{1,\ldots,12\}^2}$ de taille $12\times12$, vérifiant l'équivalence suivante:
\begin{equation}
\begin{cases}
m_{i,j}=1\Leftrightarrow\Delta(i-1,j-1,1)\in S\vert_{[0,12]^2}\\
m_{i,j}=0\Leftrightarrow\Delta(i-1,j-1,1)\notin S\vert_{[0,12]^2}
\end{cases}
\end{equation}

\subsubsection{$\mathfrak{D}^2$ et le $4\times4$-groupement}

Pour commencer on peut calculer très simplement le nombre d'éléments de $\V(\Delta(x-1,y-1,1))\cap S$ (pour $(x,y)\in[2,11]^2$) en comptant le nombre de cases voisines de $(x,y)$ égales à $1$ dans la représentation de $S$. On obtient alors l'algorithme~\ref{algorithmA}.

\begin{algorithm}
\caption{Calcul de $\#(\V(\Delta(x-1,y-1,1))\cap S)$}
\label{algorithmA}
\begin{algorithmic}[1]
\REQUIRE $M\in\M_{12}(\{0,1\})$ et $(x,y)\in\{2,\ldots,11\}^2$
\ENSURE $s=\#(\V(\Delta(x-1,y-1,1))\cap S)$
\STATE $s\Leftarrow0$
\FORALL{$(u,v)\in\{-1,0,1\}^2$}
\IF{$m_{x+u,y+v}=1$}
\STATE $s\Leftarrow s+1$
\ENDIF
\ENDFOR
\end{algorithmic}
\end{algorithm}

La représentation de $\mathfrak{D}(S)$ est alors obtenue en ne gardant que les cases pour lesquelles ce nombre est égal à $9$, c'est l'algorithme~\ref{algorithmB}.

\begin{algorithm}
\caption{Calcul de la représentation de $\mathfrak{D}(S)$}
\label{algorithmB}
\begin{algorithmic}[1]
\REQUIRE $M\in\M_{12}(\{0,1\})$
\ENSURE $M'\in\M_{12}(\{0,1\})$ telle que $M'$ représente $\mathfrak{D}(S)\vert_{[1,11]^2}$
\FORALL{$(x,y)\in\{1,\ldots,12\}^2$}
\STATE $m'_{x,y}\Leftarrow0$
\ENDFOR
\FORALL{$(x,y)\in\{2,\ldots,11\}^2$}
\IF{$\#(\V(\Delta(x-1,y-1,1))\cap S)=9$}
\STATE $m'_{x,y}\Leftarrow0$
\ENDIF
\ENDFOR
\end{algorithmic}
\end{algorithm}

À présent en posant $S'=\mathfrak(D)^2(S)$ l'algorithme~\ref{algorithmC} permet de vérifier si $S'\vert_{[2,10]^2}$ est $4\times4$-groupé (relativement à l'origine $(2,2)$).

\begin{algorithm}
\caption{Détermine si $S'\vert_{[2,10]^2}$ est $4\times4$-groupé par rapport à l'origine $(2,2)$}
\label{algorithmC}
\begin{algorithmic}[1]
\REQUIRE $M'\in\M_{12}(\{0,1\})$
\ENSURE \TRUE\ si $S'\vert_{[2,10]^2}$ est $4\times4$-groupé, \FALSE\ sinon
\FORALL{$(u,v)\in\{2,6\}^2$}
\FORALL{$(x,y)\in\{0,\ldots,3\}^2$}
\IF{$m'_{x+u,y+v}\neq m'_{u,v}$}
\RETURN\FALSE
\ENDIF
\ENDFOR
\ENDFOR
\RETURN\TRUE
\end{algorithmic}
\end{algorithm}

Considérons l'ensemble des restrictions possibles au carré $[0,12]^2$ de complexes unitaires $4\times4$-groupés par rapport à l'origine $(0,0)$. Il est clair qu'il y en a $2^9$ et en les notant $G_i$ pour $0\leq i<2^9$ on peut les paramétrer par:
\begin{multline}
\forall (x,y)\in\{0,\ldots,3\}^2,\forall (u,v)\in\{0,1,2\}^2\colon\\
\Delta(3u+x,3v+x)\in G_i\Leftrightarrow2^{3u+v}\leq i\bmod{2^{3u+v+1}}.
\end{multline}
L'algorithme~\ref{algorithmD} est utilisé pour générer la représentation de $G_i$.

\begin{algorithm}
\caption{Calcul de la représentation de $G_i$}
\label{algorithmD}
\begin{algorithmic}[1]
\REQUIRE $i\in\{0,\ldots,2^9-1\}$
\ENSURE $M=G_i$
\FORALL{$(u,v)\in\{0,1,2\}^2$}
\FORALL{$(x,y)\in\{0,\ldots,3\}^2$}
\IF{$2^{3u+v}\leq i\bmod{2^{3u+v+1}}$}
\STATE $m_{3u+x,3v+y}\Leftarrow1$
\ELSE
\STATE $m_{3u+x,3v+y}\Leftarrow0$
\ENDIF
\ENDFOR
\ENDFOR
\end{algorithmic}
\end{algorithm}

En utilisant les quatre algorithmes précédents (\ref{algorithmA},~\ref{algorithmB},~\ref{algorithmC} et~\ref{algorithmD}) on peut donner l'algorithme~\ref{algorithmE} qui permet de vérifier si $\mathfrak{D}^2$ préserve la propriété de $4\times4$-groupement (par rapport à l'origine $(2,2)$) pour tout complexe $4\times4$-groupé (par rapport à l'origine $(0,0)$).

\begin{algorithm}
\caption{Vérifier si $\mathfrak{D}^2$ préserve le $4\times4$-groupement}
\label{algorithmE}
\begin{algorithmic}[1]
\ENSURE \TRUE\ si $\mathfrak{D}^2$ préserve la propriété de $4\times4$-groupement, \FALSE\ sinon
\FORALL{$i\in\{0,\ldots,2^9\}$}
\STATE $M\Leftarrow G_i$
\STATE $M\Leftarrow\mathfrak{D}(M)$
\STATE $M\Leftarrow\mathfrak{D}(M)$
\IF{$M$ est $4\times4$-groupé par rapport à l'origine $(2,2)$}
\RETURN\TRUE
\ELSE
\RETURN\FALSE
\ENDIF
\ENDFOR
\end{algorithmic}
\end{algorithm}

\subsubsection{Découpage en complexes linéaires}

Sur le modèle de l'algorithme~\ref{algorithmA}, l'algorithme~\ref{algorithmF} permet de calculer le nombre d'éléments de $\V^*(\Delta(x-1,y-1,1))\cap S$.

\begin{algorithm}
\caption{Calcul de $\#(\V^*(\Delta(x-1,y-1,1))\cap S)$}
\label{algorithmF}
\begin{algorithmic}[1]
\REQUIRE $M\in\M_{12}(\{0,1\})$ et $(x,y)\in\{4,\ldots,9\}^2$
\ENSURE $s=\#(\V^*(\Delta(x-1,y-1,1))\cap S)$
\STATE $s\Leftarrow0$
\FORALL{$(u,v)\in\{-1,0,1\}^2$ tel que $\vert u+v\vert=1$}
\IF{$m_{x+u,y+v}=1$}
\STATE $s\Leftarrow s+1$
\ENDIF
\ENDFOR
\end{algorithmic}
\end{algorithm}

Il existe en tout huit isométries laissant le carré $[3,9]^2$ invariant (l'identité et la symétrie par rapport au centre, les deux symétries par rapport aux diagonales, les deux symétries par rapport aux médianes et les deux rotations d'angles $\frac{\pi}{2}$ et $-\frac{\pi}{2}$). On les notera $\phi_i$ pour $1\leq i\leq8$. L'algorithme~\ref{algorithmG} permet de vérifier si deux complexes restreints à $[3,9]^2$ sont égaux à une isométrie près.

\begin{algorithm}
\caption{Déterminer si deux complexes sont égaux à une isométrie près}
\label{algorithmG}
\begin{algorithmic}[1]
\REQUIRE $M,M'\in\M_{12}(\{0,1\})$
\ENSURE \TRUE\ si les deux grilles restreintes à $[3,9]^2$ sont égales à une isométrie près, \FALSE\ sinon
\FORALL{$i\in\{1,\ldots,8\}$}
\STATE $s\Leftarrow\TRUE$
\FORALL{$(x,y)\in\{3,9\}^2$}
\IF{$m_{x,y}\neq m_{\phi_i(x,y)}$}
\STATE $s\Leftarrow\FALSE$
\ENDIF
\ENDFOR
\IF{$s$}
\RETURN\TRUE
\ENDIF
\ENDFOR
\RETURN\FALSE
\end{algorithmic}
\end{algorithm}

Pour finir, l'algorithme~\ref{algorithmH} vérifie que pour tout complexe $S$ $4\times4$-groupé, $(S\setminus\mathfrak{D}(S))\vert_{[3,9]^2}$ ne contient que des cubes à deux voisins tangents, et dans le même temps se charge de dresser la liste des configurations possibles de $(\mathfrak{D}(S)\setminus\mathfrak{D}^2(S))\vert_{[3,9]^2}$ lorsque $S$ parcourt l'ensemble des complexes dyadiques unitaires $4\times4$-groupés par rapport à l'origine $(0,0)$.

\begin{algorithm}
\caption{Vérification que $(S\setminus\mathfrak{D}(S))\vert_{[3,9]^2}$ est sulciforme en cycles et liste des configurations possibles de $(\mathfrak{D}(S)\setminus\mathfrak{D}^2(S))\vert_{[3,9]^2}$ pour tout complexe $S$ $4\times4$-groupé}
\label{algorithmH}
\begin{algorithmic}[1]
\ENSURE $s=\TRUE$ si pour tout complexe $S$ $4\times4$-groupé:
\begin{equation} \forall\delta\in(S\setminus\mathfrak{D}(S))\vert_{[4,8]^2}\colon\#(\V^*(\delta)\cap\left(S\setminus\mathfrak{D}(S))\vert_{[3,9]^2}\right)=2
\end{equation}
\FALSE\ sinon. $L$ contient la liste des configurations possibles à une isométrie près de $(\mathfrak{D}(S)\setminus\mathfrak{D}^2(S))\vert_{[3,9]^2}$
\STATE $L\Leftarrow\emptyset$
\STATE $s\Leftarrow\TRUE$
\FORALL{$i\in\{0,\ldots,2^9\}$}
\STATE $M\Leftarrow G_i$
\STATE $M'\Leftarrow M\setminus\mathfrak{D}(M)$
\IF{$\exists\delta\in S'\vert_{[4,8]^2}\colon\#(\V^*(\Delta(x-1,y-1,1))\cap S')\neq2$}
\STATE $s\Leftarrow\FALSE$
\ENDIF
\STATE $M''\Leftarrow\mathfrak{D}(M)\setminus\mathfrak{D}^2(M)$
\IF{$M''\notin L$ à une restriction à $[3,9]^2$ et une isométrie près}
\STATE $L\Leftarrow L\cup\{M''\}$
\ENDIF
\ENDFOR
\end{algorithmic}
\end{algorithm}

\subsection{Implémentation en C}

Ce programme est la traduction en langage C des algorithmes donnés plus haut.

\begin{lstlisting}
#include <stdio.h>
#include <stdlib.h>
#include <assert.h>

/* Taille des groupements considérés (ici ¤$4\times 4$¤) */
#define N 4

/* Nombre de groupements dans une grille (ici ¤$3\times 3$¤) */
#define M 3

/* Une grille de taille ¤$12\times 12$¤ */
typedef int GRID[M*N][M*N];

/* Un masque pour remplir les grilles ¤$4\times 4$¤-groupées */
typedef int MASK[M*M];

/* Le tableau des ¤$2^9$¤ masques possibles */
MASK masks[1<<(M*M)];

/* La liste des grilles trouvées, et le nombre d'éléments qu'elle contient */
#define TAB_GRID_MAX 100
GRID tab_grid[TAB_GRID_MAX];
int tab_grid_count=0;

/* Affiche une grille restreinte à ¤$[2,10]^2$¤ */
void print_big_GRID(GRID g){
  int x,y;
  for (x=N/2;x<N*(M-1)+N/2;x++){
    for (y=N/2;y<N*(M-1)+N/2;y++)
      if (g[x][y])
        printf("# ");
      else
        printf("· ");
    printf("\n");
  }
}

/* Affiche une grille restreinte à ¤$[3,9]^2$¤ */
void print_small_GRID(GRID g){
  int x,y;
  for (x=N-1;x<N*(M-1)+1;x++){
    for (y=N-1;y<N*(M-1)+1;y++)
      if (g[x][y])
        printf("# ");
      else
        printf("· ");
    printf("\n");
  }
}

/* Copie la grille src dans dest */
void copy_GRID(GRID src,GRID dest){
  int x,y;
  for (x=0;x<M*N;x++)
    for (y=0;y<M*N;y++)
      dest[x][y]=src[x][y];
}

/* Remplit le tableau de masques */
void fill_masks(void){
  int i,j;
  for (i=0;i<1<<(M*M);i++)
    for (j=0;j<M*M;j++)
      masks[i][j]=(i & (1<<j))!=0;
}

/* Remplit une grille selon un masque */
void fill_GRID(MASK m,GRID g){
  int X[9]={0,N,2*N,0,N,2*N,0,N,2*N};
  int Y[9]={0,0,0,N,N,N,2*N,2*N,2*N};
  int i,x,y,u,v;
  for (i=0;i<M*M;i++){
    x=X[i];
    y=Y[i];
    for (u=0;u<N;u++)
      for (v=0;v<N;v++)
        g[x+u][y+v]=m[i];
  }
}

/* Calcul de l'image de ¤$(x,y)\in[0,12]^2$¤ par l'homothétie ¤$\phi_{index}$¤ */
void isometric_transform(int *x, int *y,int index){
  int z;
  switch(index){
    case (1):
      *x=N*M-1-*x;
      break;
    case (2):
      *y=N*M-1-*y;
      break;
    case (3):
      *x=N*M-1-*x;
      *y=N*M-1-*y;
      break;
    case (4):
      z=*x;
      *x=*y;
      *y=z;
      break;
    case (5):
      z=*x;
      *x=N*M-1-*y;
      *y=z;
      break;
    case (6):
      z=*x;
      *x=*y;
      *y=N*M-1-z;
      break;
    case (7):
      z=*x;
      *x=N*M-1-*y;
      *y=N*M-1-z;
  }
}

/* Copie l'image de la grille src par l'homothétie ¤$\phi_{index}$¤ dans dest */
void copy_isometric_transform_GRID(GRID src,GRID dest,int index){
  int x,y,u,v;
  for (x=0;x<M*N;x++)
    for (y=0;y<M*N;y++){
      u=x;
      v=y;
      isometric_transform(&u,&v,index);
      dest[u][v]=src[x][y];
    }
}

/* Détermine si 2 grilles restreintes à ¤$[2,9]^2$¤ sont égales à une isométrie près */
int compare_isometric_transformed_GRID(GRID g1,GRID g2){
  int s,i,x,y;
  GRID g;
  for (i=0;i<8;i++){
    s=1;
    copy_isometric_transform_GRID(g1,g,i);
    for (x=N-1;x<N*(M-1)+1;x++)
      for (y=N-1;y<N*(M-1)+1;y++)
        s=s && (g[x][y]==g2[x][y]);
    if (s)
      return 1;
  }
  return 0;
}

/* Ajoute une grille à la liste si elle n'y est pas déjà à une isométrie près */
void insert_GRID(GRID g){
  int i;
  for (i=0;i<tab_grid_count;i++)
    if ( compare_isometric_transformed_GRID(g,tab_grid[i]))
      return;
  assert(tab_grid_count<TAB_GRID_MAX);
  copy_GRID(g,tab_grid[tab_grid_count]);
  tab_grid_count++;
}

/* Calcule ¤$\#(\mathcal{V}(\Delta(x,y,1))\cap T)$¤ */
int neighbor_count(int x,int y,GRID g){
  int r=0,u,v;
  for (u=-1;u<=1;u++)
    for (v=-1;v<=1;v++)
      if (g[x+u][y+v])
        r++;
  return r;
}

/* Calcule ¤$\#(\mathcal{V}^*(\Delta(x,y,1))\cap T)$¤ */
int star_neighbor_count(int x,int y,GRID g){
  int r=0,u,v;
  for (u=-1;u<=1;u++)
    for (v=-1;v<=1;v++)
      if (((u==0)^(v==0)) && (g[x+u][y+v]))
        r++;
  return r;
}

/* Calcule ¤$\mathfrak{D}(T)$¤ */
void dig_GRID(GRID g){
  int x,y;
  GRID h;
  copy_GRID(g,h);
  for (x=1;x<M*N-1;x++)
    for (y=1;y<M*N-1;y++)
      if (neighbor_count(x,y,g)!=9)
        h[x][y]=0;
  copy_GRID(h,g);
}

/* Vérifie si une grille est bien ¤$4\times 4$¤-groupée par rapport à l'origine ¤$(2,2)$¤ */
void check_GRID_is_4_4_group(GRID g){
  int u,v,x,y;
  for (x=N;x<N*(M-1);x++)
    for (y=N;y<N*(M-1);y++)
      for (u=-N/2;u<=N/2;u+=N)
        for (v=-N/2;v<=N/2;v+=N)
          if (g[x+u][y+v] ^ g[N+u][N+v]){
            printf("Ce complexe n'est pas 4x4-groupé:\n");
            print_big_GRID(g);
            abort();
          }
}

/* Vérifie si une grille ne contient que des points d'ordre 2 */
void check_GRID_is_cyclic(GRID g){
  int k,x,y;
  for (x=N;x<N*(M-1);x++)
    for (y=N;y<N*(M-1);y++){
      k=star_neighbor_count(x,y,g);
      if ((k<1) && g[x][y]){
        printf("Ce complexe n'est pas une restriction de complexe sulciforme cyclique:\n");
        print_small_GRID(g);
        abort();
      }
    }
}

/* Différence ensembliste entre 2 grilles */
void substract_GRID(GRID g1,GRID g2){
  int x,y;
  for (x=0;x<N*M;x++)
    for (y=0;y<N*M;y++)
      g2[x][y]=g2[x][y] && !g1[x][y];
}

/* Fonction principale */
void main(){
  GRID g,h;
  int i;
  fill_masks();
  for (i=0;i<1<<(M*M);i++){
    fill_GRID(masks[i],g);

    copy_GRID(g,h);
    dig_GRID(g);
    substract_GRID(g,h);
    check_GRID_is_cyclic(h);

    copy_GRID(g,h);
    dig_GRID(g);
    check_GRID_is_4_4_group(g);

    substract_GRID(g,h);
    insert_GRID(h);
  }
  printf("Les deux propriétés\ntestées sont vérifiées.\n");
  printf("Il y a %i complexes\ndans la liste.\n\n\n\n\n",tab_grid_count);
  for (i=0;i<tab_grid_count;i++){
    printf("\nComplexe no %i:\n",i+1);
    print_small_GRID(tab_grid[i]);
  }
}
\end{lstlisting}

\subsection{Résultats donnés par le programme}

Notre programme se compile sans message d'erreur ni d'avertissement avec une version récente de \verb|gcc|. Il s'exécute en moins d'une seconde et imprime les résultats regroupés dans la table~\ref{tableoutput}.

\begin{table}
\setlength\columnsep{90pt}
\begin{multicols}{3}
\small
\begin{verbatim}
Les deux propriétés
testées sont vérifiées.
Il y a 14 complexes
dans la liste.





Complexe no 1:
· · · · · · 
· · · · · · 
· · · · · · 
· · · · · · 
· · · · · · 
· · · · · · 

Complexe no 2:
· · · · · · 
· · · · · · 
· · # # · · 
· · # # · · 
· · · · · · 
· · · · · · 

Complexe no 3:
· · · · · · 
· · · · · · 
# # # # · · 
# # # # · · 
· · · · · · 
· · · · · · 

Complexe no 4:
· · # # · · 
· · # # · · 
# # # # · · 
# # # # · · 
· · · · · · 
· · · · · · 

Complexe no 5:
· · · # · · 
· · · # · · 
· · · # · · 
# # # # · · 
· · · · · · 
· · · · · · 

Complexe no 6:
· · # # · · 
· · # # · · 
· · # # · · 
· · # # · · 
· · # # · · 
· · # # · · 

Complexe no 7:
· · # # · · 
· · # # · · 
# # # # · · 
# # # # · · 
· · # # · · 
· · # # · · 

Complexe no 8:
· · · # · · 
· · · # · · 
· · · # · · 
# # # # · · 
· · # # · · 
· · # # · · 

Complexe no 9:
· · · # · · 
· · · # · · 
· · · # · · 
· · · # · · 
· · · # · · 
· · · # · · 

Complexe no 10:
· · # # · · 
· · # # · · 
# # # # # # 
# # # # # # 
· · # # · · 
· · # # · · 

Complexe no 11:
· · · # · · 
· · · # · · 
· · · # # # 
# # # # # # 
· · # # · · 
· · # # · · 

Complexe no 12:
· · · # · · 
· · · # · · 
· · · # # # 
· · · # # # 
· · · # · · 
· · · # · · 

Complexe no 13:
· · # · · · 
· · # · · · 
# # # · · · 
· · · # # # 
· · · # · · 
· · · # · · 

Complexe no 14:
· · · · · · 
· · · · · · 
· · · · · · 
· · · # # # 
· · · # · · 
· · · # · · 
\end{verbatim}
\end{multicols}
\caption{Résultats imprimés par le programme\label{tableoutput}}
\end{table}

\subsection{Fin de la démonstration du lemme}

On a regroupé dans le tableau suivant les quatorze différentes configurations possibles à une isométrie près des cubes de $U$ sur un voisinage de taille $6\times6$. Les traits en pointillés délimitent un voisinage de taille $4\times4$ qui constitue les bords du motif de pavage permettant de recouvrir $U$ (puisqu'on a vu que $\mathfrak{D}^2$ préserve le $4\times4$-groupement). C'est à dire que quel que soit le complexe original $S$ $4\times4$-groupé, toute restriction de $U$ à $z+[-1,5]^2$ pour $z\in(4\mathbb{Z})^2$ est égale --- à une isométrie près --- à l'un des motifs du tableau, et toute restriction de $U$ à $z+[0,4]^2$ est égale à l'un des motifs du tableau délimité par les traits en pointillés. La couche d'épaisseur $1$ autour du carré en pointillés va être utilisée pour calculer le nombre de voisins des cubes dyadiques à l'intérieur des pointillés.

\begin{center}
\fbox{\includegraphics[width=0.8\linewidth]{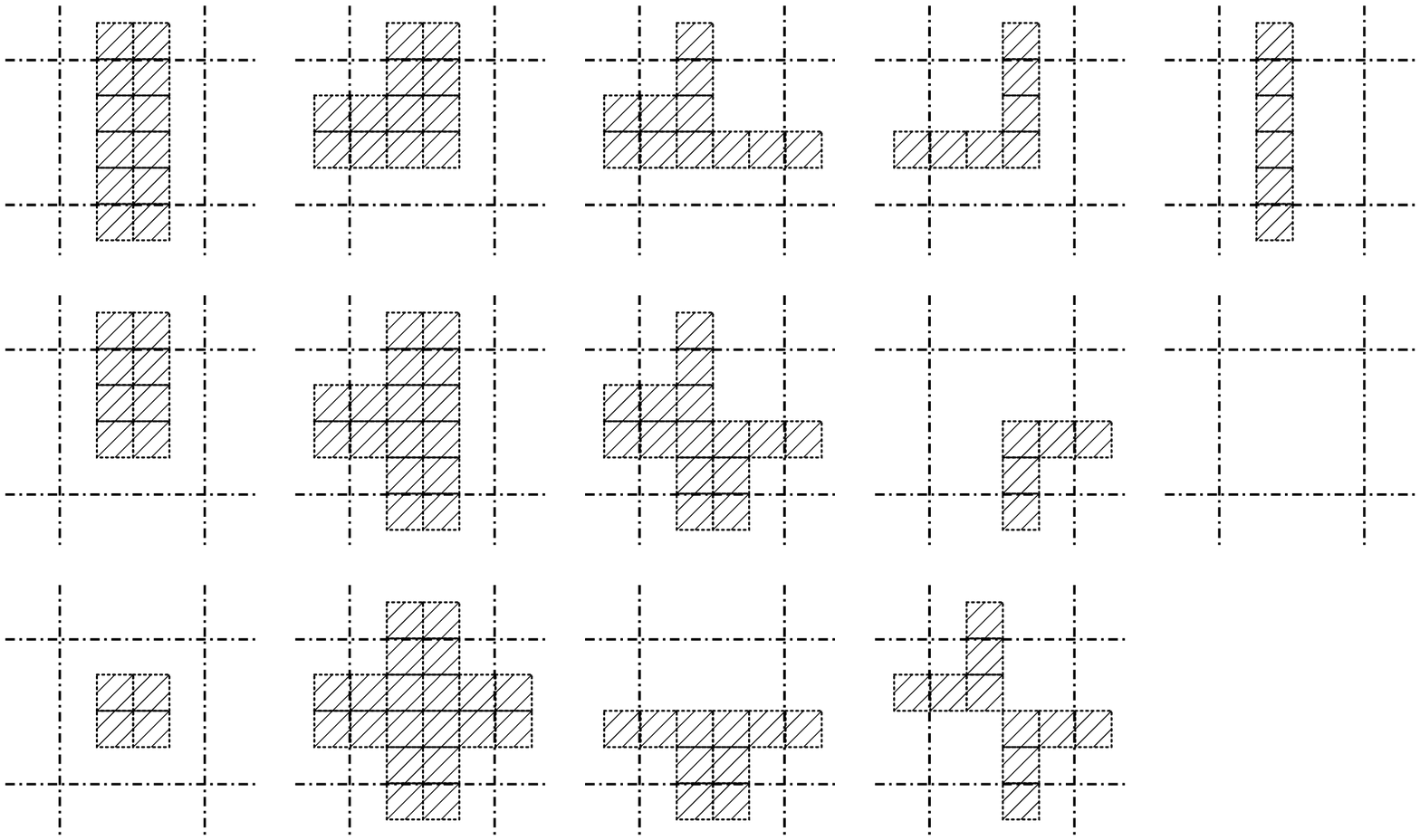}}
\end{center}

Il nous reste encore à extraire de $U$ les deux complexes $U_1$ et $U_2$ annoncés. Considérons le découpage présenté dans le tableau suivant: les hachures verticales représentent $U_1$ les hachures verticales $U_2$. Les traits en gras définissent une partition de $U_2$ en sillons, ou les endroits où des cubes de $U_2$ touchent des cubes de $U_1$.

\begin{center}
\fbox{\includegraphics[width=0.8\linewidth]{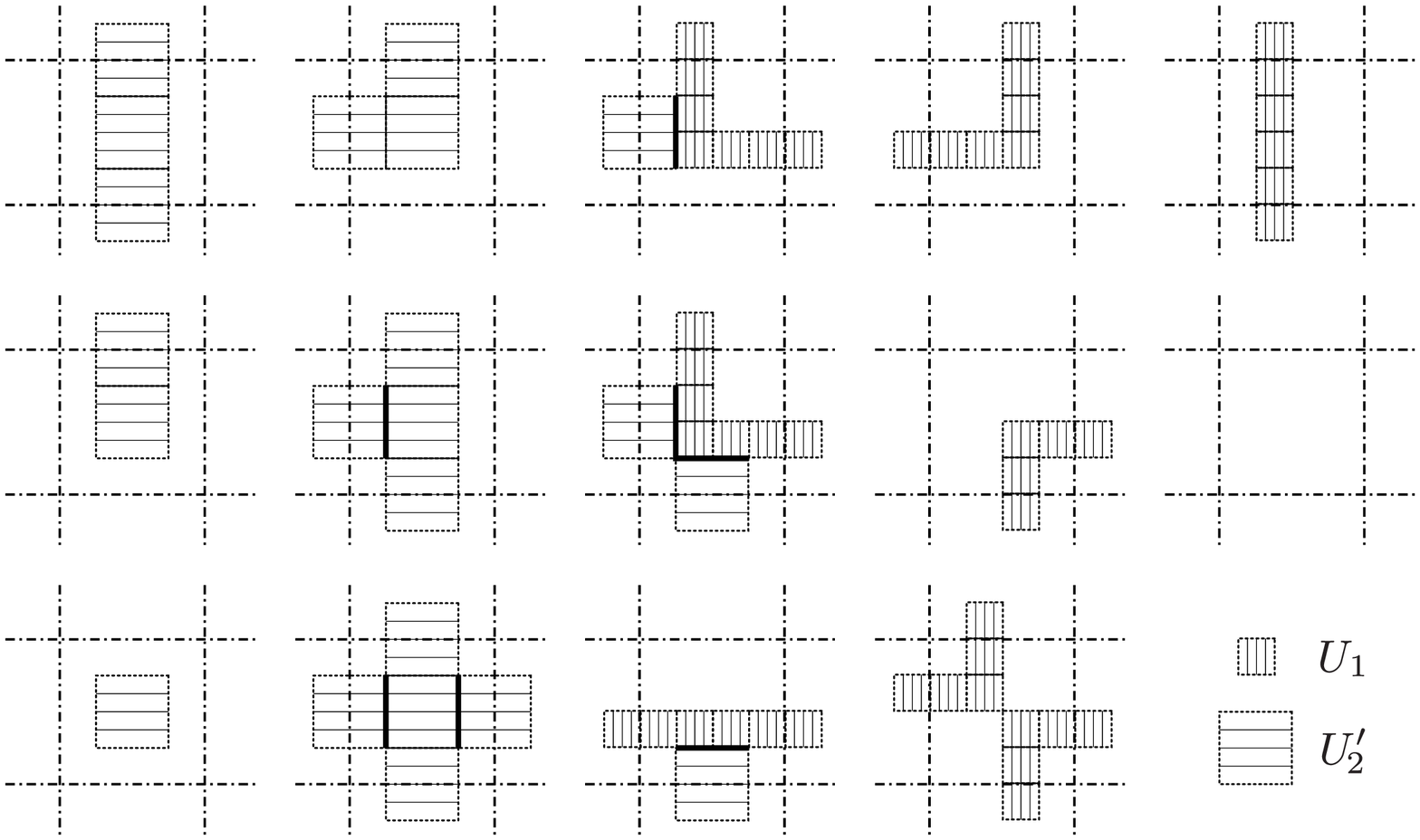}}
\end{center}

Si l'on observe les restrictions de $U_1$ au motif de pavage de taille $4\times4$ matérialisé par les traits en pointillé, on constate que $U_1$ est bien sulciforme en cycles (c'est à dire que tout cube de $U_1$ a deux voisins tangents, éventuellement à l'extérieur du carré). En outre $U_2$ est bien $2\times2$-groupé (par rapport à une origine placée au centre du carré du motif) et si on observe les partitions de son $2\times2$-groupement $U'_2$ délimitées par les traits en gras, on constate bien que chacun des cubes des sous-complexes a au plus deux voisins tangents, donc que $U'_2$ est quasi-sulciforme. De plus, les seuls cubes de $U'_2$ qui touchent des cubes de $U_1$ sont des extrémités de $U'_2$ (c'est à dire qui ont au plus un voisin tangent dans $U'_2$).

Une écriture formelle des explications qu'on vient de donner serait vraisemblablement fastidieuse et guère plus convaincante, on en restera là et on considèrera que les trois complexes $T$, $U_1$ et $U_2$ qu'on a extraits de $S$ vérifient bien toutes les propriétés annoncées dans le lemme~\ref{lemmalaboureur}.

\bibliography{Polyedres}

\begin{thebibliography}{Mat95}

\bibitem[Bas59]{bastiani}
A.~Bastiani.
\newblock {Polyèdres convexes dans les espaces vectoriels topologiques}.
\newblock Sémin. de Topologie et de Géométrie différentielle Ch. Ehresmann 1
  (1957/58), No.19, 46 p., 1959.

\bibitem[Dav03]{david}
G.~David.
\newblock {Limits of Almgren quasiminimal sets}.
\newblock In {\em Harmonic Analysis at Mount Holyoke: Proceedings of an
  Ams-Ims-Siam Joint Summer Research Conference on Harmonic Analysis, June
  25-July 5, 2001, Mount Holyoke College, South Hadley, Ma}, volume~32.
  American Mathematical Society, 2003.

\bibitem[DP07]{depauw:acr}
T.~De~Pauw.
\newblock {Approximating compact rectifiable surfaces in Hausdorff measure and
  in Hausdorff distance by locally acyclic surfaces having the same boundary}.
\newblock 2007.

\bibitem[Kir34]{kirszbraun}
M.D. Kirszbraun.
\newblock {Uber die zusammenziehenden und Lipschitzschen Transformationen}.
\newblock {\em Fund. Math}, 22:77--108, 1934.

\bibitem[KM40]{kreinmilman}
M.~Krein and D.~Milman.
\newblock {On extreme points of regular convex sets}.
\newblock {\em Studia Math}, 9:133--138, 1940.

\bibitem[Mat95]{mattila}
P.~Mattila.
\newblock {\em {Geometry of Sets and Measures in Euclidean Spaces: Fractals and
  Rectifiability}}.
\newblock Cambridge University Press, 1995.

\bibitem[Rei60]{reifenberg}
E.R. Reifenberg.
\newblock {Solution of the Plateau problem for m-dimensional surfaces of
  varying topological type}.
\newblock {\em Acta Mathematica}, 104(1):1--92, 1960.

\end{thebibliography}

\end{document}